\documentclass[final]{elsarticle}

\usepackage{hyperref}
\usepackage{a4wide}
\usepackage{amsmath}
\usepackage{amsthm}
\usepackage{amssymb}
\usepackage{amsfonts}
\usepackage{mathrsfs}
\usepackage{stmaryrd}
\usepackage{booktabs}
\usepackage{graphicx}
\usepackage{paralist}
\usepackage{xcolor}
\usepackage{xspace}
\usepackage{enumitem}
\usepackage{units}
\usepackage{subcaption}
\usepackage{cleveref}
\usepackage{siunitx}
\usepackage{multirow}
\usepackage{bm}

\newcommand{\dumux}{DuMu\textsuperscript{x}\xspace}

\newcommand{\eg}{e.g.\@\xspace}
\newcommand{\ie}{i.e.\@\xspace}

\newcommand{\domain}{\Omega}
\newcommand{\interface}{\gamma}
\newcommand{\surface}{\Gamma}

\newcommand{\fracTips}{\hat{\surface}}

\newcommand{\aperture}{a}
\newcommand{\jump}[1]{\llbracket #1 \rrbracket}
\newcommand{\crossSection}{\epsilon}
\newcommand{\fracLength}{L}
\newcommand{\fracLocalGrad}[1]{\grad_{#1} \, }
\newcommand{\fracLocalDiv}[1]{\nabla_{#1} \, \scal }
\newcommand{\N}{\mathbf{N}}
\newcommand{\fracNet}{\mathcal{F}}
\newcommand{\fracture}{f}

\newcommand{\condInterfaceSet}{\mathcal{H}}
\newcommand{\blockInterfaceSet}{\mathcal{B}}
\newcommand{\sdPartitionSet}{I}
\newcommand{\neighborDomainSet}{\mathcal{N}}
\newcommand{\numSubSets}{n}
\newcommand{\numFractureSides}{m}

\newcommand{\juncIdx}{\mathrm{j}}
\newcommand{\isIdx}{\mathrm{is}}
\newcommand{\fracIdx}{\mathrm{f}}
\newcommand{\pmIdx}{\mathrm{b}}
\newcommand{\diriIdx}{\mathrm{D}}
\newcommand{\neumIdx}{\mathrm{N}}
\newcommand{\indexTupleSet}{\mathcal{T}}
\newcommand{\couplIdx}{\mathrm{c}}
\newcommand{\internalIdx}{\mathrm{int}}
\newcommand{\externalIdx}{\mathrm{ext}}
\newcommand{\codimOneIndexMap}{\mathcal{I}}
\newcommand{\eq}{\mathrm{eq}}

\newcommand{\tpfa}{\textsc{tpfa}\xspace}
\newcommand{\tpfaDfm}{\tpfa-\textsc{dfm}\xspace}
\newcommand{\mpfa}{\textsc{mpfa}\xspace}
\newcommand{\mpfaDfm}{\mpfa-\textsc{dfm}\xspace}

\newcommand{\boxScheme}{\textsc{box}\xspace}
\newcommand{\boxDfm}{\boxScheme-\textsc{dfm}\xspace}
\newcommand{\eboxDfm}{\textsc{ebox}-\textsc{dfm}\xspace}
\newcommand{\eboxMortarDfm}{\textsc{ebox}-\textsc{mortar}-\textsc{dfm}\xspace}

\newcommand{\mesh}{\mathcal{M}}
\newcommand{\element}{E}

\newcommand{\elemFacet}{\varsigma}
\newcommand{\meshFacetSet}{\mathcal{X}}
\newcommand{\cvSet}{\mathcal{K}}
\newcommand{\cv}{K}
\newcommand{\otherCV}{L}
\newcommand{\otherFace}{\Xi}
\newcommand{\scv}{\kappa}
\newcommand{\scvSet}{\mathcal{C}}
\newcommand{\face}{\sigma}
\newcommand{\faceSet}{\mathcal{S}}
\newcommand{\discFlux}{F}
\newcommand{\basisFunc}{\varphi}
\newcommand{\projection}{\Pi}
\newcommand{\errorNorm}{\varepsilon}
\newcommand{\discLength}{\eta}
\newcommand{\mortarGridFactor}{\psi}
\newcommand{\mortarElemSet}{\mathcal{E}}
\newcommand{\someEntity}{\iota}

\newcommand{\head}{h}
\newcommand{\mortar}{\Lambda}

\newcommand{\darcyVel}{\mathbf{q}}
\newcommand{\someVariable}{v}
\newcommand{\mortarHelperFlux}{\Psi}

\newcommand{\perm}{\mathbf{K}}

\newcommand{\source}{q}
\newcommand{\permAngle}{\Theta}

\newcommand{\I}{\mathbf{I}}
\newcommand{\n}{\mathbf{n}}
\newcommand{\R}{\mathbb{R}}              
\newcommand{\eqbydef}{:=}                
\newcommand{\scal}{\cdot}
\newcommand{\meas}[1]{\lvert{#1}\rvert}  
\newcommand{\grad}{{\nabla}}             
\newcommand{\divergence}{{\nabla \scal}} 
\newcommand{\fracTangGrad}[1]{\grad_{#1}}
\newcommand{\fracTangDiv}[1]{\nabla_{#1} \scal}
\newcommand{\position}{\mathbf{x}}
\newcommand{\opPlaceHolder}{\left( * \right)} 
\newcommand{\surfIntegrand}{\mathrm{d}\Gamma}
\newcommand{\volIntegrand}{\mathrm{d}V}
\newcommand{\Cos}[1]{\cos\left( #1 \right)}

\newcommand{\exact}{\star} 
\newcommand{\minDim}{\delta}

\journal{Journal of Computational Physics}

\biboptions{authoryear}

\begin{document}

\begin{frontmatter}

\title{Comparison of cell- and vertex-centered finite-volume schemes for flow in fractured porous media}

\author[addressIWS]{Dennis Gl\"aser}
\ead{dennis.glaeser@iws.uni-stuttgart.de}
\author[addressIWS]{Martin Schneider}
\ead{martin.schneider@iws.uni-stuttgart.de}
\author[addressIWS]{Bernd Flemisch}
\ead{bernd.flemisch@iws.uni-stuttgart.de}
\author[addressIWS]{Rainer Helmig}
\ead{rainer.helmig@iws.uni-stuttgart.de}

\cortext[mycorrespondingauthor]{Corresponding author}
\address[addressIWS]{Institute for Modelling Hydraulic and Environmental Systems,
                     University of Stuttgart,
                     Pfaffenwaldring 61,
                     70569 Stuttgart, Germany}

\begin{abstract}
Flow in fractured porous media is of high relevance in a variety of geotechnical
applications, given the fact that they ubiquitously occur in nature and that they
can have a substantial impact on the hydraulic properties of rock. As a response
to this, an active field of research has developed over the past decades, focussing
on the development of mathematical models and numerical methods to describe and
simulate the flow processes in fractured rock. In this work, we present a comparison
of different finite volume-based numerical schemes for the simulation of flow in
fractured porous media by means of numerical experiments. Two novel vertex-centered
approaches are presented and compared to well-established numerical schemes in
terms of convergence behaviour and their performance on benchmark studies taken
from the literature. The new schemes show to produce results that are in good
agreement with those of established methods while being computationally less
expensive on unstructured simplex grids.
\end{abstract}

\begin{keyword}
porous media \sep fractures \sep discrete fractures
\sep coupling \sep finite volumes \sep box scheme
\end{keyword}

\end{frontmatter}

\section{Introduction}
\label{sec:introduction}

Fractures are a common feature in many geological formations throughout the
planet, and therefore, they are of particular importance in a wide range of geotechnical
applications from groundwater management over energy production and
storage to the disposal of waste \citep{Berkowitz2002}. Although fractures
are typically characterized by rather small apertures, they can have significant
lateral extents and form interconnected networks that have substantial
influence on the overall hydraulic and mechanical behavior of rock
\citep{Jaeger2007fundamentals}. Besides naturally ocurring fractures, they might
also be introduced into a rock mass as a consequence of human activities, which
is the case, for instance, in geothermal energy or unconventional
shale gas production techniques. In the context of radioactive waste disposal,
fractures can occur in the vicinity of emplacement tunnels as a consequence of
their excavation, and thus, they have to be taken into account in safety
assessment efforts \citep{Lisjak2016FDEM}.\\

As a response to the importance of fractures, a number of modelling approaches
and simulation tools have been developed over the past decades. Typically, the
approaches are classified into continuum fracture models and discrete fracture
matrix (dfm) models. In the first class of models, the fracture geometry is not
explicitly captured, but the overall medium is described in terms of effective properties
in so-called single-continuum models \citep{pruess1990thermohydrologic,Royer2002}
or by describing the fracture network with one or more individual overlapping continua
in dual- or multi-continuum models
\citep{Warren1963,kazemi1976numerical,Lawrence199MINC,zimmerman1993numerical}.
In dfm models, the geometries of the most important fractures are resolved, while
small scale fractures might still be accounted for by means of effective properties
or additional continua. This allows for an explicit representation of fractures
whose sizes are comparable to that of the domain and which are difficult to upscale,
while maintaining the computational efficiency of continuum
approaches where applicable. In such models, the fractures are incorporated either by
describing them as heterogeneities in the rock, resolving the interior of
the fractures with the computational mesh \citep{Matthai2007}, or
as lower-dimensional interfaces embedded in the bulk medium. The latter
approach has received significant attention in the past and a variety of models
and numerical schemes have been developed, which can be further subdivided into
conforming and non-conforming models. In the former, the computational grid used
in the bulk medium is constrained such that the grid element faces are aligned
with the fractures, while in the latter, bulk medium and fractures can be
discretized independently. Examples for conforming approaches can be found, for
instance,
in~\cite{Karimi2004TpfaDfm,martin2005modeling,Sandve2012MpfaDfm,Ahmed2015MpfaDfm,Ahmed2017MpfaFpsDfm,brenner2016gradient,Boon2018,nordbotten2019unified},
while non-conforming models are the subject
of~\citet{Schwenck2015,Flemisch2016,Tene2017pEDFM,koppel2018lagrange,schadle20193d}.
While mesh generation is much simpler for non-conforming methods, they might struggle
to capture relevant phenomena in the case of low-permeable fractures that act as
barriers for flow.\\

In this work, we consider conforming approaches, and compare several cell-
and vertex-centered numerical schemes in terms of accuracy and computational
efficiency. Some of the schemes were presented in the references mentioned above,
while two of the vertex-centered schemes, to the authors' knowledge, have not yet been
introduced. Therefore, they receive special attention in this work.
All schemes have been implemented into the same software framework,
namely \dumux~\citep{Dumux,Kochetal2020Dumux}, and all computations have been
carried out on the same machine, which allows for an easier comparison
of the computational cost. We investigate the convergence behavior of the
schemes against analytical and equi-dimensional reference solutions, and apply
them to benchmark cases published in~\citet{Flemisch2018Benchmarks,Berre2021Benchmarks},
which enables us to compare the results to those of a large number of numerical schemes.\\

This paper is structured as follows. In~\cref{sec:equiModel}, we
present the equi-dimensional mathematical model for flow in fractured porous media,
where the fractures are treated as heterogeneities. Subsequently, the mixed-dimensional
approximation of the model is presented in~\cref{sec:mixedModel}. We briefly describe
the considered numerical schemes in~\cref{sec:discretization}, and provide a more
detailed description of the novel vertex-centered schemes.
\Cref{sec:examples} contains the numerical examples, and we conclude the paper
in~\cref{sec:conclusions} with a summary and outlook.

\section{Equi-dimensional model}
\label{sec:equiModel}
\begin{figure}[ht]
    \centering
    \includegraphics[width=0.9\textwidth]{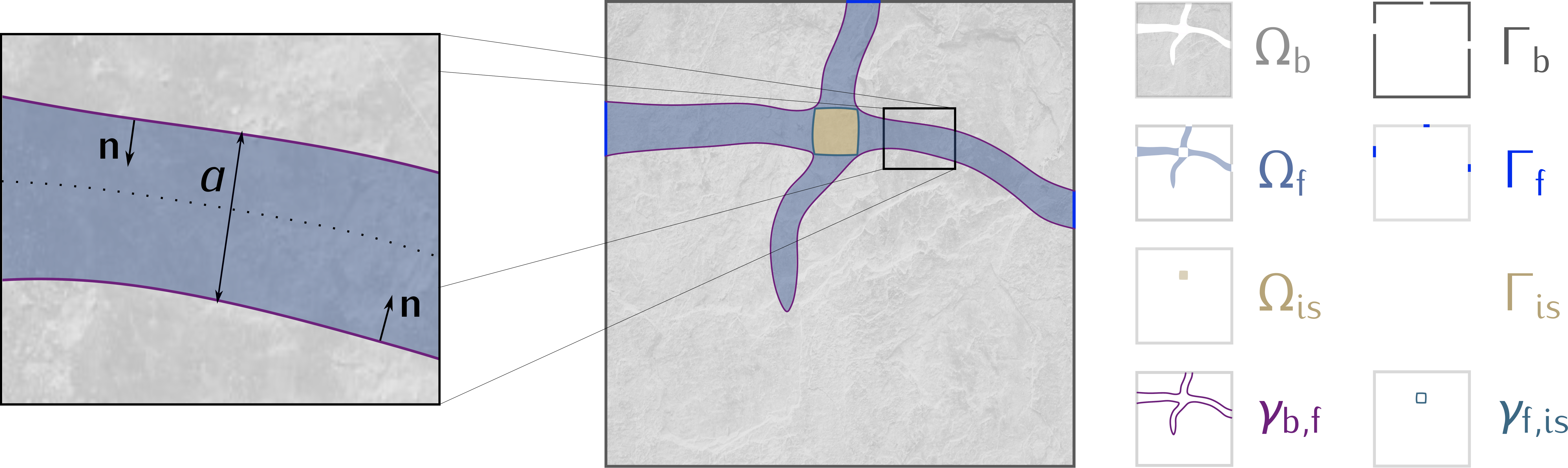}
    \caption{\textbf{Equi-dimensional domain decomposition}. The domain $\domain$ is decomposed into the disjoint partitions $\domain_\fracIdx$ for the fractures, $\domain_\isIdx$ for the fracture intersections and $\domain_\pmIdx$ for the surrounding bulk porous medium. The interfaces between the domains are denoted by $\interface_{\pmIdx, \fracIdx}$ and $\interface_{\fracIdx, \isIdx}$ while the external boundaries of the sub-domains are denoted by $\surface_\pmIdx$, $\surface_\fracIdx$ and $\surface_\isIdx$, respectively. Furthermore, the close-up on the left side illustrates the fracture aperture $\aperture$. In a two-dimensional setting, junctions of fracture intersections do not occur and are not visualized here. Please also note that the thickness of the fractures is exaggerated in this figure for illustrative purposes.}
    \label{fig:modeldomain_equi}
\end{figure}
Let us consider a domain $\domain \subset \R^3$ with boundary $\partial \domain$.
Furthermore, let $\domain_\pmIdx$, $\domain_\fracIdx$, $\domain_\isIdx$ and $\domain_\juncIdx$ be
disjoint partitions of $\domain$ representing the bulk porous medium, the fractures, the intersections
of fractures and the intersections of intersections, respectively. For the sake of readability, we
will use the term junctions to refer to the intersections of intersections.
The external subdomain boundaries are denoted by
$\surface_\pmIdx \eqbydef \partial \domain_\pmIdx \cap \partial \domain$,
$\surface_\fracIdx \eqbydef \partial \domain_\fracIdx \cap \partial \domain$ and
$\surface_\isIdx \eqbydef \partial \domain_\isIdx \cap \partial \domain$.
For simplicity, we assume that junctions do not occur on the external boundary, \ie it is
$\surface_\juncIdx \eqbydef \partial \domain_\juncIdx \cap \partial \domain = \emptyset$.
The interfaces between subdomains are defined as
$\interface_{k, l} = \partial \domain_k \cap \partial \domain_l$, where
$k \neq l$ and $k, l \in \{ \pmIdx, \fracIdx, \isIdx, \juncIdx \}$.
We furthermore define the normal vectors on the interfaces $\interface_{k, l}$ to be pointing outwards
of the domains $\domain_k$. A two-dimensional illustration of this decomposition is
depicted in~\cref{fig:modeldomain_equi}. Let us further decompose the external sub-domain boundaries
$\surface_\pmIdx = \surface_\pmIdx^{\diriIdx} \cup \surface_\pmIdx^{\neumIdx}$,
$\surface_\fracIdx = \surface_\fracIdx^{\diriIdx} \cup \surface_\fracIdx^{\neumIdx}$ and
$\surface_\isIdx = \surface_\isIdx^{\diriIdx} \cup \surface_\isIdx^{\neumIdx}$ to
indicate whether Dirichlet (superscript $\diriIdx$) or Neumann (superscript $\neumIdx$)
boundary conditions are imposed.
We require $\meas{\surface_l^{\diriIdx}} > 0$ for at least one subdomain
$l \in \{ \pmIdx, \fracIdx, \isIdx \}$.\\

The decomposition of the domain as described above
might lead to the subdomains being composed of several disconnected sets.
For instance, in the situation illustrated in~\cref{fig:modeldomain_equi}, the bulk domain is
partitioned into three and the fracture domain into four subsets.
We make this apparent in the notation by equipping the subdomains with the superscript $i$ indicating
the index within the partition of a subdomain, \ie it is
$\domain_k = \bigcup_{i \in \sdPartitionSet_k} \domain_k^i$, where
$\sdPartitionSet_k = \{0, 1, \dots, \numSubSets_{k} \}$ and
$k \in \{ \pmIdx, \fracIdx, \isIdx, \juncIdx \}$.
Consequently, let us introduce the interfaces between subsets
$\interface_{k, l}^{i, j} = \partial \domain_k^i \cap \partial \domain_l^j$,
$k \neq l$, $i \in \sdPartitionSet_k$ and $j \in \sdPartitionSet_l$, where we note that
$\interface_{k, l} = \bigcup_{\substack{i \in \sdPartitionSet_k \\ j \in \sdPartitionSet_l}} \interface_{k, l}^{i, j}$.
In this setting, we state the equi-dimensional formulation for single-phase flow:
\begin{subequations}
  \label{eq:prob_equi}
  \begin{align}
    \darcyVel_k^i + \perm_k^i \grad \head_k^i &= 0, \label{eq:prob_equi_darcy}                    \\
    \divergence \darcyVel_k^i &= \source_k^i,
    &&\mathrm{in} \, \domain_k^i, \label{eq:prob_equi_mass}                                       \\
    \head_k^i &= \head_l^j, \label{eq:prob_equi_if_p}                                             \\
    \darcyVel_k^i \scal \n &= \darcyVel_l^j \scal \n,
    &&\mathrm{on} \, \interface_{k, l}^{i, j}, \label{eq:prob_equi_if_flux}                       \\
    \head_k^i &= \head_k^\diriIdx,
    &&\mathrm{on} \, \surface_k^\diriIdx \cap \partial \domain_k^i, \label{eq:prob_equi_diri}     \\
    \darcyVel_k^i \scal \n &= f_k^i,
    &&\mathrm{on} \, \surface_k^\neumIdx \cap \partial \domain_k^i, \label{eq:prob_equi_neumann}
  \end{align}
\end{subequations}
where $i \in \sdPartitionSet_k$, $j \in \sdPartitionSet_l$,
$k,l \in \{ \pmIdx, \fracIdx, \isIdx, \juncIdx \}$ and $l \neq k$.
\Cref{eq:prob_equi_darcy} states that Darcy's Law is used in all subdomains to relate
the flux $\darcyVel$ (given in \si{\cubic\meter\per\second\per\meter\squared}) to gradients
in the hydraulic head $\head$ (given in \si{\meter}), with $\perm$ being the hydraulic
conductivity given in \si{\meter\per\second}. \Cref{eq:prob_equi_mass} states the conservation
of mass in each subdomain, with $\source$ denoting sources or sinks in \si{\per\second}.
\Cref{eq:prob_equi_if_p,eq:prob_equi_if_flux} state the interface conditions between different
subdomains, which comprise of the continuity of the hydraulic head and fluxes. Finally,~\cref{eq:prob_equi_diri,eq:prob_equi_neumann}
state the Dirichlet and Neumann boundary conditions on the external boundaries.\\

In the subsequent section, we present the mixed-dimensional formulation that approximates
the model given in~\eqref{eq:prob_equi}. The notation that we will use is strongly
influenced by the presentations given
in~ \citet{Flemisch2018Benchmarks,nordbotten2019unified}.

\section{Mixed-dimensional model}
\label{sec:mixedModel}
\begin{figure}[h]
    \centering
    \includegraphics[width=0.9\textwidth]{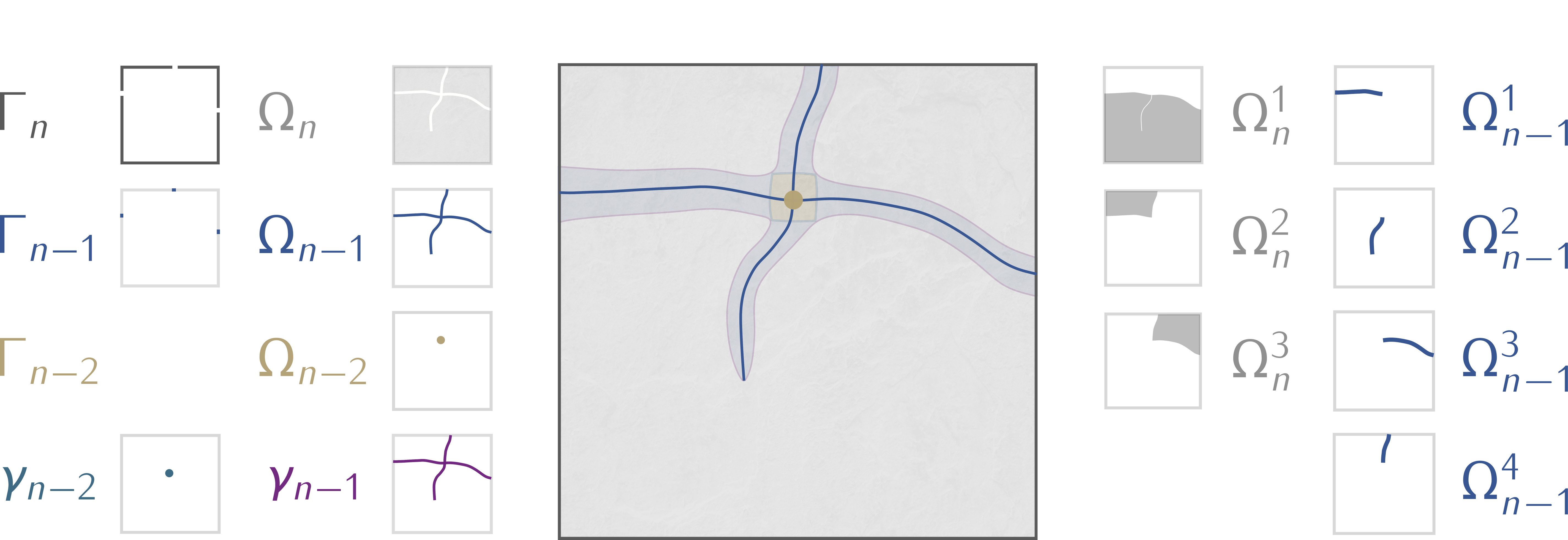}
    \caption{\textbf{Mixed-dimensional domain decomposition}. The domain $\domain$, illustrated in~\cref{fig:modeldomain_equi} and depicted transparent here, is approximated by a mixed-dimensional decomposition consisting of an $n$-dimensional subdomain $\domain_n$ for the bulk porous medium and the lower-dimensional subdomains $\domain_{n-1}$ for the fractures and $\domain_{n-2}$ for the intersections of fractures. Each subdomain is potentially split into disconnected subsets $\domain_d^i$, which are illustrated on the right. The interfaces between $d$- and ($d-1$)-dimensional subdomains are denoted by $\interface_{d-1}$, while the external boundaries of the sub-domains are denoted by $\surface_d$. Note that in two-dimensional settings, junctions of intersections do not occur which is why they are not visualized here.}
    \label{fig:modeldomain_mixeddim}
\end{figure}
We now want to approximate the subdomains of the equi-dimensional setting for the fractures,
fracture intersections and junctions of intersections, \ie $\domain_\fracIdx$, $\domain_\isIdx$ and
$\domain_\juncIdx$, by lower-dimensional geometries. More precisely, the fractures are approximated
by $2$-dimensional entities that are constructed such that they follow the planes that describe the
center of the fractures. Correspondingly, the intersections of fractures are approximated by
$1$-dimensional segments that describe center lines of the intersection regions between fractures and
the junctions are approximated by $0$-dimensional points located at their center.
In this mixed-dimensional setting, we will from now on refer to the sub-domains with the subscript
$d$ indicating the dimension, \ie the original domain $\domain$ is split into a $3$-dimensional bulk
domain $\domain_3$, as well as the lower-dimensional domains $\domain_2$, $\domain_1$ and $\domain_0$
for the fractures, their intersections and junctions of intersections, respectively.
The external boundaries of the subdomains are denoted by
$\surface_d = \partial \domain_d \cap \partial \domain$, where we note that
$\partial \domain_0 = \surface_0 = \emptyset$.
Please note that in this mixed-dimensional setting, there can also occur immersed boundaries
(fracture tips) $\fracTips_{2} = \partial \domain_{2} \setminus \surface_{2}$ and
$\fracTips_{1} = \partial \domain_{1} \setminus \surface_{1}$, and we extend the notation to
the highest-dimensional domain by defining $\fracTips_3 = \emptyset$.
A two-dimensional illustration of this decomposition of $\domain$ is given
in~\cref{fig:modeldomain_mixeddim}.\\

A first simplification with respect to the equi-dimensional setting is that we only consider
interfaces between subdomains of codimension one, which we denote by
$\interface_d = \partial \domain_{d+1} \cap \domain_d$, $0 \leq d < 3$,
\ie $\interface_2$, $\interface_1$ and $\interface_0$ denote the interface between the bulk domain
and the fractures, the interface between the fractures and their intersections and the interface
between intersections and junctions. The normal vectors on these interfaces are defined such that
they point outwards of the higher-dimensional domains.\\

As for the equi-dimensional case, the subdomains
$\domain_d$ might comprise of several disconnected subsets, that is,
$\domain_d^i$, $i \in \sdPartitionSet_d$, is a subset of $\domain_d$, and we denote by
$\interface_d^{i,j} = \partial \domain_{d+1}^i \cap \domain_d^j$,
$i \in \sdPartitionSet_{d+1}$, $j \in \sdPartitionSet_d$,
an interface between two such subsets of codimension one.
Let us further introduce for each $\domain_d^i$, $0 \leq d < 3$,
the set of indices of the neighboring higher-dimensional subsets
$\domain_{d+1}^j$:
\begin{equation}
    \neighborDomainSet_d^{i, \uparrow} = \{j \in \sdPartitionSet_{d+1}: \interface_d^{j,i} \neq \emptyset \}.
    \label{eq:neighborSetHigher}
\end{equation}

\paragraph{Cross section-averaged quantities}
At each point $\position \in \domain_d$, $d < 3$, let
$\crossSection_d = \crossSection_d \left( \position \right)$
denote the cross section that describes the geometry of the equi-dimensional counterpart
around that position. That is, $\crossSection_{2}$ is a segment orthogonal to $\domain_2$,
whose measure is the fracture aperture $\aperture$, \ie $\meas{\crossSection_{2}} = \aperture$.
Correspondingly, $\meas{\crossSection_{1}}$ is the cross-sectional area of an intersection
region of two or more fractures. Thus, we note that $\meas{\crossSection_d}$ is given in
$\si{\meter^{3-d}}$.
Let us now define the cross section-averaged hydraulic head as
\begin{equation}
      \head_d^i = \frac{1}{\meas{\crossSection_d^i}}
                  \int_{\crossSection_d^i}
                      \head_{k}^i \,
                  \mathrm{d} x,
    \label{eq:averagedHead}
\end{equation}
where $\head_k^i$ refers to the hydraulic head in the equi-dimensional setting, and
we consider the tuples
$\left(d, k \right) \in \indexTupleSet \eqbydef \{ \left( 2, \fracIdx \right),
                                                   \left( 1, \isIdx   \right) \}$.
Correspondingly, we denote by $\source_d^i$ the cross section average of the sources
$\source_k^i$ appearing in~\cref{eq:prob_equi}. It will become clear in the following
paragraphs why the tuple $\left( 0, \juncIdx \right)$ is not considered in $\indexTupleSet$.

\paragraph{Cross section-integrated balance equations}
Integration of the mass balance equation~\cref{eq:prob_equi_mass} over the cross sections
$\crossSection_d$, $1 \leq d < 3$, yields the PDEs describing the flow along the lower-dimensional domains
$\domain_d$, for which a detailed derivation can be found in~\citet{martin2005modeling}.
This leads to the lower-dimensional mass balances
\begin{equation}
    \fracLocalDiv{d} \left( \meas{\crossSection_d^i} \darcyVel_d^i \right)
    = \meas{\crossSection_d^i} \source_d^i
      + \jump{\darcyVel_{d+1} \scal \n}^i_d,
    \label{eq:averageMassBalance}
\end{equation}
with $\darcyVel_d^i$ being the cross section-averaged flux, and
$\fracLocalDiv{d} \opPlaceHolder$
denoting the divergence operator in tangential direction of $\domain_d$.
The last term in~\cref{eq:averageMassBalance} describes the mass transfer with
the surrounding higher-dimensional domain and appears as an additional source/sink
term. It sums up all fluxes leaving or entering the neighboring ($d+1$)-dimensional
subdomain, \ie we have:
\begin{equation}
    \jump{\darcyVel_{d+1} \scal \n}^i_d
    = \sum_{j \in \neighborDomainSet_d^{i, \uparrow}}
        \sum_{l = 1}^{\numFractureSides_d^{i, j}}
            \left( \darcyVel_{d+1}^j \scal \n \right)_l,
    \label{eq:fluxJump}
\end{equation}
where we note that $\numFractureSides^{i, j}_d \in \{1, 2 \}$.
That is because for a point on the lower-dimensional
subdomain $\domain_d$, there might be two sides of it that are associated with the
same higher-dimensional subset $\domain_{d+1}^j$.
Such situations occur, for instance, around immersed fracture tips,
as illustrated by subdomain $\domain_{n-1}^2$ in~\cref{fig:modeldomain_equi},
which is fully immersed in $\domain_n^1$.\\

Under the assumption that $\perm_k$ is constant over the cross section $\crossSection_d$,
$\left(d, k \right) \in \indexTupleSet$, the cross section-averaged flux $\darcyVel_d^i$ is
given by
\begin{equation}
    \darcyVel_d^i = - \perm_d^i \fracLocalGrad{d} \head_d^i,
    \label{eq:averageFlux}
\end{equation}
where $\fracLocalGrad{d} \opPlaceHolder$ is the gradient operator in tangential
direction of $\domain_d$, and $\perm_d^i$ is the tangential part of $\perm_k^i$,
that is, $\perm_d^i = \perm_k^i - \N_d \, \perm_k^i$,
with $\N_d$ being the projection into the space normal to $\domain_d$.
For the fracture domain it is $\N_2 = \n_2 \otimes \n_2$, where $\n_2$ is a vector
normal to $\domain_2$ and $\otimes$ is the dyadic produc.
For intersections, the normal space is defined by two normal vectors to $\domain_1$
that are themselves orthogonal to each other: $\N_1 = \n_{1, 1} \otimes \n_{1, 2}$
Correspondingly, let us denote the normal part by
$\perm_d^{i, \perp} = \n \otimes \n \, \perm_k^i$.
Please note that we assume the normal part $\perm_d^{i, \perp}$ to be isotropic
within the normal space. This is a necessary requirement for the mixed-dimensional
coupling condition~\eqref{eq:prob_mixed_if_flux} presented later.

\paragraph{Junctions of fracture intersections}
In the model presented in this work, we assume continuity of the hydraulic head
and fluxes across junctions of intersections:
\begin{subequations}
    \label{eq:junctionConditions}
    \begin{align}
        \head_1^i
        &= \head_1^j, &&\mathrm{on} \, \domain_0^k, \quad k \in \sdPartitionSet_0, \quad i, j \in \neighborDomainSet_0^{k, \uparrow}, \label{eq:junctionCondHead} \\
        0 &= \sum_{i \in \neighborDomainSet_0^{k, \uparrow}} \darcyVel_1^i \scal \n, &&\mathrm{on} \, \domain_0^k, \quad k \in \sdPartitionSet_0. \label{eq:junctionCondFlux}
    \end{align}
\end{subequations}
Thus, we assume that the junctions are significantly more permeable than the
adjacent fracture intersection branches. This assumption is mainly motivated by
limitations of our current implementation, and other authors have used different
conditions that capture possible jumps in the hydraulic head across junctions of
intersections \citep{nordbotten2019unified,keilegavlen2020porepy}.

\paragraph{Mixed-dimensional formulation}
Until now, we have considered a three-dimensional domain, taking into account flow
in the bulk medium, the fractures and fracture intersections, while assuming
continuity of the hydraulic head and the flux across junctions of fracture intersections.
In order to generalize the model, let $n \in \{2, 3 \}$ denote the dimension of
the bulk domain, and $\minDim$, $1 \leq \minDim \leq n-1$, the dimension of the
lowest-dimensional subdomain along which flow is taken into account. Moreover,
we use~\cref{eq:averageMassBalance} to also describe flow in the bulk domain by
using the following extensions to the notation:
$\fracLocalDiv{n} \opPlaceHolder \eqbydef \divergence \opPlaceHolder$,
$\fracLocalGrad{n} \opPlaceHolder \eqbydef \grad \opPlaceHolder$,
$\head_n \eqbydef \head_\pmIdx$, $\darcyVel_n \eqbydef \darcyVel_\pmIdx$
$\perm_n \eqbydef \perm_\pmIdx$,
$\meas{\crossSection_n} = 1$
and
$\jump{\darcyVel_{n+1} \scal \n }_n \eqbydef 0$.
The mixed-dimensional problem formulation then reads, for $\minDim \leq d \leq n$:
\begin{subequations}
  \label{eq:prob_mixed_mass_darcy}
  \begin{align}
    \darcyVel_d^i + \perm_d^i \fracTangGrad{d} \head_d^i &= 0, \label{eq:prob_mixed_darcy}          \\
    \fracTangDiv{d} \left( \meas{\crossSection_d^i} \darcyVel_d^i \right)
    &= \meas{\crossSection_d^i} \source_d^i
    + \jump{ \darcyVel_{d+1} \scal \n }_d^i,
    &&\mathrm{in} \,\, \domain_d^i, \label{eq:prob_mixed_balance}
  \end{align}
\end{subequations}
together with the boundary conditions
\begin{subequations}
  \label{eq:prob_mixed_bcs}
  \begin{align}
    \head_d^i &= \head_d^\diriIdx,
    &&\mathrm{on} \,\, \surface_d^\diriIdx \cap \partial \domain_d^i, \label{eq:prob_mixed_diri}    \\
    \darcyVel_d^i \scal \n &= f_d^i,
    &&\mathrm{on} \,\, \surface_d^\neumIdx \cap \partial \domain_d^i, \label{eq:prob_mixed_neumann} \\
    \darcyVel_d^i \scal \n &= 0,
    &&\mathrm{on} \,\, \fracTips_d \cap \partial \domain_d^i, \label{eq:prob_mixed_neumann_tips}
  \end{align}
\end{subequations}
the conditions~\eqref{eq:junctionConditions} at junctions of the lowest-dimensional domain:
\begin{subequations}
  \label{eq:prob_mixed_junctions}
  \begin{align}
       \head_\minDim^i
       &= \head_\minDim^j,
       &&\mathrm{on} \, \domain_{\minDim-1}^k, \quad k \in \sdPartitionSet_{\minDim-1}, \quad i, j \in \neighborDomainSet_{\minDim-1}^{k, \uparrow}, \label{eq:prob_mixed_junction_p} \\
       0 &= \sum_{i \in \neighborDomainSet_{\minDim-1}^{k, \uparrow}} \darcyVel_\minDim^i \scal \n,
       &&\mathrm{on} \, \domain_{\minDim-1}^k, \quad k \in \sdPartitionSet_{\minDim-1}, \label{eq:prob_mixed_junction_flux}
  \end{align}
\end{subequations}
and the interface conditions between subdomains,
with $\minDim \leq \tilde{d} < n$:
\begin{subequations}
  \label{eq:prob_mixed_if}
  \begin{align}
    - \n^T \perm_{\tilde{d}}^{j, \perp} \n \,
      \frac{ \head_{\tilde{d}}^{j} -
             \head_{\tilde{d}+1}^i
           }{\fracLength_{\tilde{d}}^{j, \perp}}
    &= \darcyVel_{\tilde{d}+1}^i \scal \n,
    &&\mathrm{on} \,\, \interface_{\tilde{d}}^{i, j}
                       \in \blockInterfaceSet_{\tilde{d}}, \label{eq:prob_mixed_if_flux} \\
    \head_{\tilde{d}}^j &= \head_{\tilde{d}+1}^i,
    &&\mathrm{on} \,\, \interface_{\tilde{d}}^{i, j}
                       \in \condInterfaceSet_{\tilde{d}}. \label{eq:prob_mixed_if_p}
  \end{align}
\end{subequations}
As can be seen in the above conditions, we distinguish two types of interfaces.
In the set $\condInterfaceSet_{\tilde{d}}$ we collect the interfaces for which
the permeability in the lower-dimensional neighboring domain is much larger than in
the higher-dimensional neighbor. In this case, the jump in hydraulic head across
the lower-dimensional domain may be neglected, which is expressed by
condition~\eqref{eq:prob_mixed_if_p}. In particular, we want to use this condition
on interfaces to open fractures or intersections, which are often modeled by
using an aperture-dependent tangential permeability \citep{Matthai2007,ucar2018three},
but for which a normal permeability
$\perm_{\tilde{d}}^{j, \perp}$, $j \in \sdPartitionSet_{\tilde{d}}$, cannot be defined.
The remaining interfaces are collected in the set $\blockInterfaceSet_{\tilde{d}}$,
on which the continuity of the flux and the hydraulic head is enforced via
condition~\eqref{eq:prob_mixed_if_flux}, where the normal flux inside the fracture
or fracture intersections is expressed by means of the finite difference given on
the left hand side. Therein, $\fracLength_{\tilde{d}}^{j, \perp}$
is a characteristic length scale associated with the distance of $\domain_{\tilde{d}}$
to the corresponding interface in the equi-dimensional setting. For the fractures,
this is half of the aperture, \ie
$\fracLength_{n-1}^{j, \perp} = \aperture^j/2 = \meas{\crossSection_{n-1}^j}/2$,
while for fracture intersections (when $n=3$), it might be different for
each neighboring fracture branch. In this work, we use
$\fracLength_{1}^{j, \perp} = \sqrt{\meas{\crossSection_{1}}}/2$, but more complex
approaches can be used if information on the intersection geometry is available
(see \eg \citet{Walton2017Elimination}). Note that condition~\eqref{eq:prob_mixed_if_flux}
corresponds to the one given in \citet{martin2005modeling} for $\xi = 1.0$,
where $\xi$ is a parameter of the model presented therein.

\section{Discretization schemes}
\label{sec:discretization}
In this work, we compare several finite-volume based approaches for the discretization
of the model problem stated in~\eqref{eq:prob_mixed_mass_darcy}-\eqref{eq:prob_mixed_if}
on the basis of a set of test cases. Several schemes are taken into account, some of which
have been presented in the literature, and some of which contain, to the authors' knowledge,
novel features in the way they are presented in this work. In the following, we want to provide
a short overview over the considered schemes, and subsequently, we will provide more details
on the vertex-centered approaches introduced in this work. All schemes are implemented into the
open-source simulation framework \dumux~\citep{Dumux,Kochetal2020Dumux}, and the source
code to the numerical examples presented in this work can be found at
\href{https://git.iws.uni-stuttgart.de/dumux-pub/glaeser2020b}{git.iws.uni-stuttgart.de/dumux-pub/glaeser2020b}.
Note that a monolithic solution strategy is chosen, that is, all discrete equations
are assembled into a single, block-structured, system matrix.

\paragraph{\tpfaDfm}
This scheme is based on a cell-centered finite-volume scheme using the two-point flux
approximation (\tpfa), and was first introduced in~\citet{Karimi2004TpfaDfm}.
Further details can be found in~\citet{Sandve2012MpfaDfm}.

\paragraph{\mpfaDfm}
Also a cell-centered scheme, it differs from the \tpfaDfm scheme in the way the
fluxes between the grid cells are computed. In particular, it uses the multi-point
flux approximation (\mpfa) technique presented in \citet{Aavatsmark2002} for
unfractured media, and extensions to fractured porous media as presented
in~\citet{Sandve2012MpfaDfm, Ahmed2015MpfaDfm}. Moreover, it was used in the context
of isothermal and non-isothermal two-phase flow in fractured porous media
in~\citet{Glaeser2017MpfaDfmTwoP,glaser2019hybrid}.

\paragraph{\boxDfm}
Based on the vertex-centered \boxScheme scheme (see \eg~\citet{Hackbusch89,helmig97}),
the \boxDfm scheme was introduced in~\citet{Reichenberger2006mixed} for two-phase
flow in fractured porous media. The reduction to single-phase flow is straightforward,
and we refer to~\citet{Reichenberger2006mixed,Tatomir2012} for further details.
It should be mentioned that this scheme does not exactly
solve~\eqref{eq:prob_mixed_mass_darcy}-\eqref{eq:prob_mixed_if}. Instead, continuity
of the hydraulic head across the fractures (as expressed in~\cref{eq:prob_mixed_if_p})
is assumed. Therefore, no additional degrees of freedom with respect to the standard
\boxScheme scheme are introduced to describe the state inside the fracture.
The effect of the fractures is incorporated by providing additional
connectivity within the single mass balance equation that is solved for the bulk medium.
This means, in particular, that~\cref{eq:prob_mixed_balance} is not
solved for the fracture domain, and that the condition~\cref{eq:prob_mixed_if_flux}
is not fulfilled at matrix-fracture interfaces. Moreover, the implementation as
provided in \dumux does not allow for flow along fracture intersections to be
considered.

\paragraph{\eboxDfm}
This scheme is an extension of the \boxScheme scheme in the sense that additional
degrees of freedom are introduced in and around the fractures, which allow for an incorporation
of the conditions~\eqref{eq:prob_mixed_if} at the matrix-fracture and fracture-intersection
interfaces and the solution of the mass balance equations~\eqref{eq:prob_mixed_balance} for the
lower-dimensional domains, and which make it possible to capture the jump in fluxes
and hydraulic head that can occur across low-permeable fractures. This scheme will
be presented in more detail in~\cref{sec:eboxDfm}.

\paragraph{\eboxMortarDfm}
We introduce the extended box scheme, the \eboxDfm scheme, to overcome the drawbacks
of the \boxDfm scheme with respect to the incorporation of the interface conditions
given in~\cref{eq:prob_mixed_if}. However, the conditions~\eqref{eq:prob_mixed_if_p}
can only be incorporated weakly, as will be seen in~\cref{sec:eboxDfm}. To this end,
the \eboxMortarDfm scheme was developed, in which an additional mortar variable
representing the normal flux, is introduced on the interfaces between domains of
different dimensionality. This enables a strong incorporation of the
conditions~\eqref{eq:prob_mixed_if}. Details on this scheme will be provided
in~\cref{sec:eboxDfmMortar}.\\

In the following sections, we want to outline the basic concepts of the \eboxDfm
and the \eboxMortarDfm schemes. \Cref{sec:mesh} presents the construction of the
finite-volume discretizations of the subdomains, which is common to both schemes,
and subsequently, the schemes are presented in \cref{sec:eboxDfm,sec:eboxDfmMortar}.

\subsection{Computational mesh}
\label{sec:mesh}
\begin{figure}[ht]
    \centering
    \includegraphics[width=0.99\textwidth]{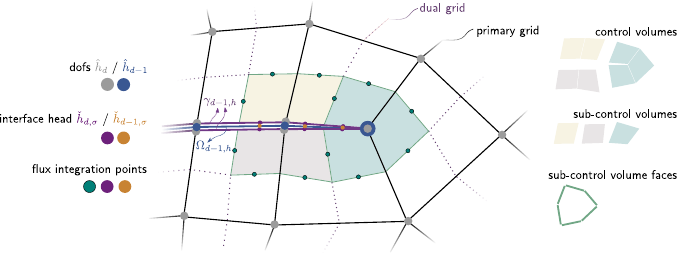}
    \caption{\textbf{Illustration of the \eboxDfm scheme}. Shown is an exemplary configuration around a fracture tip, depicting the distribution of the degrees of freedom, the subdivision of the control volumes into sub-control volumes, and the points at which fluxes are assembled. For illustration purposes, the fracture surfaces have been moved away from each other, separating the interface $\interface_{d-1}$ from the lower-dimensional domain $\domain_{d-1}$. They actually geometrically coincide in the discrete representation.}
    \label{fig:eboxDfm}
\end{figure}
Let us introduce the primary discretizations $\mesh_d$, $\minDim \leq d \leq n$,
with elements $\element \in \mesh_d$ such that
$\bar{\domain}_{d, h} \equiv \bigcup_{\element \in \mesh_d} \bar{\element}$ are discrete approximations of $\domain_d$.
Correspondingly, let $\mesh_d^i$ be the subset of $\mesh_d$ associated with the discrete subdomain
$\domain_{d, h}^i$. 
Moreover, let us denote with $\meshFacetSet_d$ the set of faces of the primary discretization,
and define the subsets $\meshFacetSet_\element$ such that for each element $\element \in \mesh_d$
it is $\partial \element = \bigcup_{\elemFacet \in \meshFacetSet_\element} \bar{\elemFacet}$. We require
the grids to be conforming such that the elements of the lower-dimensional domains coincide with
faces of the next higher-dimensional domains, that is, we have $\forall \element \in \mesh_d$,
$d < n$: $\exists \elemFacet \in \meshFacetSet_{d+1}$ such that $\elemFacet \equiv \element$. This is
illustrated in~\cref{fig:eboxDfm}.\\

In addition to the primary grid, let us introduce the tesselations $\cvSet_d$, composed of control
volumes $\cv \in \cvSet_d$ with measure $\meas{\cv} > 0$ such that
$\bigcup_{\cv \in \cvSet_d} \bar{\cv} \equiv \bar{\domain}_{d, h}$.
As for the primary discretization, we denote by $\cvSet_d^i$ the subset associated with the
discrete subdomain $\domain_{h, d}^i$.
While for cell-centered schemes it is $\cvSet_d \equiv \mesh_d$, the \boxScheme discretization
is defined on a dual grid such that the elements and control volumes do not coincide. The dual grid
is constructed by connecting the barycenters of the cells and edges (and in 3d also the faces) of the
primary discretizations, as illustrated in~\cref{fig:eboxDfm} for a two-dimensional setting, and the
control volumes are defined by the emerging bounded regions around the grid vertices. In contrast to
the standard \boxScheme scheme, the control volumes in the \boxDfm scheme might be split into two
or more control volumes by intersecting fractures or intersections, which can be seen in the
illustration given in~\cref{fig:eboxDfm} for the two central control volumes. Each control volume has
a unique degree of freedom associated with it, and correspondingly, additional degrees of freedom
are introduced for each side of a fracture or intersection.
Please note that this is not the case at fracture tips,
at which a single degree of freedom per domain describes the hydraulic state.\\

Furthermore, we denote by $\faceSet_d$ the sets of faces of the tesselations $\cvSet_d$, and define
the subsets $\faceSet_{\cv}$ such that for a control volume $\cv \in \cvSet_d$ it is
$\partial \cv = \bigcup_{\face \in \faceSet_{\cv}} \bar{\face}$.
The corresponding unit normal vectors $\n_{\cv}^{\face}$ are defined to be pointing outwards of $\cv$.
Let us further define the subsets $\faceSet_{d, \couplIdx}$, $\faceSet_{d, \internalIdx}$ and
$\faceSet_{d, \externalIdx}$ of coupling faces that overlap with a ($d-1$)-dimensional control volume
$\cv \in \cvSet_{d-1}$, internal faces and faces on the outer domain boundary, respectively.
Moreover, we subdivide the set of coupling faces such that
$\faceSet_{d, \couplIdx} = \faceSet_{d, \condInterfaceSet} \cup \faceSet_{d, \blockInterfaceSet}$,
where $\faceSet_{d, \condInterfaceSet}$ collects all coupling faces on which
condition~\eqref{eq:prob_mixed_if_p} is enforced, while
$\faceSet_{d, \blockInterfaceSet}$ contains all faces on which~\cref{eq:prob_mixed_if_flux} holds.\\

As a result of the construction of the dual grid, the control volumes overlap with
several elements of the primary grid. We denote with $\scv$ a sub-control
volume of the control volume $\cv$, where each $\scv$ can be associated to a primary
grid element $\element$ such that $\scv = \cv \cap \element$.
Let $\scvSet_d$ be the set of sub-control volumes of a discretization together with the subset
$\scvSet_\cv$ such that $\bar{\cv} = \bigcup_{\scv \in \scvSet_\cv} \bar{\scv}$.

\subsection{Vertex-centered, discontinuous, finite-volume scheme (\eboxDfm)}
\label{sec:eboxDfm}

As in the standard \boxScheme scheme, piecewise (per primary grid element) linear basis functions are used
to approximate the hydraulic head and its gradients within the primary grid cells. For instance,
the discrete approximation of the hydraulic head at a position $\position \in \domain_{d, h}^i$
is expressed by:
\begin{equation}
  \tilde{\head}_{d}^i \left( \position \right)
  =
  \sum_{\cv \in \cvSet_d^i}
    \hat{\head}_{\cv} \, \basisFunc_{\cv} \left( \position \right),
  \label{eq:discreteHeadGlobal}
\end{equation}
where $\basisFunc_\cv$ refers to the basis function that corresponds to the degree of freedom
associated with the control volume $\cv$, and $\hat{\head}_\cv$ is the nodal value of the
hydraulic head at that degree of freedom. The basis functions are defined as for continuous,
linear finite elements, that is,
$\basisFunc_\cv \left( \hat{\position}_\cv \right) = 1$, with $\hat{\position}_\cv$ being
the position of the degree of freedom associated with control volume $\cv$ (which is the
position of the associated grid vertex), and
$\basisFunc_\cv \left( \hat{\position}_\otherCV \right) = 0$, with
$\cv, \otherCV \in \cvSet_d$, $\cv \neq \otherCV$.
Furthermore, the basis functions $\basisFunc_\cv$ are only non-zero within primary
grid elements that the control volume $\cv$ overlaps with.
Let us denote with $\cvSet_\element$ the set of control volumes for which it is
$\cv^\circ \cap \element^\circ \neq \emptyset$, with $\cv \in \cvSet_d$ and $\element \in \mesh_d$.
The discrete hydraulic head and its gradient, evaluated at a position
$\position \in \element$, $\element \in \mesh_d^i$, can then be written as:
\begin{subequations}
  \label{eq:discreteApproximations}
  \begin{align}
  \tilde{\head}_{d}^i \left( \position \right)
  =
  \sum_{\cv \in \cvSet_\element}
    \hat{\head}_{\cv} \, \basisFunc_{\cv} \left( \position \right), \label{eq:discHead} \\
  \fracLocalGrad{d} \tilde{\head}_{d}^i \left( \position \right)
  =
  \sum_{\cv \in \cvSet_\element}
    \hat{\head}_{\cv} \, \grad \basisFunc_{\cv} \left( \position \right). \label{eq:discGrad}
  \end{align}
\end{subequations}

\paragraph{Discrete mass balance equation}
Integration of~\cref{eq:prob_mixed_balance} over a control volume $\cv \in \cvSet_d^i$, applying the Gauss
divergence theorem and splitting the resulting surface integral over the boundary of $\cv$ into a sum
over the faces $\face \in \faceSet_{\cv}$, yields:
\begin{equation}
    \sum_{\face \in \faceSet_{\cv}}
        \int_{\face}
            \meas{\crossSection_d^i} \left( \darcyVel_d^i \scal \n_{\cv}^{\face} \right)
        \, \surfIntegrand
    =
    \int_{\cv}
        \meas{\crossSection_d^i} \source_{d}^i
        + \jump{ \darcyVel_{d+1} \scal \n  }_d^i
    \, \volIntegrand.
    \label{eq:massbal_cv}
\end{equation}
Moreover, the volume integrals in~\cref{eq:massbal_cv} can be split into a sum
over the sub-control volumes:
\begin{equation}
    \sum_{\face \in \faceSet_{\cv}}
        \int_{\face}
            \meas{\crossSection_d^i} \left( \darcyVel_d^i \scal \n_{\cv}^{\face} \right)
        \, \surfIntegrand
    =
    \sum_{\scv \in \scvSet_\cv}
        \int_{\scv}
            \meas{\crossSection_d^i} \source_{d}^i
            + \jump{ \darcyVel_{d+1} \scal \n  }_d^i
        \, \volIntegrand.
    \label{eq:massbal_scv}
\end{equation}
Let us now use the notation $\opPlaceHolder_\scv$ and $\opPlaceHolder_\face$
to denote quantities that are evaluated for a sub-control volume $\scv$ or a
face $\face$, where we now omit the subscript related to the dimension, as the
discrete entities carry the notion of the dimensionality of the grid on which
they are defined. Following the same argumentation, we omit the superscript $i$
refering to the subset $\domain_{d,h}^i \subset \domain_{d,h}$ on which a variable is defined.
We will use this notation with any discrete entity, and it means that the quantity
is taken as constant on that entity. With this, we can write the discrete balance
equation for a control volume $\cv \in \cvSet_d^i$ as follows:
\begin{equation}
    \sum_{\face \in \faceSet_{\cv}}
        \discFlux_\cv^\face
    =
    \sum_{\scv \in \scvSet_\cv}
        \meas{\scv}
        \meas{\crossSection_\scv}
        \source_\cv
        + \jump{ \discFlux_{d+1} }_\scv,
    \label{eq:massbal_discrete}
\end{equation}
where we have introduced the discrete flux
\begin{equation}
    \discFlux_\cv^\face
    \approx
    - \int_\face
        \meas{\crossSection_d^i} \left(\n_\cv^\face \right)^T
                                 \perm_d^i
                                 \fracLocalGrad{d} \head_d^i
      \, \surfIntegrand
    =
    \int_\face
        \meas{\crossSection_d^i} \left( \darcyVel_d^i \scal \n_{\cv}^\face \right)
    \surfIntegrand
    \label{eq:discFluxIntro}
\end{equation}
as an approximation of the flux terms appearing in the sum on the left hand side
of~\cref{eq:massbal_scv}, and where $\jump{ \discFlux_{d+1} }_\scv$
is the jump of the discrete flux in the adjacent higher-dimensional domain across
the sub-control volume $\scv$. That is, it is the sum over the discrete fluxes on
all faces $\otherFace \in \faceSet_{\cv^\uparrow}$, embedded in the control
volumes ${\cv^\uparrow} \in \cvSet_{d+1}$,
which overlap with $\scv$ such that $\scv \equiv \otherFace$.

\paragraph{Discrete flux expressions on internal faces}
Using the expression for the discrete gradient given in~\cref{eq:discGrad}, we
define the discrete flux on internal faces $\face \in \faceSet_{d, \internalIdx}$,
embedded in the control volume $\cv \in \cvSet_d^i$, by:
\begin{equation}
    \discFlux_\cv^\face
    =
    - \meas{\face} \meas{\crossSection_\scv}
      \left(\n_\cv^\face \right)^T
      \perm_\scv
      \fracLocalGrad{d}
      \tilde{\head}_d^i \left( \position_\face \right)
    \label{eq:discFluxInternal}
\end{equation}
where $\position_\face$ is the center point of the face $\face$.
Here, $\scv$ refers to the sub-control volume adjacent
to the face $\face$, that is, $\face \subset \overline{\scv}$.
Please note that we define the permeability and the
cross-sectional measure to be constant per primary grid cell. Thus, on the
sub-control volumes we inherit these properties, that is, we have
$\meas{\crossSection_\scv} = \meas{\crossSection_\element}$ and
$\perm_\scv = \perm_\element$ for
$\element \in \mesh_d^i$,
$\scv \in \scvSet_\cv$,
$\cv \in \cvSet_d^i$,
$\scv \subset \element$.

\paragraph{Discrete flux expressions on coupling faces}
On coupling faces $\face \in \faceSet_{d, \blockInterfaceSet}$,
we want to enforce the condition~\eqref{eq:prob_mixed_if_flux}. To this end, we
substitute the discrete flux by the finite difference given on the left hand
side of~\cref{eq:prob_mixed_if_flux}:
\begin{equation}
    \discFlux_\cv^\face
    =
    - \meas{\face} \meas{\crossSection_\scv}
      \left( \n_\cv^\face \right)^T
      \perm_{\scv^\downarrow}^\perp
      \n_\cv^\face
      \frac{
              \check{\head}_{d-1, \face} - \check{\head}_{d, \face}
           }{
              \fracLength_{\scv^\downarrow}
           }.
    \label{eq:discFluxCouplBlock}
\end{equation}
Here, ${\scv^\downarrow} \in \scvSet_{d-1}$ refers to the lower-dimensional control volume
that coincides with the face $\face$, \ie $\face \equiv {\scv^\downarrow}$, and
$\check{\head}_{d-1, \face} = \tilde{\head}_{d-1} \left( \position_\face \right)$ and
$\check{\head}_{d, \face} = \tilde{\head}_{d} \left( \position_\face \right)$ denote
the discrete hydraulic heads evaluated at the center $\position_\face$ of the face
after~\cref{eq:discHead}.\\

The \eboxDfm scheme does not allow for a strong incorporation of the
conditions~\eqref{eq:prob_mixed_if_p} for the degrees of freedom that live on the
interfaces to highly-permeable or open fractures or fracture intersections, while
guaranteeing that the discrete mass balance equation~\eqref{eq:massbal_discrete}
is fulfilled for the corresponding control volumes. Therefore, the condition can
only be incorporated weakly, for instance, via a penalty term.
However, this approach requires an adequate choice of the penalty parameter, where
for large values, one typically runs into the problem that the conditioning of the
linear system deteriorates and the solution of the problem becomes cumbersome, while
for small values the desired condition of continuity of the hydraulic head might
be significantly violated.
In this work, we want to investigate a different approach.
To this end, let us introduce the map
\begin{equation}
    \codimOneIndexMap : \cvSet_d \rightarrow \cvSet_{d-1},
        \quad \mathrm{s.t.} \,\,
        \hat{\position}_\cv = \hat{\position}_{\codimOneIndexMap \left( \cv \right)}, \,
        \cv \in \cvSet_d,
    \label{eq:codimOneIndexMap}
\end{equation}
which maps to each control volume $\cv \in \cvSet_d$ that touches an interface
to the adjacent ($d-1$)-dimensional domain, the lower-dimensional control volume
$\otherCV = \codimOneIndexMap \left( \cv \right)$ whose
degree of freedom is located at the same geometric location. We then define the
discrete flux on coupling faces $\face \in \faceSet_{d, \condInterfaceSet}$ on
which the condition~\eqref{eq:prob_mixed_if_p} is to be used by
\begin{equation}
    \discFlux_\cv^\face
    =
    - \meas{\face} \meas{\crossSection_\scv}
      \left(\n_\cv^\face \right)^T
      \perm_\scv
      \left[
        \sum_{ \substack{\cv \in \cvSet_\element \\ \hat{\position}_\cv \notin \domain_{d-1, h}} }
            \hat{\head}_{\cv} \, \grad \basisFunc_{\cv} \left( \position_\face \right)
        +
        \sum_{ \substack{\cv \in \cvSet_\element \\ \hat{\position}_\cv \in \domain_{d-1, h}} }
            \hat{\head}_{\codimOneIndexMap \left( \cv \right)} \,
            \grad \basisFunc_{\cv} \left( \position_\face \right)
      \right].
    \label{eq:discFluxCouplCond}
\end{equation}
Thus, for the assembly of the coupling fluxes, we substitute the nodal
hydraulic heads for the degrees of freedom located on an interface with the
nodal values at those locations on the lower-dimensional domain.

\subsection{Vertex-centered, discontinuous, finite-volume scheme with flux mortars (\eboxMortarDfm)}
\label{sec:eboxDfmMortar}

\begin{figure}[h]
    \centering
    \includegraphics[width=0.9\textwidth]{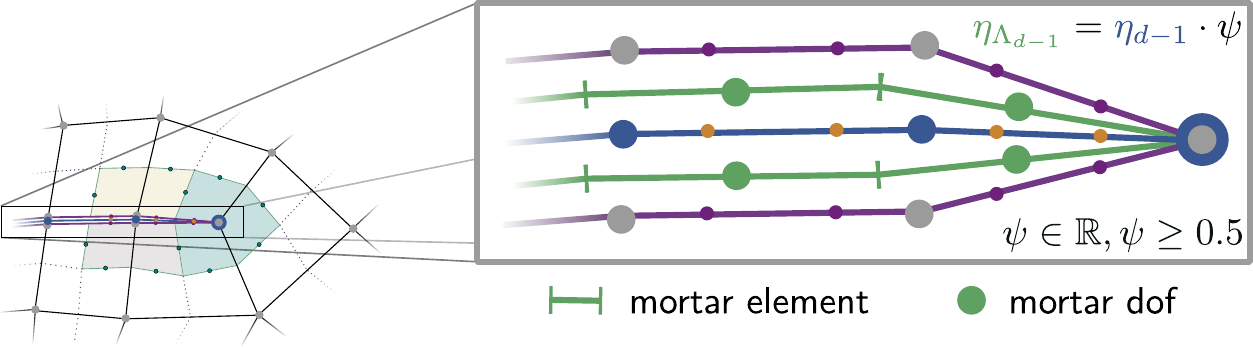}
    \caption{\textbf{Illustration of the \eboxMortarDfm scheme}. In addition to the setting illustrated in~\cref{fig:eboxDfm}, the mortar variable $\mortar_{d-1}$ is introduced at the interfaces between the $d$- and ($d-1$)-dimensional grids. The mortar variable $\mortar$ represents the normal flux across the interface, and is defined on the ($d-1$)-dimensional grid $\mesh_{\mortar_{d-1}}$, which is potentially nonconforming to $\mesh_{d-1}$, the grid used for flow. This is illustrated in the close-up given on the right hand side, where the discretization lengths $\discLength_{\mortar_{d-1}}$ and $\discLength_{d-1}$ of $\mesh_{\mortar_{d-1}}$ and $\mesh_{d-1}$ are related by the factor $\mortarGridFactor$.}
    \label{fig:eboxMortarDfm}
\end{figure}
In the presentation of the \eboxDfm scheme, we have seen that it is not possible
to strongly incorporate the coupling condition~\eqref{eq:prob_mixed_if_p} while
maintaining local mass conservation. In order to overcome this issue, we now want
to present an approach that introduces mortar variables at the interfaces
$\interface_d$ between domains of different dimensionality, which represent the
normal fluxes $\darcyVel_{d+1} \scal \n$ across them.
To this end, let $\mortar_d$, $\minDim \leq d < n$, denote the mortar variable
defined on the interfaces between $d$- and ($d+1$)-dimensional domains. We consider
piecewise constant (per grid element) $\mortar_d$, defined on the mesh $\mesh_{\mortar_{d}}$,
which overlaps with $\mesh_d$ in a potentially nonconforming way.
That is, the characteristic discretization lengths $\discLength_d$ and
$\discLength_{\mortar_d}$ of $\mesh_d$ and $\mesh_{\mortar_d}$ might differ,
and we will use the parameter
$\mortarGridFactor = \discLength_{\mortar_d} / \discLength_d$
to describe their ratio. An illustration of this is given in~\cref{fig:eboxMortarDfm}.\\

Due to the potential non-conformity of the mortar meshes with respect to the adjacent
meshes used for the discretization of the flow equations, we need to project data
from and to the mortar grid. In particular, let us denote by
$\projection_{\mortar_d}^{\someEntity}$ the projection
of the mortar variable onto a discrete entity $\someEntity$ of the adjacent
discretizations, and by
$\projection_{\someVariable}^{\element}$ the projections of a variable $\someVariable$
defined on an adjacent subdomain onto the space of piecewise constants (per grid element)
on $\mesh_{\mortar_d}$. For example, $\projection_{\mortar_d}^\face$ is the projection
of the mortar variable onto the face $\face \in \faceSet_{d+1, \couplIdx}$, and
$\projection_{\head_d}^{\element}$ is the projection of the hydraulic head onto
the element $\element \in \mesh_{\mortar_d}$. Since we only consider projections into
spaces of piecewise constants, they can be constructed by computing the overlapping
regions between elements of the different discretizations and integrating over them.

\paragraph{Discrete flux expressions on coupling faces}
With the projections defined, the interface fluxes on coupling faces
$\face \in \faceSet_{d, \couplIdx}$, embedded in the control volume
$\cv \in \cvSet_d$,
are simply given by the projection of the corresponding mortar variable:
\begin{equation}
    \discFlux_\cv^\face = \projection_{\mortar_{d-1}}^{\face}.
    \label{eq:couplFluxMortar}
\end{equation}

\paragraph{Interface conditions}
Recall that, depending on the permeability contrast between the adjacent subdomains,
two types of interface conditions are considered, stated
in~\cref{eq:prob_mixed_if_flux,eq:prob_mixed_if_p}. Therefore, let us partition the
elements of the mortar discretization $\mesh_{\mortar_d}$ into the sets
$\mortarElemSet_{d, \condInterfaceSet}$ and $\mortarElemSet_{d, \blockInterfaceSet}$,
which collect the elements that live on interfaces on which the
conditions~\eqref{eq:prob_mixed_if_p} and \eqref{eq:prob_mixed_if_flux} should hold,
respectively. The interface conditions are then enforced by solving
\begin{subequations}
    \label{eq:ifConditionsMortar}
    \begin{align}
        \mortar_\element
        + \projection_{\meas{\crossSection_{d+1}}}^\element
          \projection_{\mortarHelperFlux_d}^\element
        &= 0,
        \quad \mathrm{on} \quad \element \in \mortarElemSet_{d, \blockInterfaceSet},
        \label{eq:ifConditionsMortarBlock} \\
        \projection_{\head_{d+1}}^\element
        - \projection_{\head_d}^\element
        &= 0,
        \quad \mathrm{on} \quad \element \in \mortarElemSet_{d, \condInterfaceSet},
        \label{eq:ifConditionsMortarCond}
    \end{align}
\end{subequations}
where we have introduced
\begin{equation}
    \mortarHelperFlux_d = \n^T
                          \perm_d
                          \n
                          \frac{
                              \head_{d} - \projection_{\head_{d+1}}^\element
                               }{
                                  \fracLength_{d}
                               }
    \label{eq:mortarHelperFluxVariable}
\end{equation}
for the sake of readability. Note that \cref{eq:ifConditionsMortarCond} is independent
of the mortar variable, which means that the associated rows of the linear system
contain zeros on the main diagonal, and only contain values in the coupling
blocks. Therefore, when using the conditions~\eqref{eq:prob_equi_if_p}, we have
experienced that the system matrices become singular when $\mortarGridFactor < 1$,
which seems to cause the system to be overconstrained. This is in agreement with
the analysis provided in~\citet{Wohlmuth2011Survey}, where Lagrange multiplier
approaches are investigated in the context of contact mechanics. They show that
for piecewise constant mortars, the system is unstable also for
$\mortarGridFactor = 1$, but one possible remedy is to use coarser meshes for
the mortar discretization. The case of $\mortarGridFactor = 1$ is briefly discussed
in the numerical examples given in~\cref{sec:examples}, but in general
we choose $\mortarGridFactor > 1$ whenever the conditions~\eqref{eq:ifConditionsMortarCond}
are used.

\section{Numerical experiments}
\label{sec:examples}

Ideally, a numerical method should be both efficient and accurate in order to
be useful, for instance, in the design process of a geotechnical engineering
application. Therefore, the test cases presented in the sequel aim at comparing
the above-mentioned schemes in terms of accuracy and efficiency.
We start with an investigation of the convergence behavior of the schemes
against an analytical solution in~\cref{sec:convAnalytic}.
Similar investigations have been the subject of
various works in the literature, where several authors
(see \eg~\cite{Sandve2012MpfaDfm,Ahmed2015MpfaDfm,Ahmed2017MpfaFpsDfm,Glaeser2017MpfaDfmTwoP})
used the analytical solution proposed in~\cite{haegland2009comparison}, although
originally derived for the equi-dimensional model problem. Consequently, it was
observed that convergence is lost when the size of the cells of the bulk
discretization approaches the aperture of the fracture, and an irreducible error
remains. This is expected as the mixed-dimensional formulation is an approximation
to the equi-dimensional problem. However, it is questionable if that irreducible
error provides a good measure of the modelling error introduced
by the dimension
reduction of the fracture. In \cite{Glaeser2017MpfaDfmTwoP} it was observed that
convergence can be maintained up to smaller cell sizes if the volume error is ruled
out, while \cite{glaeser_2020} report that the same can be achieved by integrating
the source term in the bulk elements adjacent to the fracture only over the region
that is not occupied by the fracture in the equi-dimensional setting. Such issues
seem to have little relevance to applications where the models could be applied to.
Moreover, the analytical solution is defined such that no jump in hydraulic head
occurs across the fracture, independent of the chosen permeability
contrast, which does not allow for a statement to be made whether or not a model
can be applied to low-permeable fractures. A mixed-dimensional analytical solution
is reported, for instance, in \eg~\citet{antonietti2016mimetic}, which overcomes
the issue with integration of the source term, but it also does not
feature a jump in the hydraulic head, which could appear in the case of low-permeable
fractures.\\

In this work, we want to propose an analytical solution to the mixed-dimensional
problem, in which the hydraulic head is discontinuous across the fracture for small
values of the fracture permeability, but approaches continuity for high fracture
permeabilities or small aperture values. This is subject of~\cref{sec:convAnalytic},
while in~\cref{sec:convDiscrete} we want to investigate the errors introduced by the
dimension reduction of the fracture using discrete reference solutions obtained on
equi-dimensional discretizations. Therein, we also study the performance of the
schemes for anisotropic permeability tensors in the bulk medium. In~\cref{sec:interface},
we lay special focus on the hydraulic head and fluxes at the bulk-fracture interface,
and compare the
\eboxDfm and the \eboxMortarDfm schemes when the conditions~\cref{eq:prob_mixed_if_p}
are used. Finally, we apply the schemes to
two- and three-dimensional benchmark cases in ~\cref{sec:benchmark2d}
and~\cref{sec:benchmark3d}, where, apart from the numerical results, we compare
the different schemes in terms of computational efficiency.

\subsection{Case 1: analytical reference solution}
\label{sec:convAnalytic}

Let us consider the unit square
$\domain = \domain_2 = \left(-0.5, 0.5 \right) \times \left(-0.5, 0.5 \right)$
to describe a porous medium that is intersected by a single horizontal and
lower-dimensional fracture
$\domain_1 = \left(-0.5, 0.5 \right) \times \{0\}$. The lower and upper parts
of $\domain_2$, that are separated by the fracture, are denoted by $\domain_2^1$
and $\domain_2^2$, respectively. The permeability of the bulk domain is
$\perm_2 = \I$ \si{\meter\per\second} and that of the fracture
$\perm_1 = k \, \I$ \si{\meter\per\second}. On this setting, let us define with
\begin{equation}
    \Delta H = \frac{\aperture}{k} \left( x + 0.5 \right)^2
    \label{eq:convHeadJump}
\end{equation}
a jump in hydraulic head across the fracture, for which we observe that
$\Delta H \rightarrow 0$ for $\aperture/k \rightarrow 0$ and
$x \in \left[-0.5, 0.5 \right]$.
The proposed analytical solution, denoted by $\head_d^\exact$, reads:
\begin{subequations}
  \label{eq:convTestSolution}
  \begin{align}
    \head_2^\exact \left(x, y\right) &=
            x^3 + y^3 + y\left( x + 0.5 \right)^2 - \frac{\Delta H}{2},
                &&\mathrm{for} \left(x, y\right) \in \domain_2^1,
                \label{eq:convTestH21} \\
    \head_2^\exact \left(x, y\right) &=
            x^3 + y^3 + y\left( x + 0.5 \right)^2 + \frac{\Delta H}{2},
                &&\mathrm{for} \left(x, y\right) \in \domain_2^2,
                \label{eq:convTestH22} \\
    \head_1^\exact \left(x, y\right) &= x^3,
                &&\mathrm{for} \left(x, y\right) \in \domain_1,
                \label{eq:convTestH1}
  \end{align}
\end{subequations}
with the source terms
\begin{subequations}
  \label{eq:convTestSource}
  \begin{align}
  \source_2 \left(x, y\right) &=
      -6 x -8 y + \frac{\aperture}{k}, &&\mathrm{for} \, \left(x,y\right) \in \domain_2^1, \\
  \source_2 \left(x, y\right) &=
      -6 x -8 y - \frac{\aperture}{k}, &&\mathrm{for} \, \left(x,y\right) \in \domain_2^2, \\
  \source_1 \left(x, y\right) &=
      -6 k x, &&\mathrm{for} \, \left(x,y\right) \in \domain_1.
  \end{align}
\end{subequations}
\begin{figure}[ht]
    \begin{subfigure}{0.499\textwidth}
        \centering
        \includegraphics[width=0.99\textwidth]{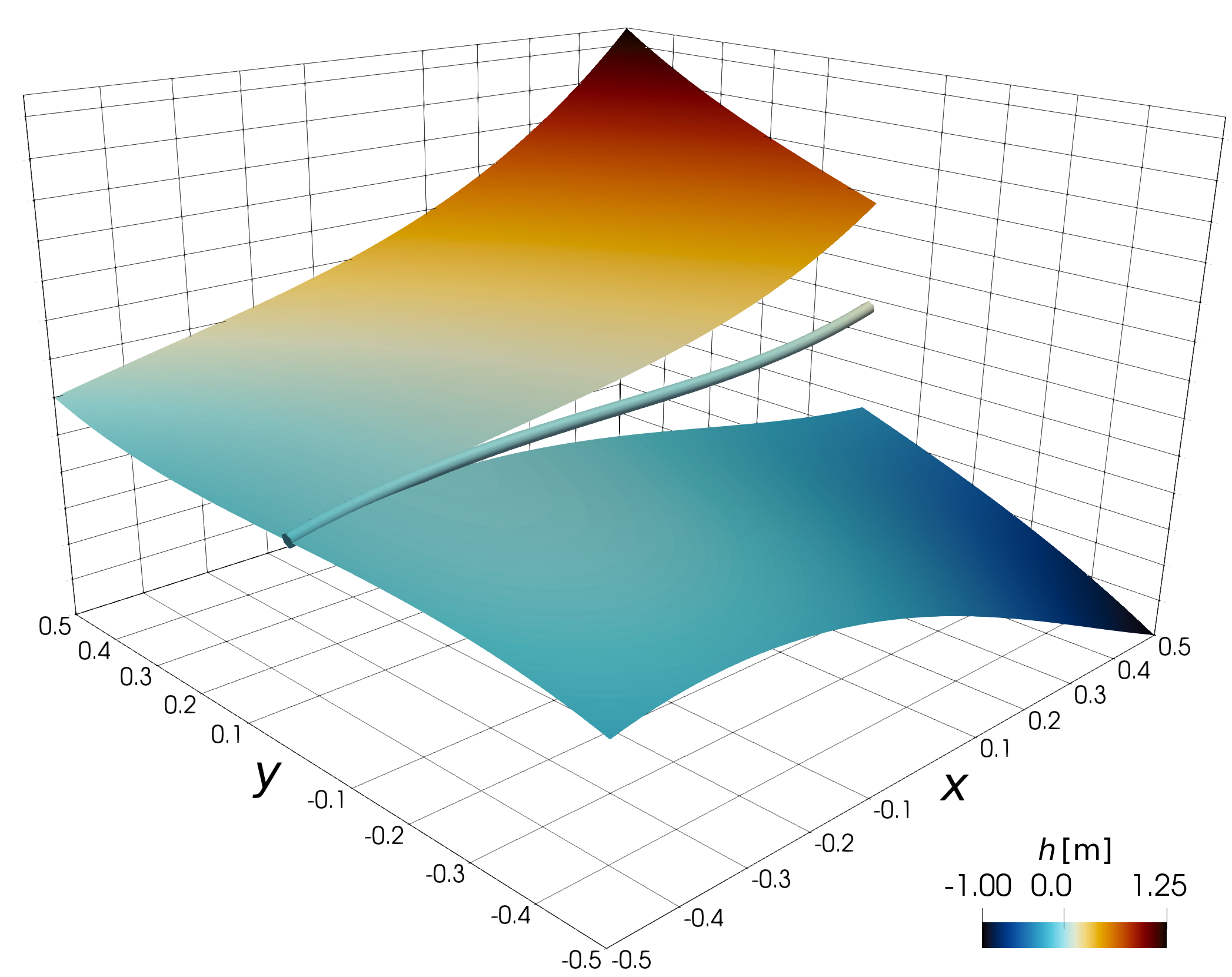}
        \caption{$k = \num{1e-4}$}
        \label{fig:analyticBarrierSol}
    \end{subfigure}
    \begin{subfigure}{0.499\textwidth}
        \centering
        \includegraphics[width=0.99\textwidth]{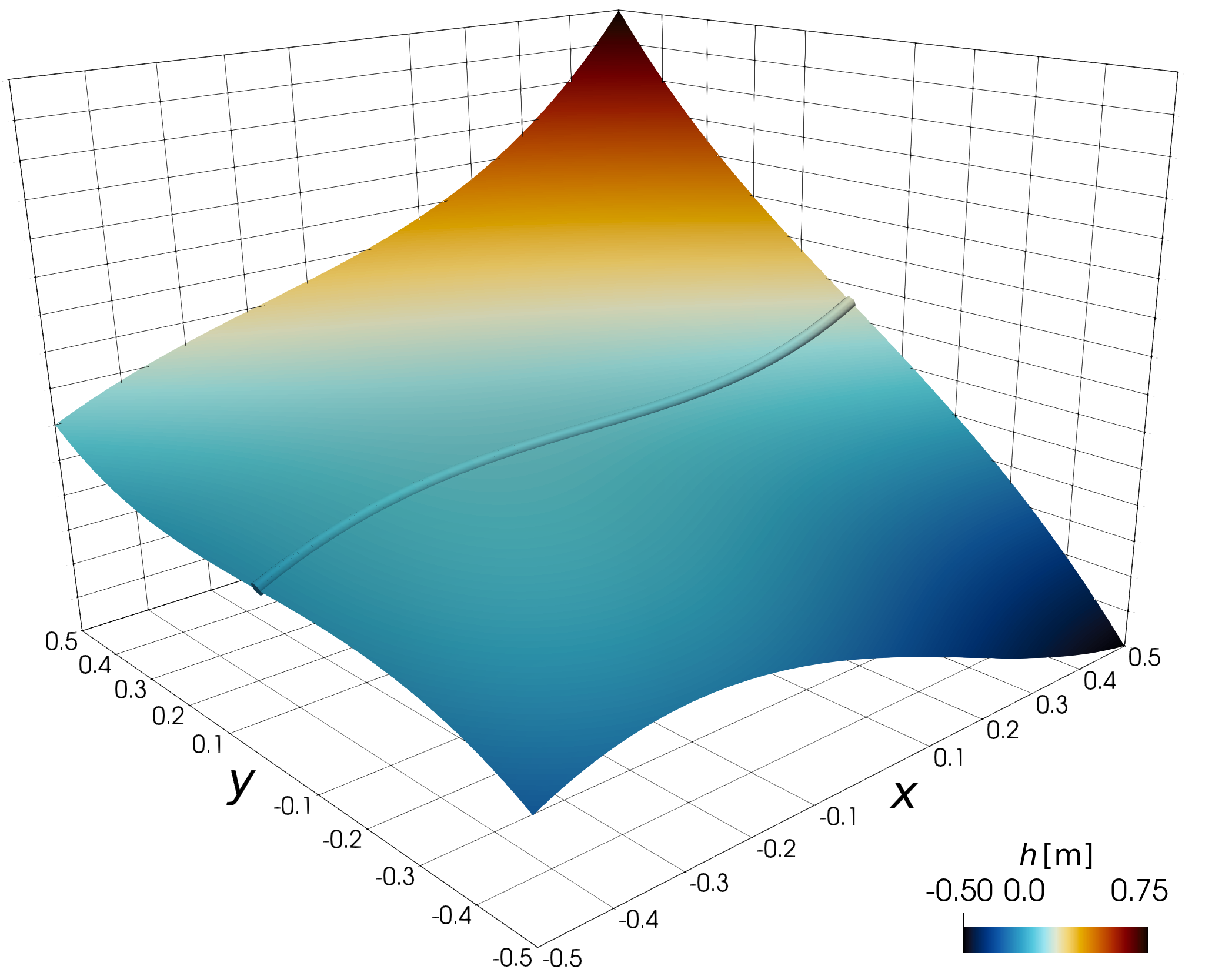}
        \caption{$k = \num{1e4}$}
        \label{fig:analyticConduitSol}
    \end{subfigure}
    \caption{\textbf{Case 1 - discrete solutions}. Depicted are the discrete solutions obtained with the \eboxDfm scheme on the finest grid of the convergence test. The solution in the fracture is illustrated by the tube at $y=0$.}
    \label{fig:convTestSols}
\end{figure}
\begin{figure}[ht]
    \begin{subfigure}{0.499\textwidth}
        \centering
        \includegraphics[width=0.99\textwidth]{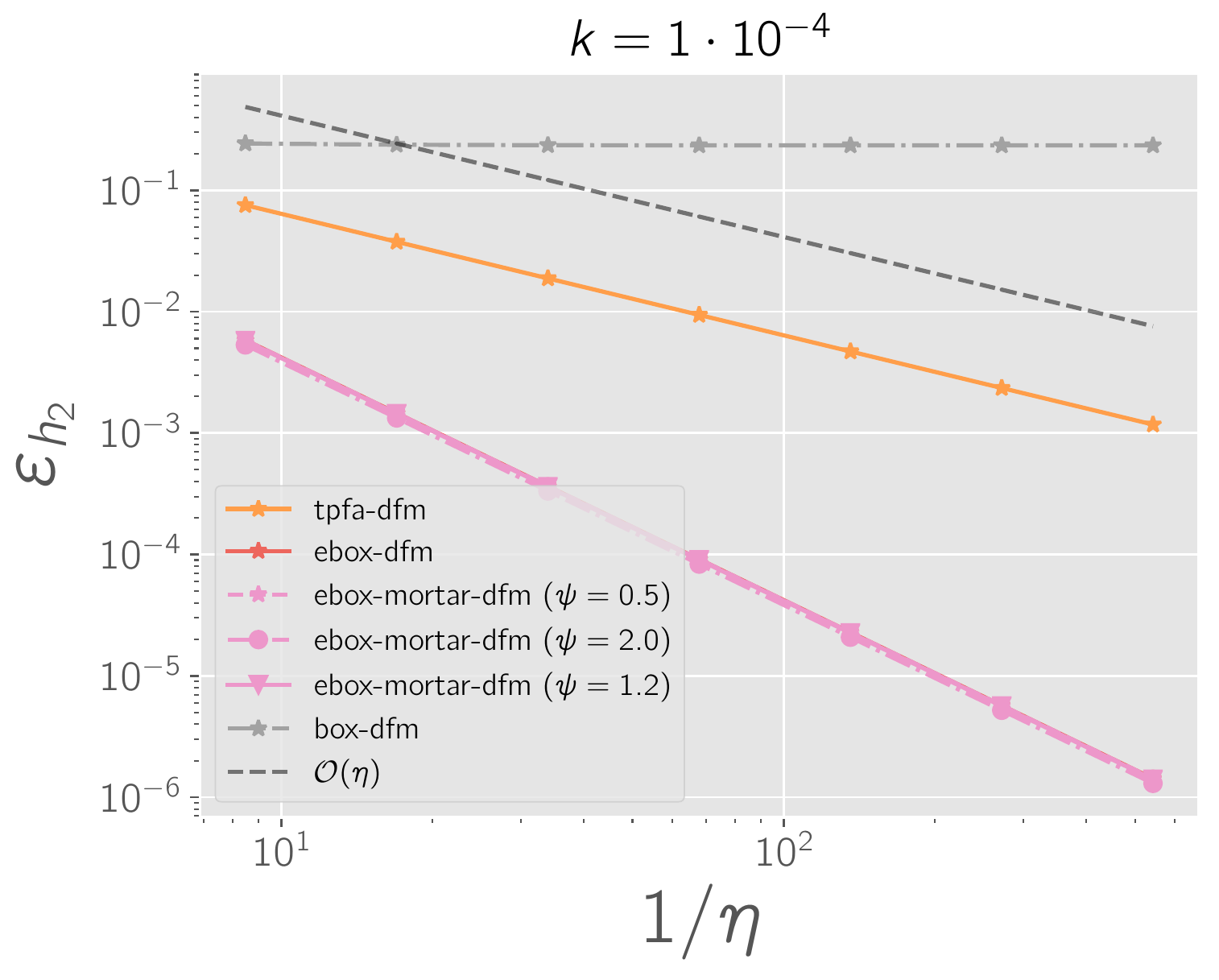}
    \end{subfigure}
    \begin{subfigure}{0.499\textwidth}
        \centering
        \includegraphics[width=0.99\textwidth]{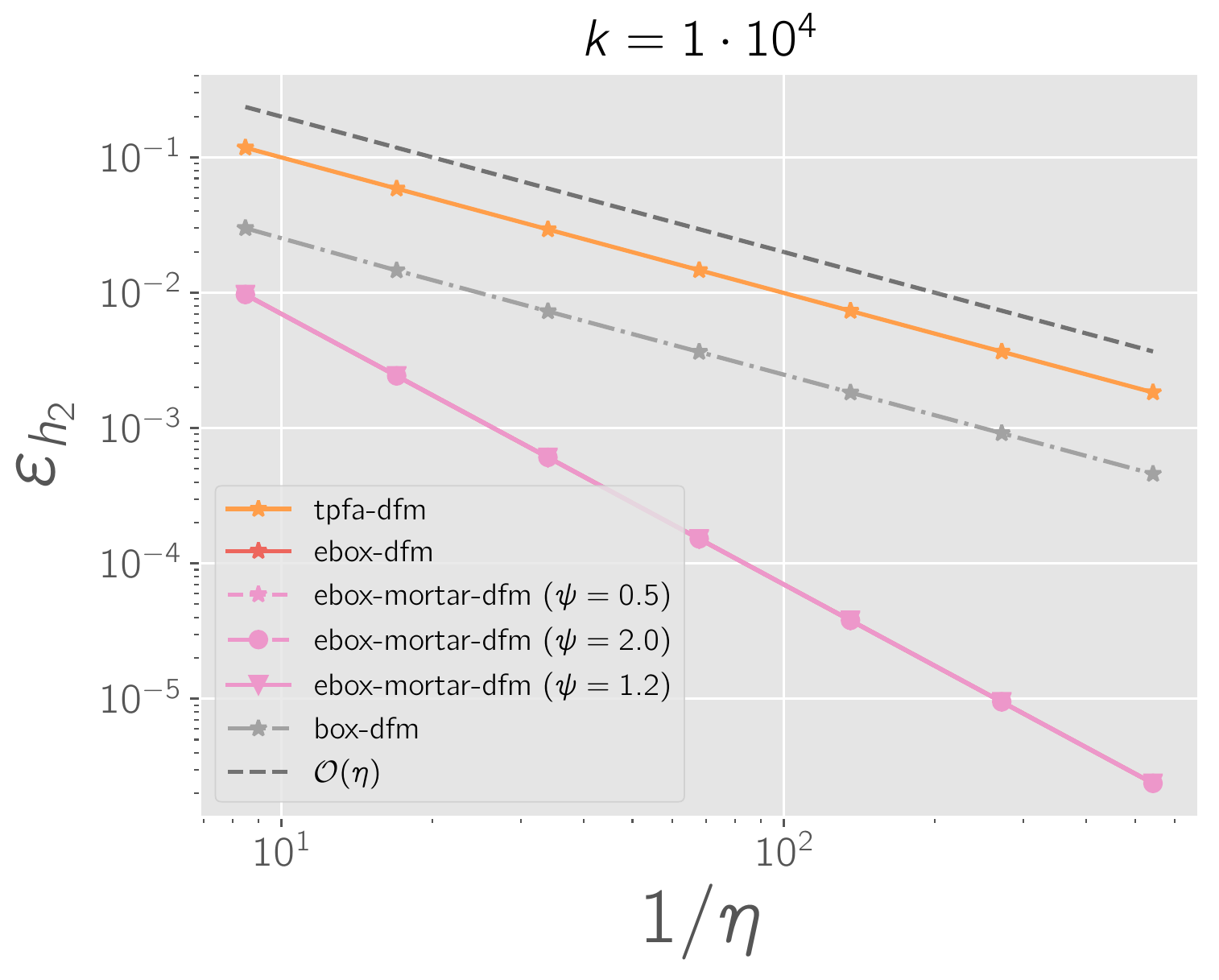}
    \end{subfigure}
    \begin{subfigure}{0.499\textwidth}
        \centering
        \includegraphics[width=0.99\textwidth]{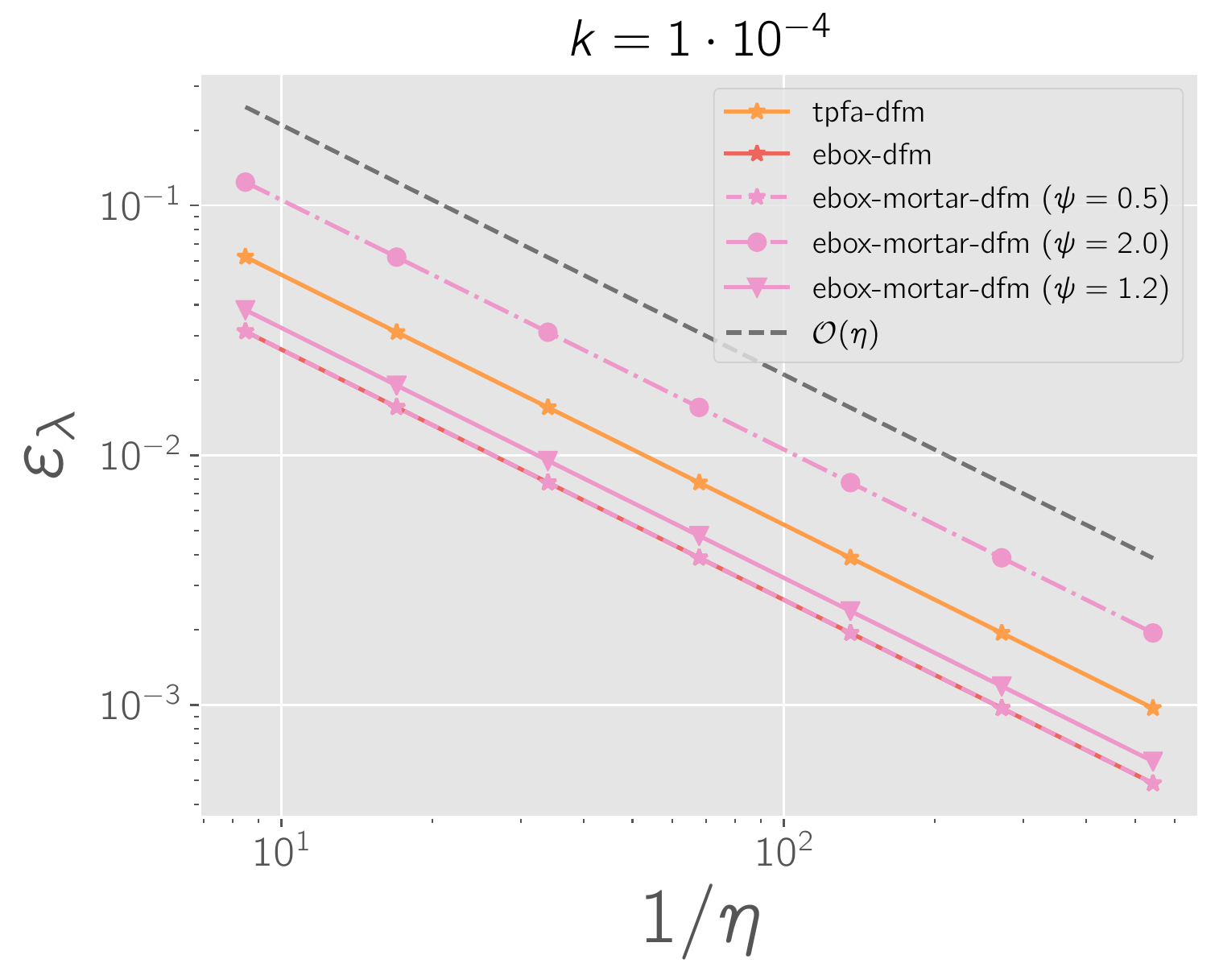}
    \end{subfigure}
    \begin{subfigure}{0.499\textwidth}
        \centering
        \includegraphics[width=0.99\textwidth]{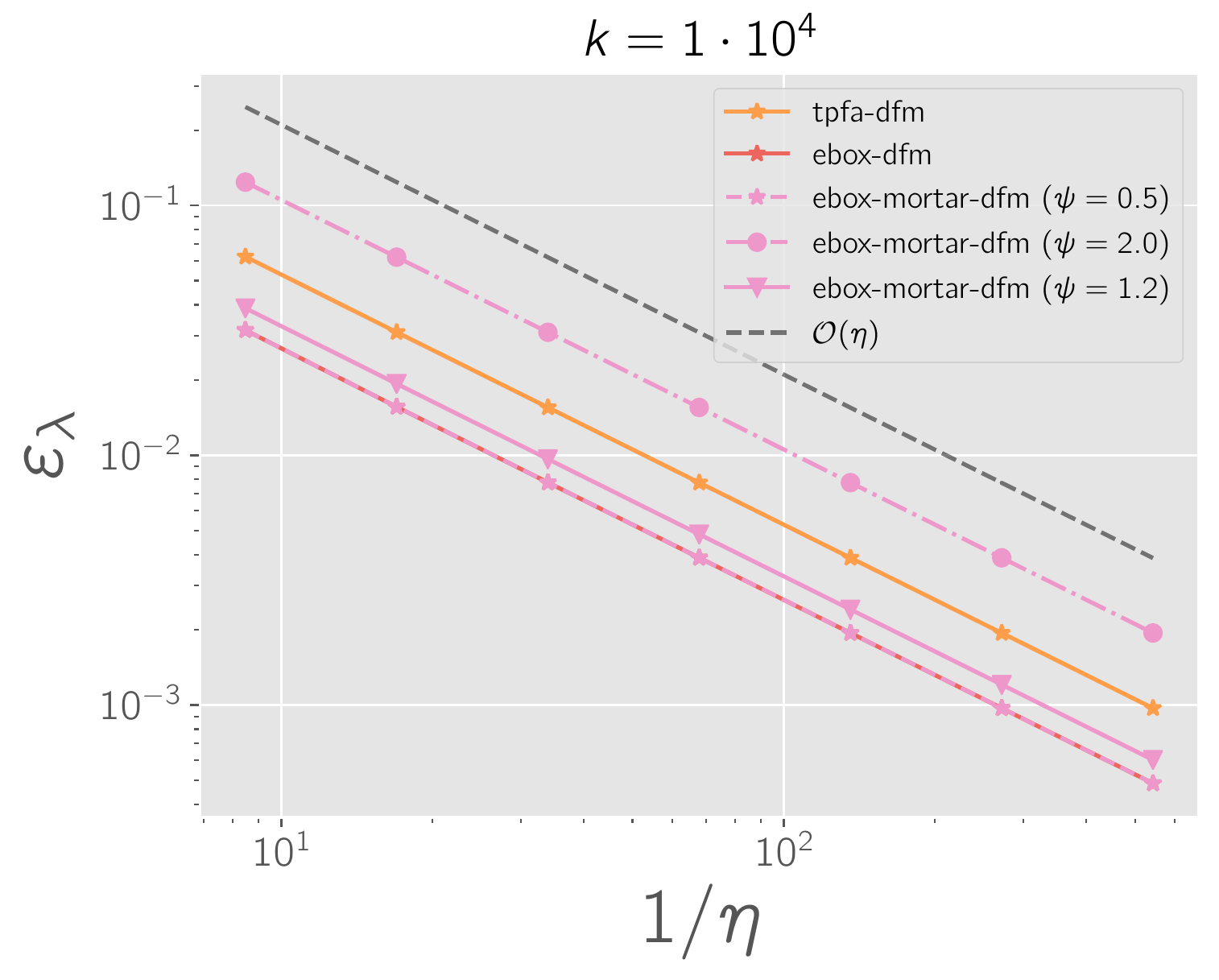}
    \end{subfigure}
    \caption{\textbf{Case 1 - errors over grid refinement}: the errors in $\head_2$ and the mass transfer $\lambda$ are plotted over grid refinement. The $x$-axis shows the inverse of the discretization length $\discLength$ and the left and right columns depict the results for $k = \num{1e-4}$ and $k = \num{1e4}$, respectively.}
    \label{fig:convAnalytic}
\end{figure}
In~\cref{fig:convTestSols}, the solution is visualized for
$\aperture = \SI{1e-4}{\meter}$
and $k \in \{\num{1e-4}, \num{1e4}\}$, which are the two cases that we consider
in the analysis. We use Dirichlet boundary conditions on all boundaries and set
$\faceSet_{d, \couplIdx} \equiv \faceSet_{d, \blockInterfaceSet}$,
that is, the conditions~\eqref{eq:prob_mixed_if_flux} are enforced on all coupling
faces. For the \eboxMortarDfm scheme we use $\mortarGridFactor = \{ 0.5, 1.2, 2.0 \}$,
that is, we consider the case that the mortar grid cells coincide with the coupling
sub-control volume faces of the bulk discretization, the case of a non-conforming
mortar grid, and a mortar grid whose elements contain exactly four bulk sub-control
volume faces.
For the quantification of the errors of the discrete solution
with respect to the analytical solution we use the norm
\begin{equation}
  \errorNorm_{\head_2}^i
  = \frac{
      \sqrt{
       \sum_{\element \in \mesh_{2}^i}
               \int_\element \left( \head - \head^\exact \right)^2
            }
    }{
    \sqrt{
        \sum_{\element \in \mesh_{2}^m}
            \int_\element \left( \head^\exact \right)^2
           }
    }
  \label{eq:errorNormConvTest1_H2}
\end{equation}
for the hydraulic head in the bulk medium, and for the
bulk-fracture transfer fluxes we use
\begin{equation}
  \errorNorm_{\lambda}^i
  = \frac{
        \sqrt{
           \sum_{\face \in \faceSet_{2, \couplIdx}^i}
                   \int_\face
                   \left( \lambda
                          - \lambda^{\exact}
                   \right)^2
               }
        }{
        \sqrt{
           \sum_{\face \in \faceSet_{2, \couplIdx}^m}
                   \int_\face
                   \left( \lambda^{\exact} \right)^2
               }
    },
  \label{eq:errorNormConvTest1_Q2}
\end{equation}
where we have introduced
$\lambda = \darcyVel_2 \scal \n |_{\interface_1}$
for the sake of readability. In the above norms, $i$ is the refinement
index, $0 < i \leq m$, and $\mesh_{2}^{i}$ denotes the corresponding mesh of the
bulk medium. The analysis is performed on structured rectangular grids, and
therefore, we do not consider the \mpfaDfm scheme in the discussion as it
reduces to the \tpfaDfm scheme on such settings, yielding the exact same results.\\

\Cref{fig:convAnalytic} shows the plots of the errors over grid refinement,
where second-order convergence is observed in $\head_2$ for the
vertex-centered \eboxDfm and \eboxMortarDfm schemes, both for a highly-permeable
and a low-permeable fracture. Furthermore, the errors and convergence behavior of
the \eboxMortarDfm scheme seems to be independent of the grid ratio
$\mortarGridFactor$. In particular, almost the same errors and rates are observed
for the case of non-matching mortar grids ($\mortarGridFactor = 1.2$).
In contrast to that, first-order convergence is observed for the \tpfaDfm scheme.
Note that second-order convergence would be seen if a norm based on point-wise
error computations in the grid element centers was used
\citep{Sandve2012MpfaDfm,Ahmed2015MpfaDfm,Glaeser2017MpfaDfmTwoP}, however, a
benefit of the vertex-centered schemes is that they allow for a straightforward
interpolation of the hydraulic head within the entire domain. Finally, first-order
convergence is observed for the \boxDfm scheme in the case of $k = \num{1e4}$, while
convergence is lost for the low-permeable fracture as a consequence of the
assumption of continuity of the hydraulic head.
\\

The interface flux $\lambda$ converges with first order for all schemes that
allow for its evaluation, and it can be seen that this holds again for all
considered values of $\mortarGridFactor$. As expected, higher errors in $\lambda$
are observed for coarser mortar discretizations, while the convergence rates
seem not to be affected by the choice of $\mortarGridFactor$. Please note that
results for the \boxDfm scheme are not shown in the plots for $\errorNorm_\lambda$,
as the scheme does not allow for an evaluation of the interface flux, which also does not
appear in its formulation.

\subsection{Case 2: discrete reference solution}
\label{sec:convDiscrete}

\begin{figure}[ht]
    \begin{subfigure}{0.3299\textwidth}
        \centering
        \includegraphics[width=0.99\textwidth]{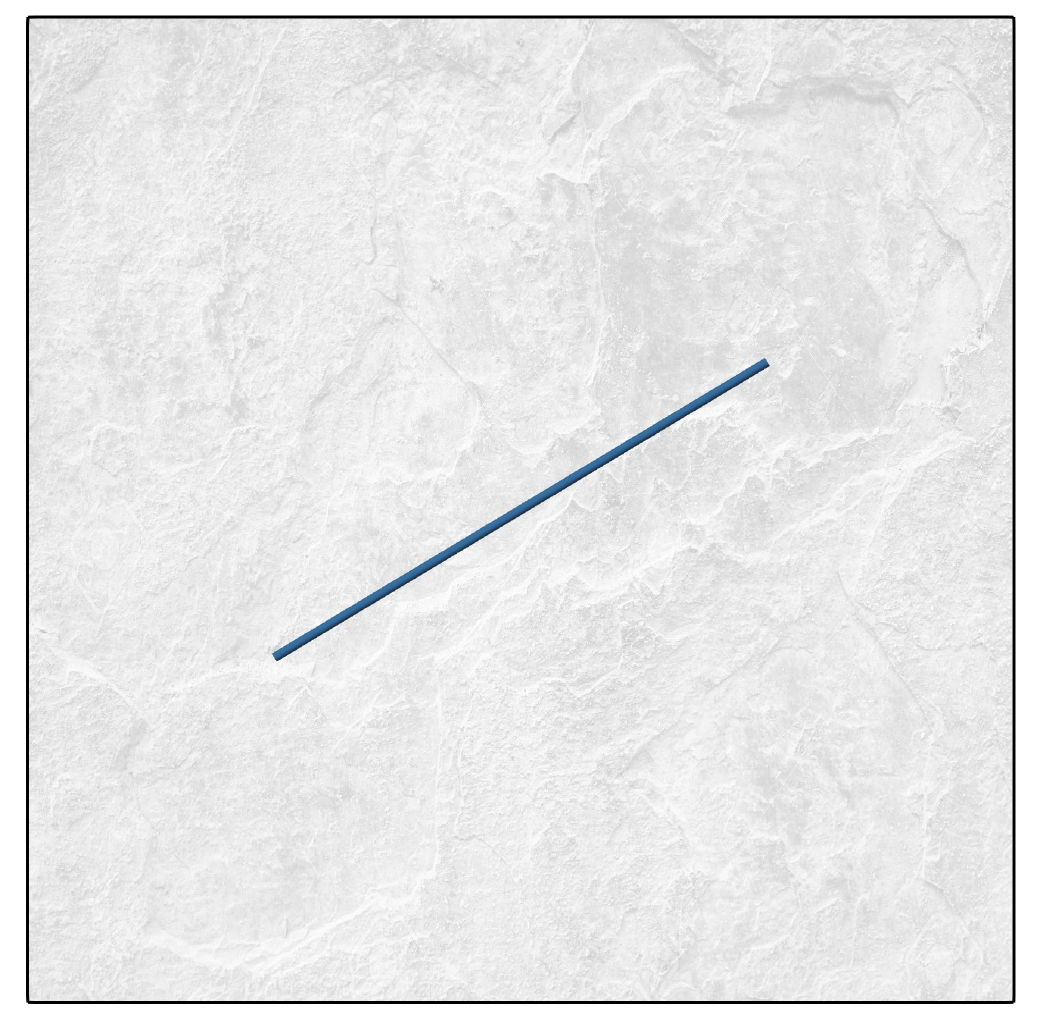}
        \caption{Case 2.1 / $\fracNet_1$}
        \label{fig:convDiscreteNetwork1}
    \end{subfigure}
    \begin{subfigure}{0.3299\textwidth}
        \centering
        \includegraphics[width=0.99\textwidth]{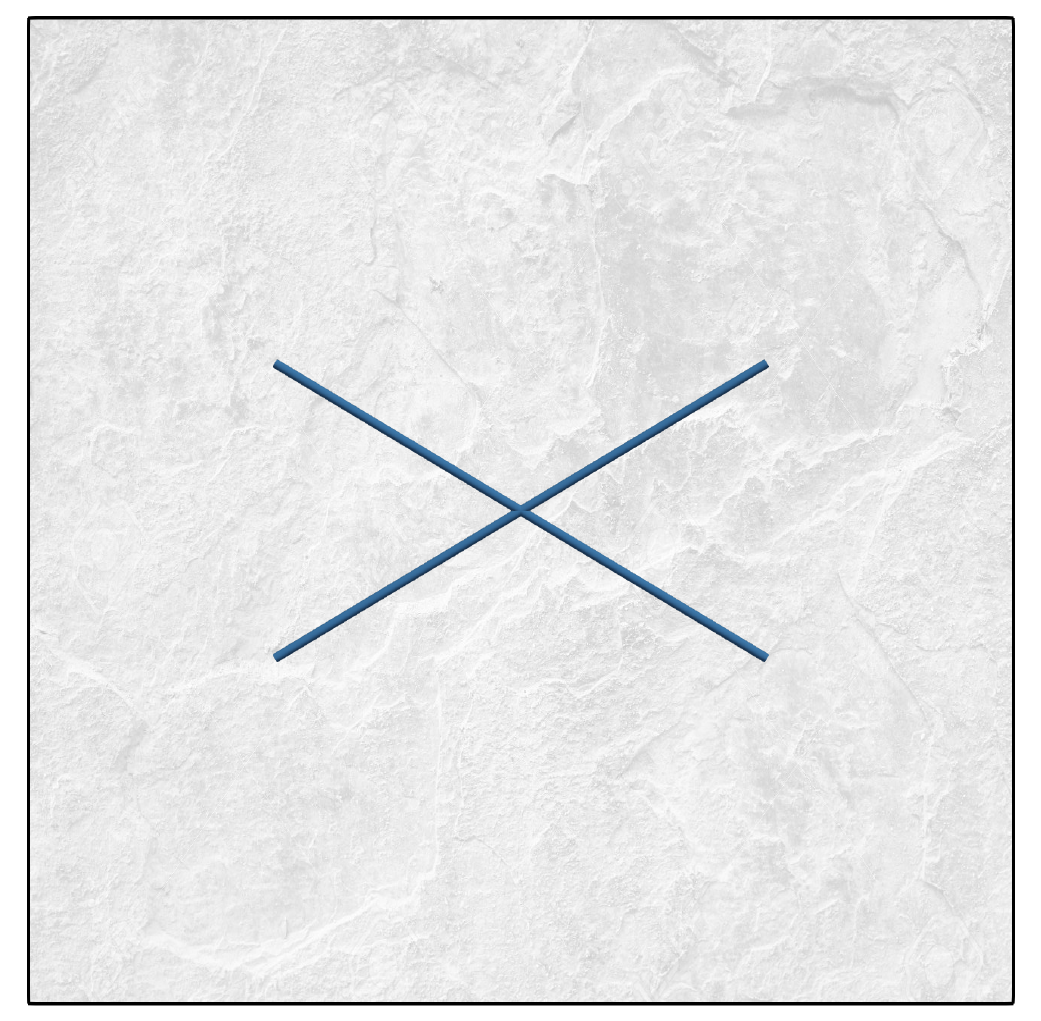}
        \caption{Case 2.2 / $\fracNet_2$}
        \label{fig:convDiscreteNetwork2}
    \end{subfigure}
    \begin{subfigure}{0.3299\textwidth}
        \centering
        \includegraphics[width=0.99\textwidth]{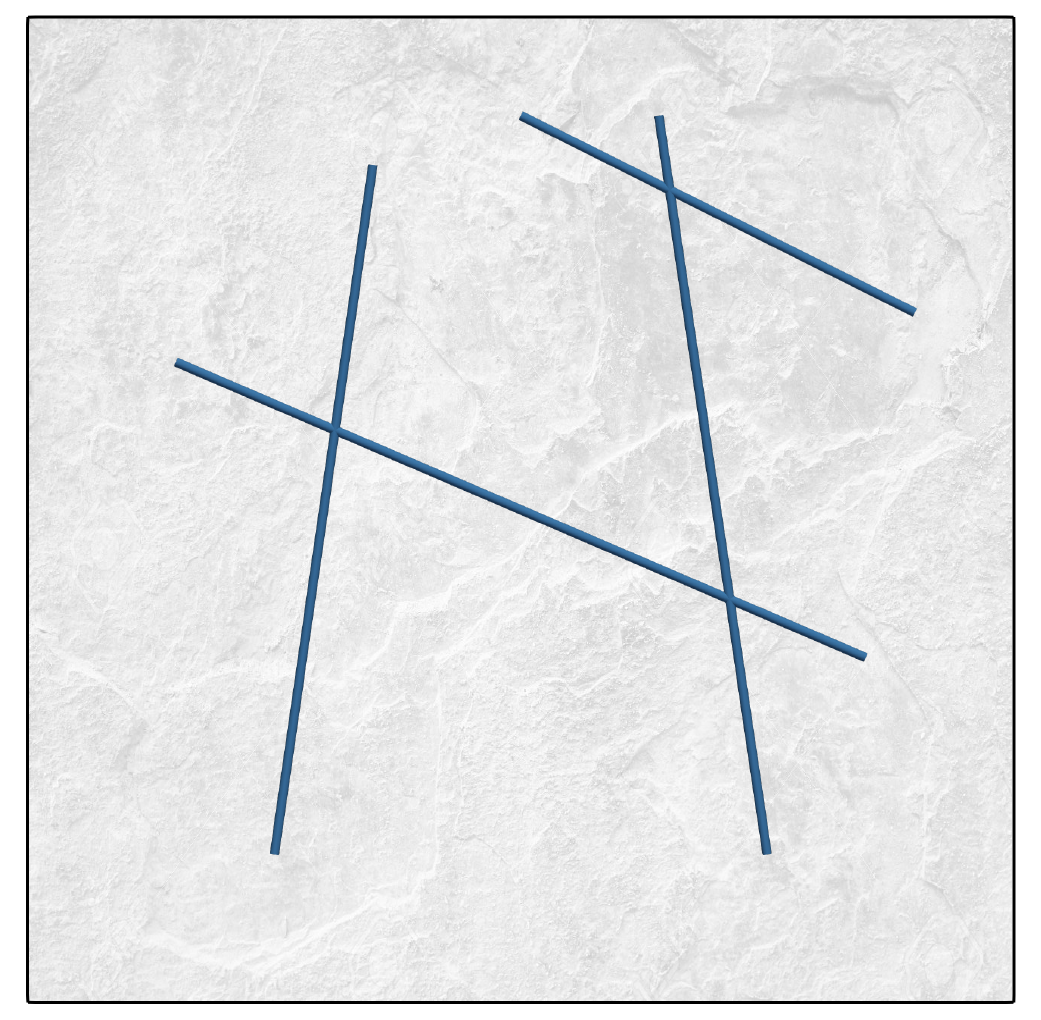}
        \caption{Case 2.3 / $\fracNet_3$}
        \label{fig:convDiscreteNetwork3}
    \end{subfigure}
    \caption{\textbf{Case 2 - fracture networks}: the figure depicts the three different fracture configurations, illustrated by the blue lines, that are investigated in the convergence tests against equi-dimensional reference solutions.}
    \label{fig:convDiscreteNetworks}
\end{figure}
As in the previous test case, we consider the domain
$\domain_2 = \left(-0.5, 0.5 \right) \times \left(-0.5, 0.5 \right)$
for the bulk medium, and take into account three different fracture
configurations $\fracNet_k$, $k \in \{1, 2, 3 \}$, composed of individual
fractures $\fracture$ such that $\fracNet_k =\bigcup_{\fracture \in \fracNet_k} f$.
An illustration of the different networks can be found in~\cref{fig:convDiscreteNetworks}.
The permeabilities in the subdomains are given by (in \si{\meter\per\second})
\begin{equation}
    \perm_\pmIdx = \perm_2 = \mathbf{R} \left( \permAngle \right)^{-1} \left( \begin{matrix} 1 & 0 \\ 0 & 1/5 \end{matrix} \right) \mathbf{R}\left( \permAngle \right), \quad
    \perm_\fracIdx = \perm_1 = k \, \I,
\end{equation}
where $\mathbf{R}$ is the two-dimensional rotation matrix in clockwise direction
around the angle $\permAngle$, and we consider both highly-permeable as well as
blocking fractures by choosing $k$ accordingly. 
An aperture of $\aperture = \SI{2e-3}{\meter}$ is used for all fractures and
Dirichlet boundary conditions are applied on the entire
outer boundary, evaluated after
\begin{equation}
    \head_2^\diriIdx = \left(1 - \left(x-0.5\right)^2 \right)
                       \left(1 - \Cos{\pi y} \right),
        \quad \mathrm{for} \, \left(x, y \right) \in \surface_2.
    \label{eq:convtest2DirichletBcs}
\end{equation}

\paragraph{Discrete error norms}
For each of the three fracture networks, we want to study the convergence behaviour
of the different mixed-dimensional schemes against equi-dimensional reference
solutions, which we obtain with the \mpfa scheme, for different combinations
of $k$ and $\permAngle$.
To this end, we create equi-dimensional discretizations of the domain geometries by
describing each fracture as a rectangle with a thickness equal to the aperture $\aperture$.
The equi-dimensional fracture domain $\domain_\fracIdx$ is then obtained by
computing the union of all fractures.
Let us denote the equi-dimensional discretizations by $\mesh^\eq$ and the corresponding
numerical solution by $\head^\eq$. We then extract the sub-meshes of $\mesh^\eq$
that describe the bulk medium and the fractures, which we denote here by
$\mesh^\eq_\pmIdx$ and $\mesh^\eq_\fracIdx$, respectively.
Let $\mesh_2^i$, $0 \leq i \leq 5$, be the $i$-th refinement of the discretization
of the bulk domain in the mixed-dimensional setting with the associated numerical
solution $\head_2^i$. The error of $\head_2^i$ with respect to the reference solution
is evaluated after
\begin{equation}
    \errorNorm_{\head_2}^i
    = \frac{\sqrt{
                \sum_{\element \in \mesh^\eq_\pmIdx}
                    \left( \projection_{\head_2^i}^{\element}
                           - \head_\element^\eq \right)^2 \meas{\element}
                }
           }{\sqrt{
                   \sum_{\element \in \mesh^\eq_\pmIdx}
                         \left( \head_\element^\eq \right)^2 \meas{\element}
           }},
    \label{eq:convtest2ErrorH2}
\end{equation}
where $\projection_{\head_2^i}^{\element}$ is the L2-projection of $\head_2^i$ onto
the element $\element \in \mesh^\eq_\pmIdx$. In order to evaluate the error in
$\head_1^i$, we generate one-dimensional reference solutions $\head_\fracture^\exact$
by averaging and projecting the equi-dimensional solution
onto fine discretizations of the individual fractures. These are denoted by
$\mesh_{\fracture}$, $\fracture \in \fracNet_k$, and allow for the computation
of the error in $\head_1$ for the $i$-th refinement and the $k$-th fracture
network after the following expression:
\begin{equation}
    \errorNorm_{\head_1}^i
    = \frac{\sum_{\fracture \in \fracNet_k}
            \sqrt{
                \sum_{\element \in \mesh_{\fracture}}
                        \left( \projection_{\head_{1}^i}^{\element}
                               - \head^\exact_{\element}
                        \right)^2 \meas{\element}
                 }
           }{
          \sum_{\fracture \in \fracNet_k}
          \sqrt{
               \sum_{\element \in \mesh_{\fracture}}
               \left( \head^\exact_{\element}
               \right)^2 \meas{\element}
               }
      }.
    \label{eq:convtest2ErrorH1}
\end{equation}
\paragraph{Discretization choices}
We use triangular grids in all simulations, and in order to reduce the cost for
the computation of the reference solution, we use
locally refined meshes towards the fractures, where we choose a characteristic discretization
length $\discLength^\eq_\fracIdx = \aperture/6$ in and around the fractures, and
use $\discLength^\eq_\surface = 5 \aperture$ at the exterior domain boundary.
For the mixed-dimensional discretizations, we generate initial meshes using
$\discLength_\fracIdx = 75\discLength^\eq_\fracIdx$, and
$\discLength_\surface = \discLength_\fracIdx \discLength_\surface^\eq / \discLength_\fracIdx^\eq$
in order to have the same coarsening ratio towards the boundaries as
in the reference grid.
All subsequent refinements are obtained by subdivision of the grid elements.
This strategy results in meshes as listed in~\cref{tab:convtest2Meshes},
where the number of elements in the reference grid as well as the coarsest and
finest mixed-dimensional grids are given. As can be seen, the ratio
$\meas{\mesh^\eq}/\meas{\mesh_2^5}$ is similar for all fracture networks,
ranging from $2$ to $2.7$. For the discretization of the mortar domain we
use $\mortarGridFactor = 1.2$, and as in the previous convergence test, we set
$\faceSet_{d, \couplIdx} \equiv \faceSet_{d, \blockInterfaceSet}$ in all simulations.
Please note that while on the structured grids used in the previous test case
the ratio $\mortarGridFactor$ between the discretization lengths used in the
mortar and fracture domains was identical along the entire fracture, this is
no longer the case for the unstructured grids used in this example. Here,
$\mortarGridFactor$ refers to the ratio of the characteristic lengths used
to generate the meshes with the mesh generation tool Gmsh \citep{Gmsh2009},
but that does not mean that
this ratio is fulfilled for all elements of the grids. The same also holds for
the grids used in~\cref{sec:interface,sec:benchmark2d}.\\

\begin{table}[t]
  \centering
  \caption{\textbf{Case 2 - meshes used in the convergence tests}. Discretization lengths used for the equi-dimensional reference grid, expressed as multiples of the fracture aperture $\aperture$. Additionally, the resulting number of elements are given for both the reference and the coarsest and finest mixed-dimensional meshes.}
  \begin{tabular}{ *{8}{l} }
  \toprule
                 & $\meas{\mesh^\eq}$ & $\meas{\mesh_2^0}$ & $\meas{\mesh_1^0}$ & $\meas{\mesh_2^5}$ & $\meas{\mesh_1^5}$\\ \midrule
    $\fracNet_1$ &  \num{351768}      & \num{172}          & \num{24}           & \num{176128}       & \num{768}         \\
    $\fracNet_2$ &  \num{1089484}     & \num{420}          & \num{48}           & \num{430080}       & \num{1536}        \\
    $\fracNet_3$ &  \num{3839860}     & \num{1386}         & \num{112}          & \num{1419264}      & \num{3584}        \\
  \bottomrule
  \end{tabular}
  \label{tab:convtest2Meshes}
\end{table}
The convergence tests are carried out for $k \in \{\num{1e-4}, \num{1e4} \}$ and
$\permAngle \in \{0, \pi/4 \}$, and all errors and rates are listed in the appendix
in~\cref{tab:convDiscreteErrors_head2_angle0,tab:convDiscreteErrors_head1_angle0,tab:convDiscreteErrors_head2_angle45,tab:convDiscreteErrors_head2_angle0},
of which we will discuss a selection in the subsequent paragraphs. In addition,
we present results obtained from further sets of simulations targeted at
investigating the dependency of the errors on a wider range of $k$ and
$\permAngle$. While using the same reference grid, we use
mixed-dimensional discretizations that are not obtained by refinement,
in contrast to the convergence tests, but instead, we generate triangulations by choosing
$\discLength_\fracIdx$ as in $\mesh_1^5$. \Cref{tab:convtest2MeshesAdditional}
lists the meshes used in those simulations, and it can be seen that it is
$\meas{\mesh^\eq}/\meas{\mesh_2} \approx 5.5$ for all fracture networks.
\begin{table}[h]
  \centering
  \caption{\textbf{Case 2 - meshes used in the additional simulations}. Number of elements in the mixed-dimensional discretizations used to investigate the dependency of the errors on $k$ and $\permAngle$.}
  \begin{tabular}{ *{4}{l} }
  \toprule
                     &  $\fracNet_1$  & $\fracNet_2$  & $\fracNet_3$ \\ \midrule
    $\meas{\mesh_2}$ &  \num{62254}   & \num{196476}  & \num{699682} \\
    $\meas{\mesh_1}$ &  \num{747}     & \num{1496}    & \num{3428}   \\
  \bottomrule
  \end{tabular}
  \label{tab:convtest2MeshesAdditional}
\end{table}

\begin{figure}[tbh]
    \begin{subfigure}{0.3299\textwidth}
        \centering
        \includegraphics[width=0.99\textwidth]{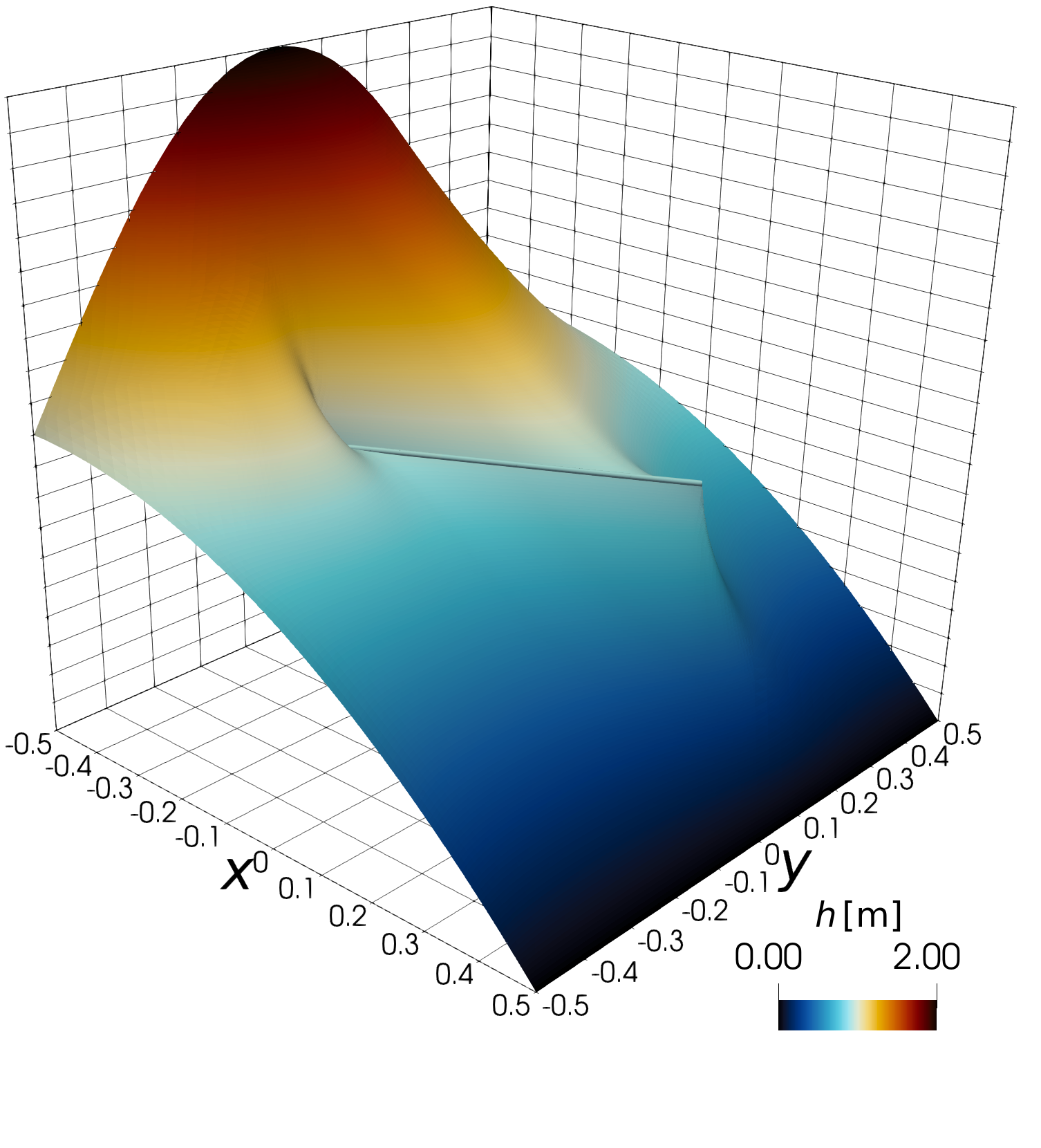}
        \caption{}
        \label{fig:convDiscreteN1ConduitSol}
    \end{subfigure}
    \begin{subfigure}{0.3299\textwidth}
        \centering
        \includegraphics[width=0.99\textwidth]{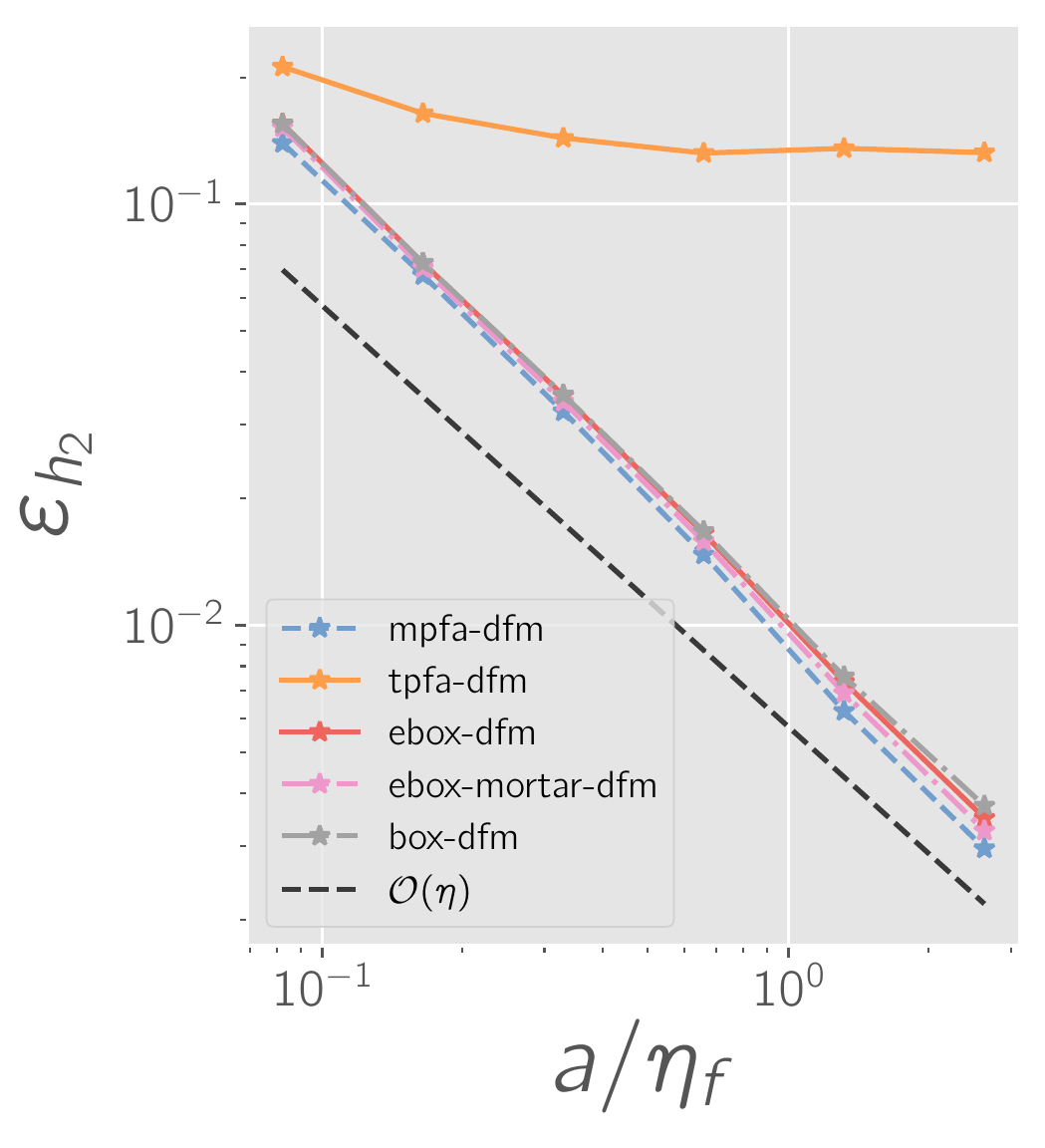}
        \caption{}
        \label{fig:convDiscreteN1ConduitBulk}
    \end{subfigure}
    \begin{subfigure}{0.3299\textwidth}
        \centering
        \includegraphics[width=0.99\textwidth]{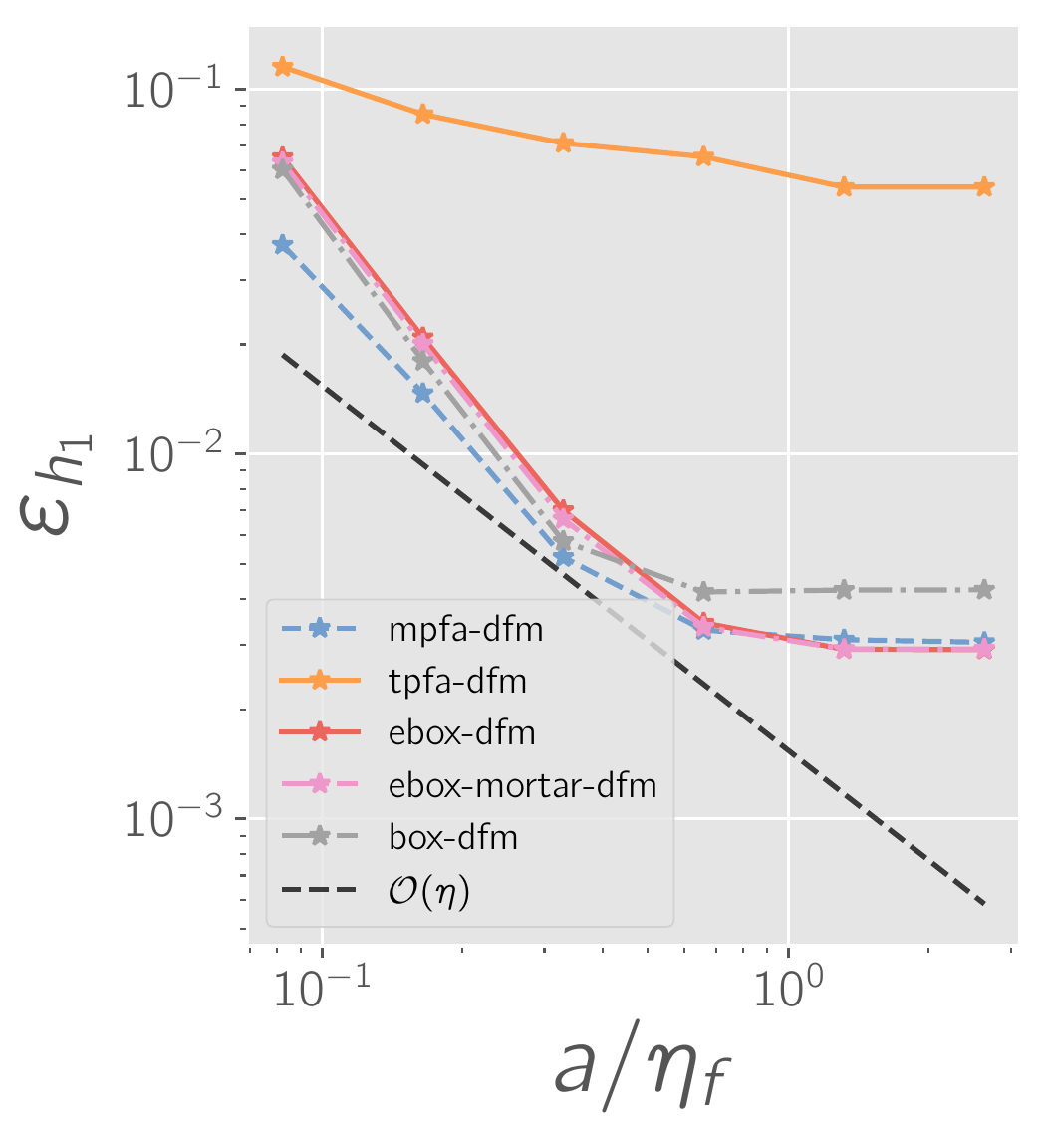}
        \caption{}
        \label{fig:convDiscreteN1ConduitFacet}
    \end{subfigure}
    \begin{subfigure}{0.3299\textwidth}
        \centering
        \includegraphics[width=0.99\textwidth]{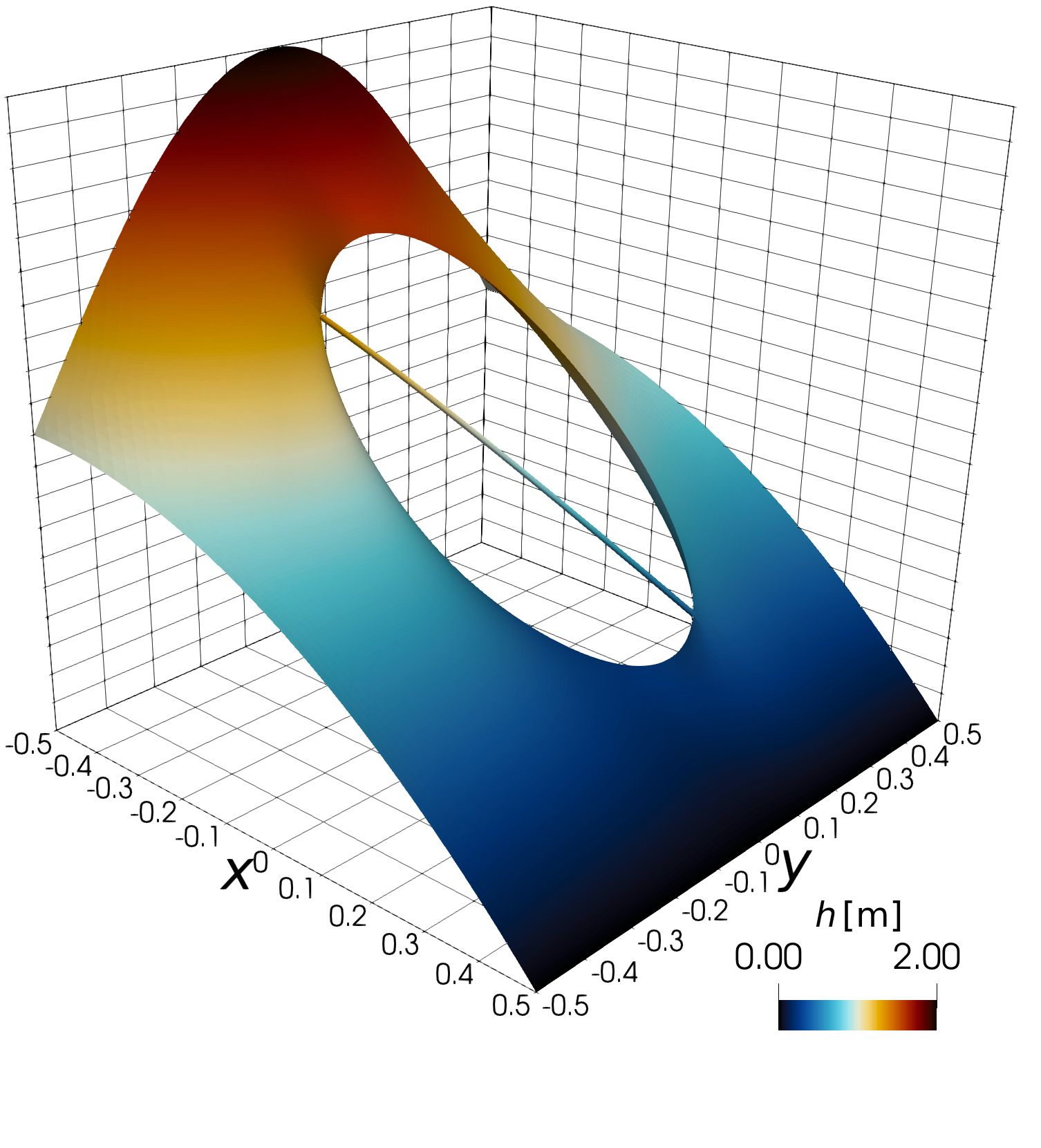}
        \caption{}
        \label{fig:convDiscreteN1BarrierSol}
    \end{subfigure}
    \begin{subfigure}{0.3299\textwidth}
        \centering
        \includegraphics[width=0.99\textwidth]{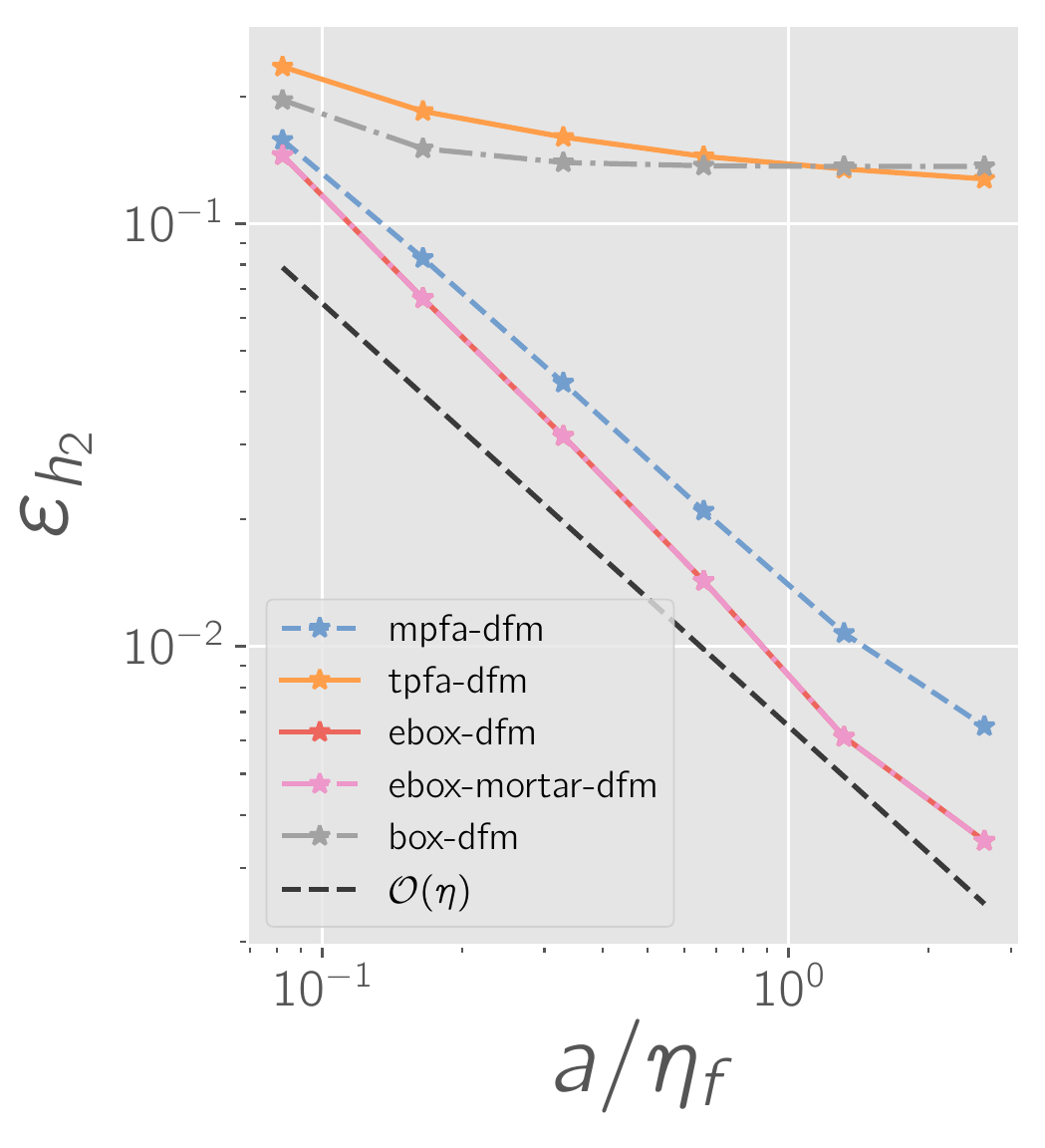}
        \caption{}
        \label{fig:convDiscreteN1BarrierBulk}
    \end{subfigure}
    \begin{subfigure}{0.3299\textwidth}
        \centering
        \includegraphics[width=0.99\textwidth]{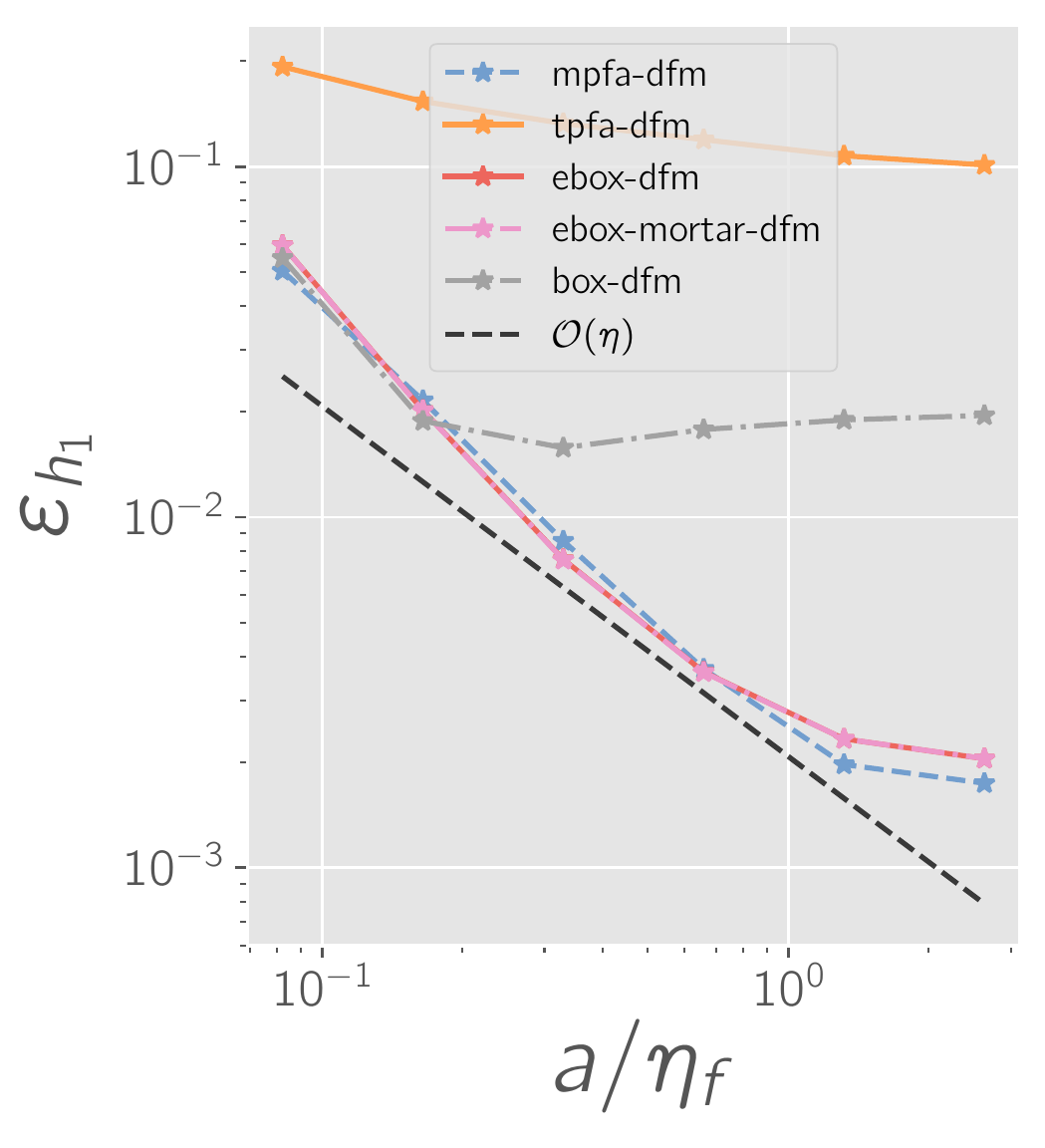}
        \caption{}
        \label{fig:convDiscreteN1BarrierFacet}
    \end{subfigure}
    \caption{\textbf{Case 2.1 - results for $\permAngle = 0$}. The upper row and lower row show the results for $k = \num{1e4}$ and $k = \num{1e-4}$, respectively. The left column depicts the distributions of the hydraulic heads in the domain, obtained with the \eboxDfm scheme on the finest grid,
    while the middle and right columns show the errors in $\head_2$ and $\head_1$ plotted over grid refinement.}
    \label{fig:convDiscreteN1}
\end{figure}
\begin{figure}[tbh]
  \begin{subfigure}{\textwidth}
    \centering
    \includegraphics[width=0.99\textwidth]{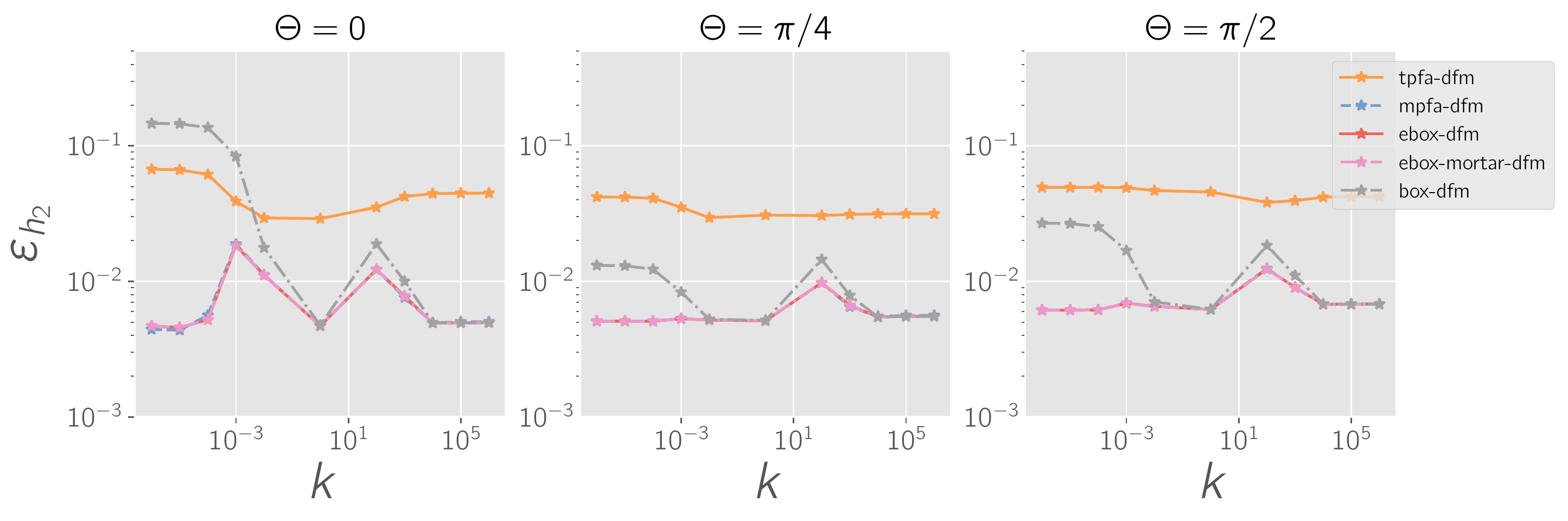}
  \end{subfigure}
  \begin{subfigure}{\textwidth}
    \centering
    \includegraphics[width=0.99\textwidth]{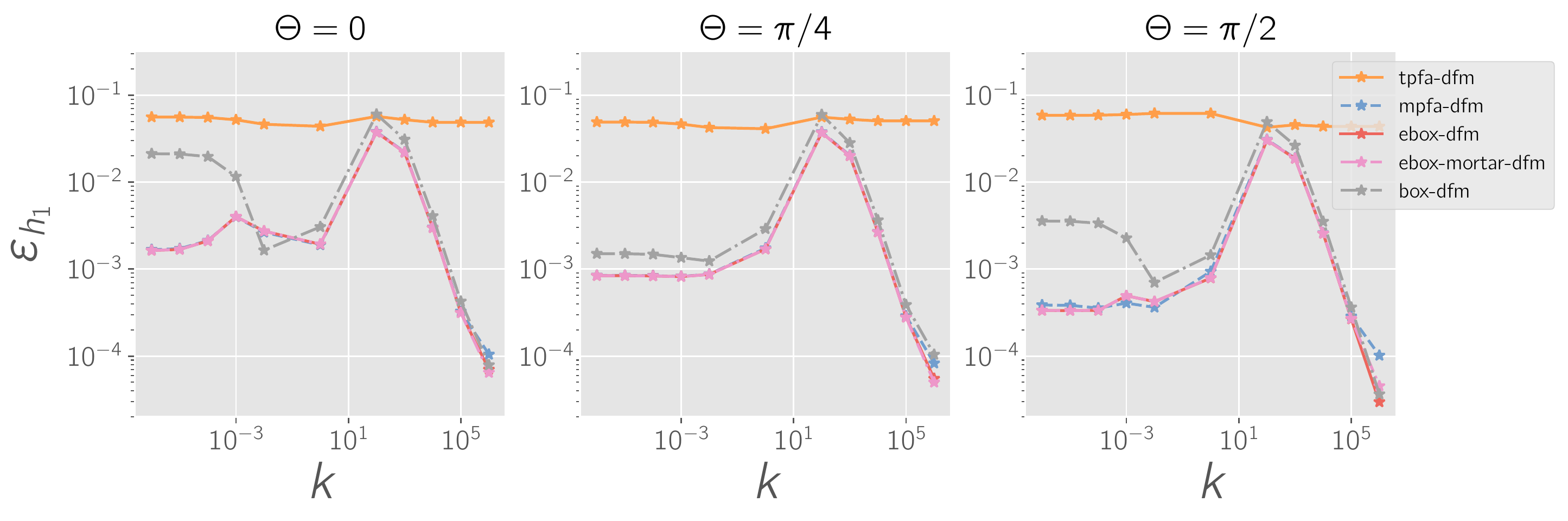}
  \end{subfigure}
  \caption{\textbf{Case 2.1 - error over $k$}. The errors in $\head_2$ and $\head_1$, obtained on the meshes as specified in~\cref{tab:convtest2MeshesAdditional}, are plotted over $k$ for various values of the permeability angle $\permAngle$.}
  \label{fig:convDiscreteN1_kplot}
\end{figure}
\paragraph{Fracture network $\fracNet_1$}
As illustrated in~\cref{fig:convDiscreteNetwork1}, this setting consists of a
single fracture placed in the center of the domain. The distributions of the
hydraulic head for $\permAngle = 0$ and both fracture permeabilities are
depicted in~\cref{fig:convDiscreteN1},
which also shows the relative errors plotted over grid refinement.
One expected observation we make, given that unstructured meshes and anisotropic
permeabilities are used, is that the \tpfaDfm scheme does not converge
in any of the considered quantities and for none of the considered parameter choices.
The error only slightly decays upon refinement, and relative errors in the order
of \SI{10}{\percent} remain on the finest grid. In contrast to that, we observe
that $\head_2$ converges with approximately first order for all remaining schemes
in the case of $k = \num{1e4}$, while the rates slightly decay upon the last refinement
for $k = \num{1e-4}$ (see~\cref{fig:convDiscreteN1BarrierBulk}).
However, the errors drop below \SI{1}{\percent} after the third refinement in
both cases. It is also noteworthy that for $k = \num{1e4}$, the \mpfaDfm,
\eboxDfm and \eboxMortarDfm schemes lead to very similar errors in each refinement,
while for $k = \num{1e-4}$, we see that the vertex-centered schemes produce
slightly smaller errors. Apart from that, only very minor differences in errors
are observed between the \eboxDfm and the \eboxMortarDfm schemes, despite the
more than two times coarser discretization of the interface fluxes in the latter.
Note that second-order convergence cannot be expected in this case due to the
projections into the space of piecewise constants on the reference grid.\\

The plots of $\errorNorm_{\head_1}$, given in~\cref{fig:convDiscreteN1ConduitFacet,fig:convDiscreteN1BarrierFacet},
show that the convergence rates deteriorate when the discretization length around the
fracture $\discLength_\fracIdx$ approaches the aperture $\aperture$.
As mentioned in the introduction to this section, this is expected since the
mixed-dimensional model is an approximation to the equi-dimensional problem.
However, rather small relative errors below \SI{0.5}{\percent} are observed
on the finest grid. Qualitatively similar results were also observed for
$\permAngle = \pi/4$, which are given in~\cref{tab:convDiscreteErrors_head2_angle45,tab:convDiscreteErrors_head1_angle45},
with the main difference being that smaller errors are obtained with the \boxDfm
scheme for $k = \num{1e-4}$. This stems from the fact that for $\permAngle = \pi/4$
the jump in hydraulic head across the fracture is significantly smaller for the
chosen fracture orientation and boundary conditions.\\

\Cref{fig:convDiscreteN1_kplot} depicts the results obtained
from a set of additional simulations carried out on the mesh as specified
in~\cref{tab:convtest2MeshesAdditional}. Again, a very similar performance of
the \mpfaDfm, \eboxDfm and \eboxMortarDfm schemes is observed, for which the error
in $\head_2$ lies between \SI{0.5}{\percent} and \SI{2}{\percent} for all values
of $k$ and $\permAngle$. As mentioned, the jump in $\head_2$ across the
fracture is largest for $\permAngle = 0$. In this case, it can be seen
from~\cref{fig:convDiscreteN1_kplot}, that the error is largest for
moderate permeability contrasts, as for $k = \num{1e-3}$ the observed error
is about three times higher as for $k \leq \num{1e-4}$.
For the \tpfaDfm scheme, the errors are about half of those observed on the
finest grid of the convergence test (see~\cref{fig:convDiscreteN1ConduitBulk}).
This seems to be related to the differences in the grids used, and indicates that
the grids obtained by refinement have properties that promote the inconsistency
of the scheme. Moreover, \cref{fig:convDiscreteN1_kplot} illustrates the
dependency of the error on $k$ for the \boxDfm scheme. While for large values of
$k$ the errors are similar to those obtained with the \eboxDfm scheme, larger
errors are observed for decreasing $k$, leading to
$\errorNorm_{\head_2} > \SI{10}{\percent}$ for $k < \num{1e-3}$ and
$\permAngle = 0$. As mentioned, the jump in $\head_2$ across the fracture is
smaller for $\permAngle \in \{\pi/4, \, \pi/2 \}$, which explains the smaller
errors in $\head_2$ for low fracture permeabilities in these cases.
Furthermore, it can be seen that the error in $\head_1$ shows the highest
values for moderate permeability contrasts and $k > 1$, while for $k < 1$, no
strong fluctuations of $\errorNorm_{\head_1}$ can be observed. Most notably,
the error strongly decreases for $k > \num{1e3}$.

\paragraph{Fracture network $\fracNet_2$}
\begin{figure}[ht]
    \begin{subfigure}{0.3299\textwidth}
        \centering
        \includegraphics[width=0.99\textwidth]{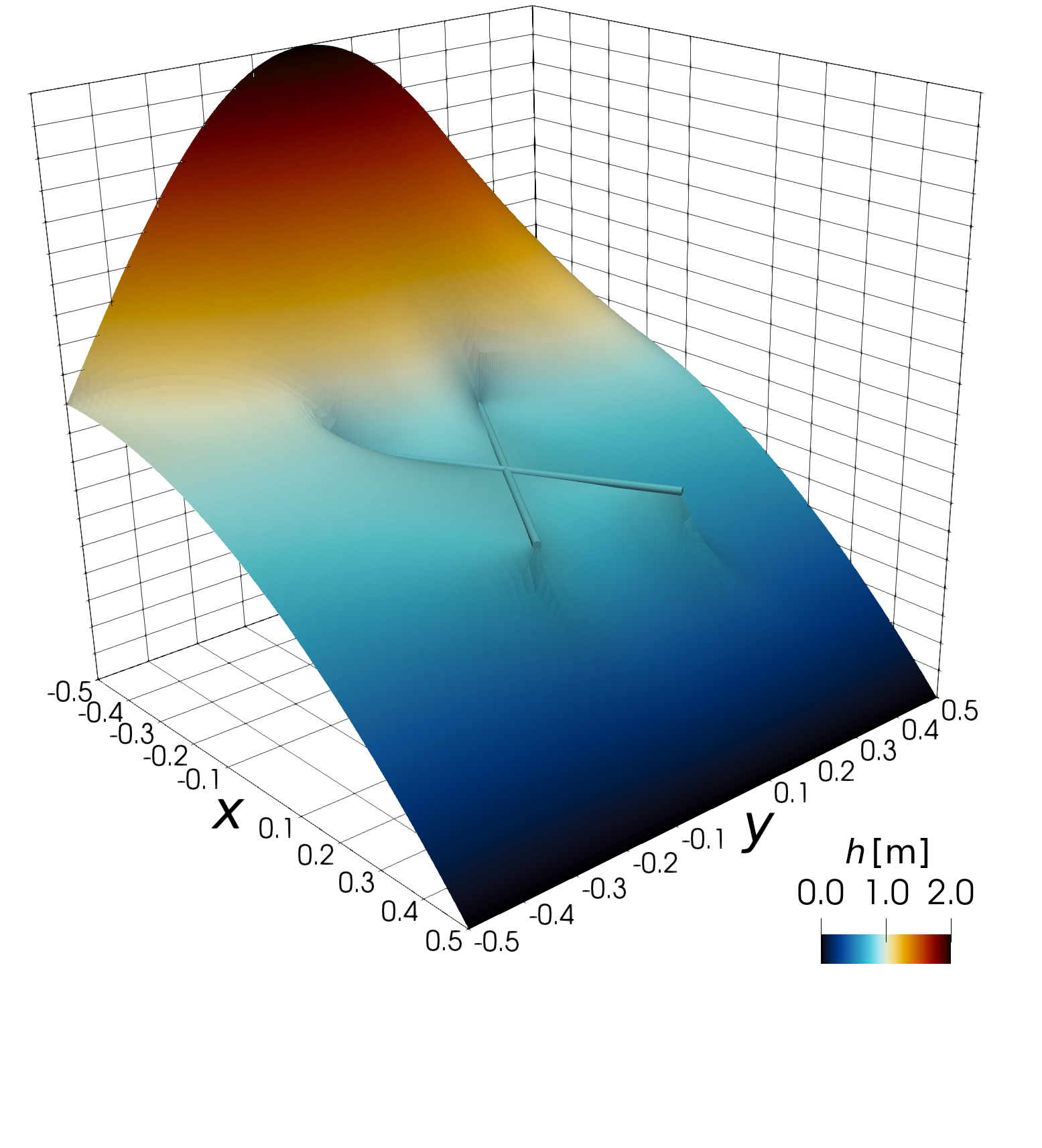}
        \caption{}
        \label{fig:convDiscreteN2ConduitSol}
    \end{subfigure}
    \begin{subfigure}{0.3299\textwidth}
        \centering
        \includegraphics[width=0.99\textwidth]{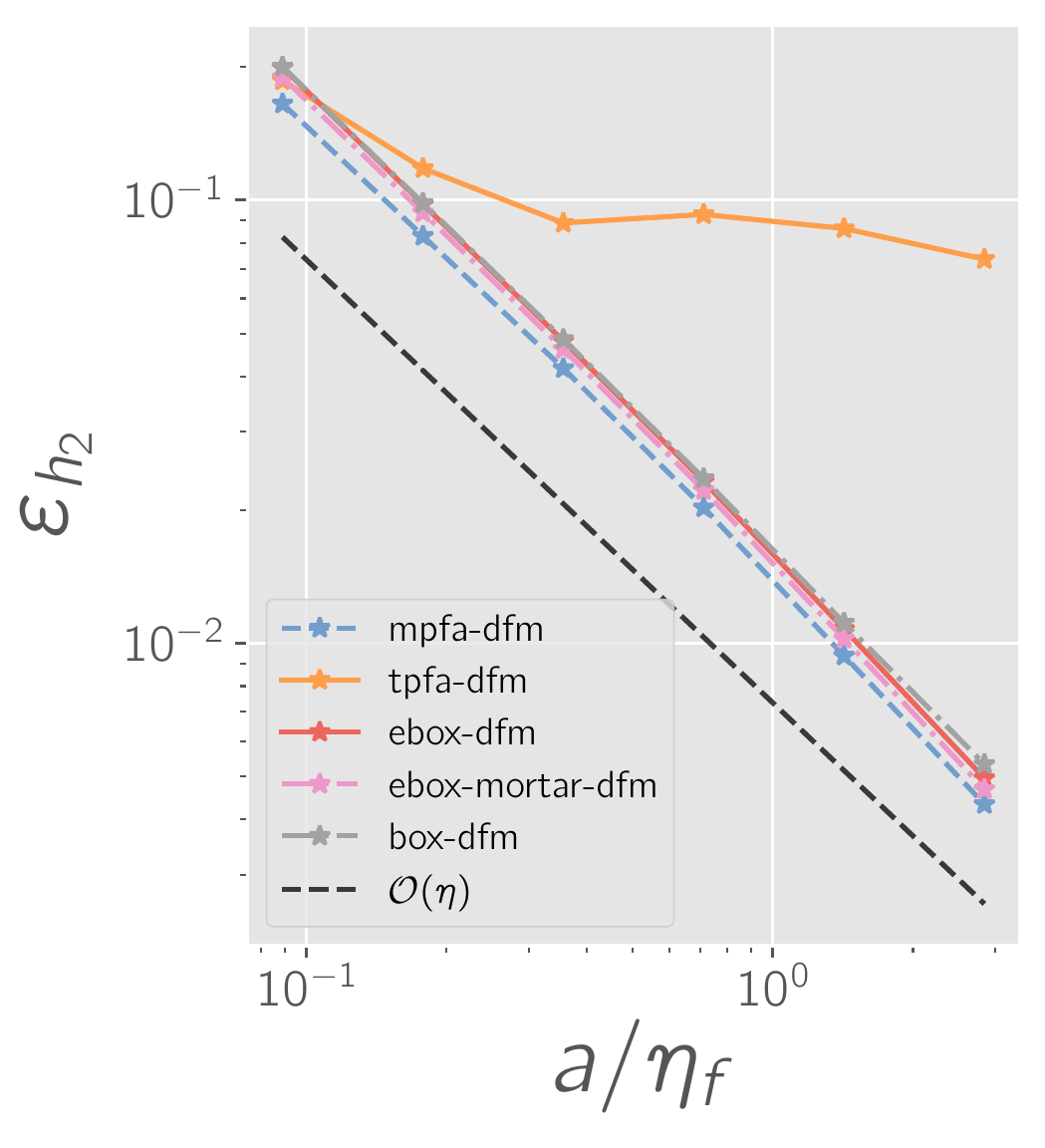}
        \caption{}
        \label{fig:convDiscreteN2ConduitBulk}
    \end{subfigure}
    \begin{subfigure}{0.3299\textwidth}
        \centering
        \includegraphics[width=0.99\textwidth]{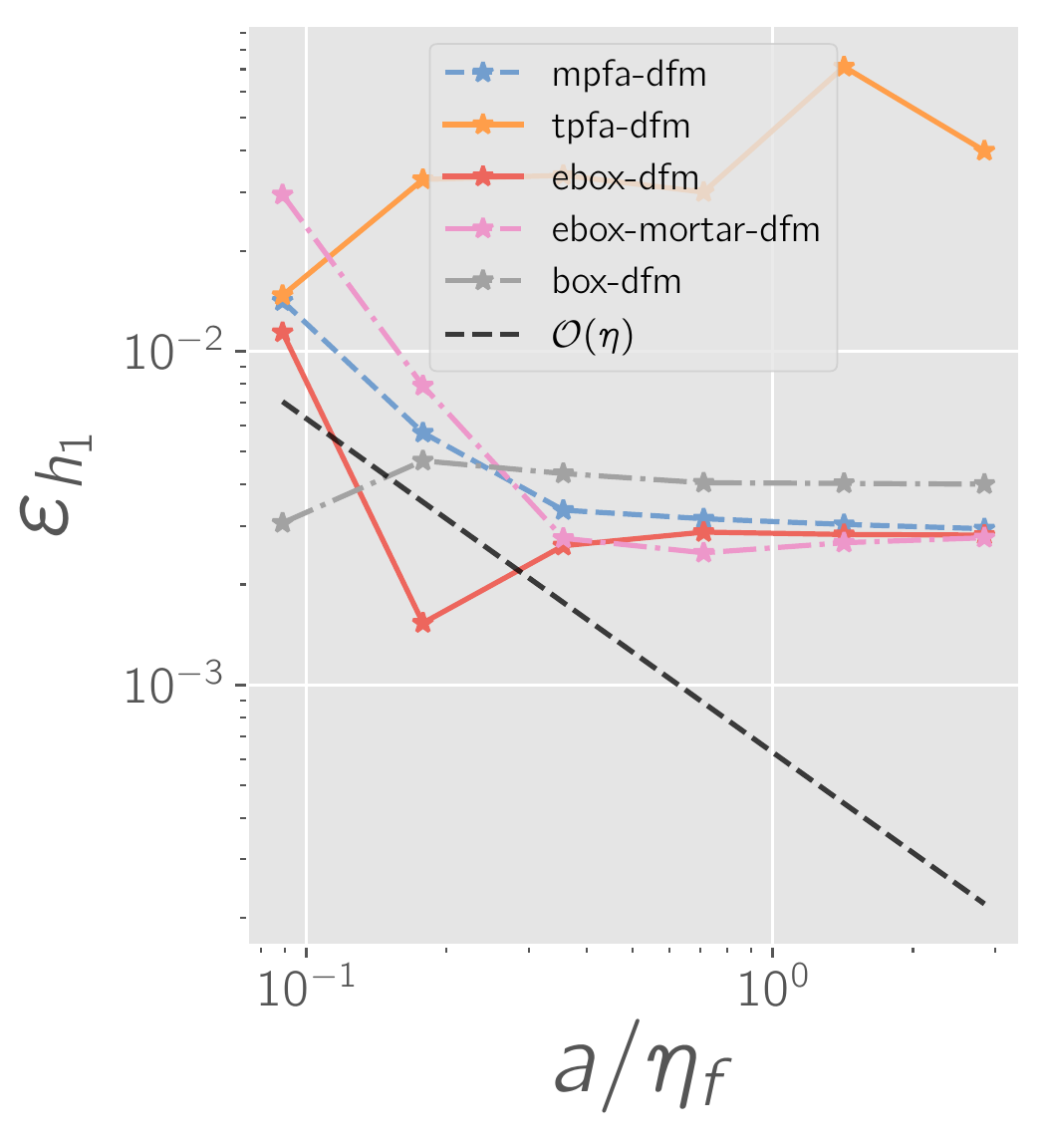}
        \caption{}
        \label{fig:convDiscreteN2ConduitFacet}
    \end{subfigure}
    \begin{subfigure}{0.3299\textwidth}
        \centering
        \includegraphics[width=0.99\textwidth]{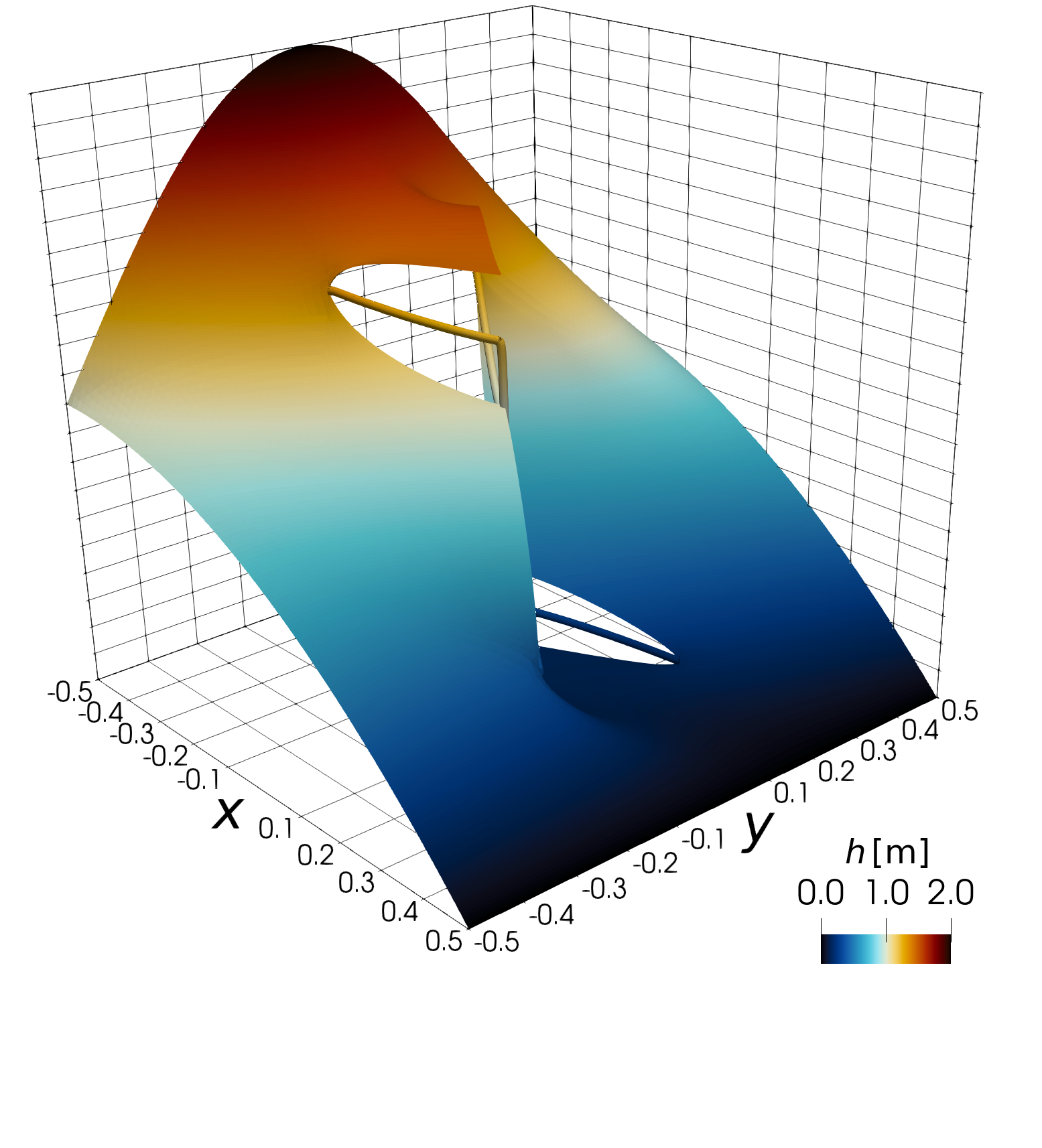}
        \caption{}
        \label{fig:convDiscreteN2BarrierSol}
    \end{subfigure}
    \begin{subfigure}{0.3299\textwidth}
        \centering
        \includegraphics[width=0.99\textwidth]{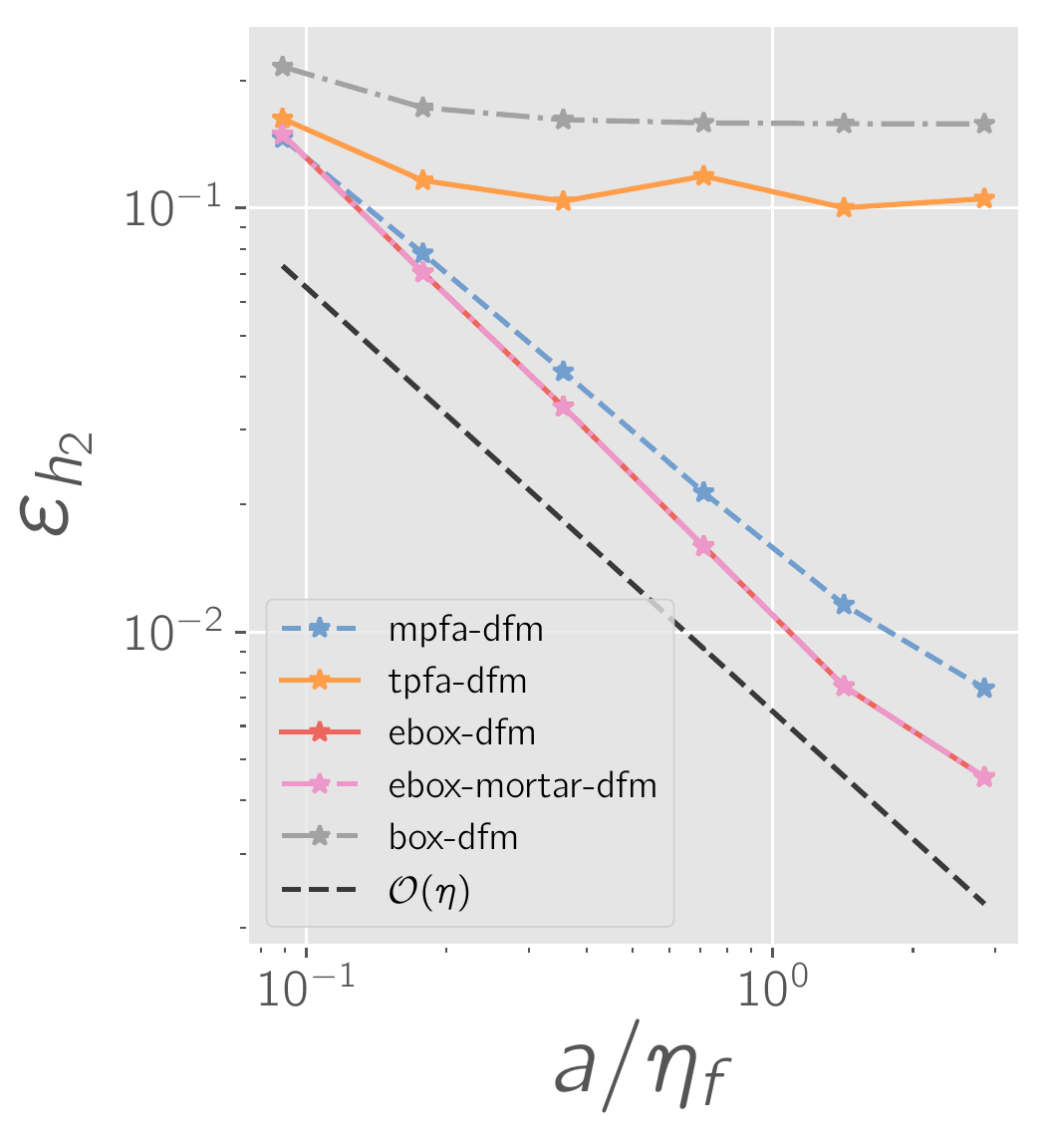}
        \caption{}
        \label{fig:convDiscreteN2BarrierBulk}
    \end{subfigure}
    \begin{subfigure}{0.3299\textwidth}
        \centering
        \includegraphics[width=0.99\textwidth]{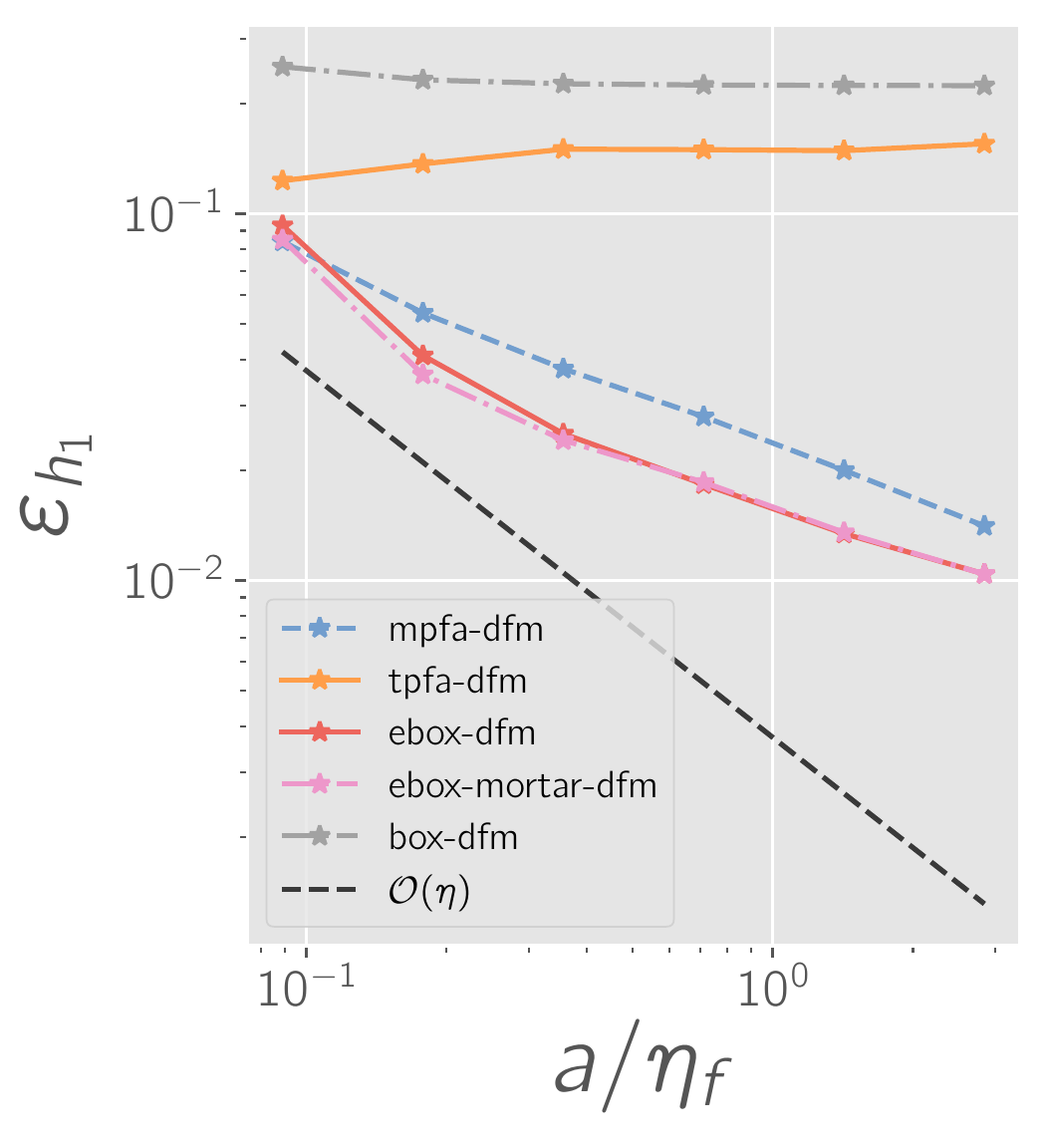}
        \caption{}
        \label{fig:convDiscreteN2BarrierFacet}
    \end{subfigure}
    \caption{\textbf{Case 2.2 - results for $\permAngle = \pi/4$}. The upper row and lower row show the results for $k = \num{1e4}$ and $k = \num{1e-4}$, respectively. The left column depicts the distributions of the hydraulic heads in the domain, obtained with the \eboxDfm scheme on the finest grid,
    while the middle and right columns show the errors in $\head_2$ and $\head_1$ plotted over grid refinement.}
    \label{fig:convDiscreteN2}
\end{figure}
The geometric complexity is increased with this fracture network as it contains an
intersection of two fractures (see~\cref{fig:convDiscreteNetwork2}). In the
mixed-dimensional model, we treat $0$-dimensional fracture intersections by the
conditions~\eqref{eq:prob_mixed_junctions}, that is, the intersection region is not
described by an individual subdomain. However, in the equi-dimensional discretization,
the intersection region is resolved and populated with grid elements, and therefore,
we assign the
permeability $\perm_\fracIdx$ to those elements in order to minimize deviations
in the setup.\\

In the following we want to focus on the discussion of the results for
$\permAngle = \pi/4$, which are again qualitatively similar to those obtained
for $\permAngle = 0$ and which can be found in~\cref{tab:convDiscreteErrors_head2_angle0,tab:convDiscreteErrors_head1_angle0}.
The solution and plots of the error norms can be seen in~\cref{fig:convDiscreteN2},
which shows that very similar errors in $\head_2$, both qualitatively and quantitatively, are
obtained when compared to those for $\fracNet_1$. Concerning $\head_1$, we again
observe that convergence is lost after a few refinements, leading to very similar
irreducible relative errors after the last refinement. Noticeable differences
to the results obtained on $\fracNet_1$ can be seen for $\head_1$ and $k = \num{1e-4}$,
where lower convergence rates are obtained already after the first refinement.
This is probably related to the fact that in the presence of fracture intersections
in low-permeable fracture networks, complex distributions of the hydraulic head within
the intersection region might develop, which are not captured by the mixed-dimensional
model. Nevertheless, the relative error is reduced to $\approx \SI{1}{\percent}$
in the last refinement, while errors below $\SI{1}{\percent}$ are again reached
for $\head_2$.\\

\begin{figure}[ht]
  \begin{subfigure}{\textwidth}
    \centering
    \includegraphics[width=0.99\textwidth]{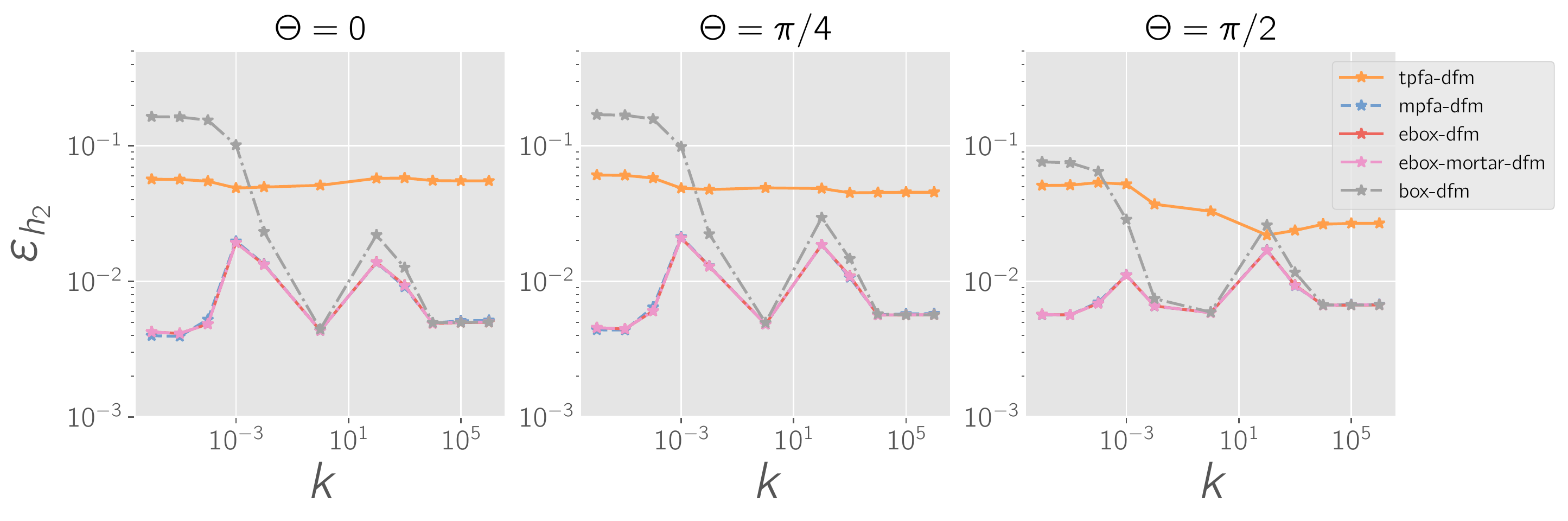}
  \end{subfigure}
  \begin{subfigure}{\textwidth}
    \centering
    \includegraphics[width=0.99\textwidth]{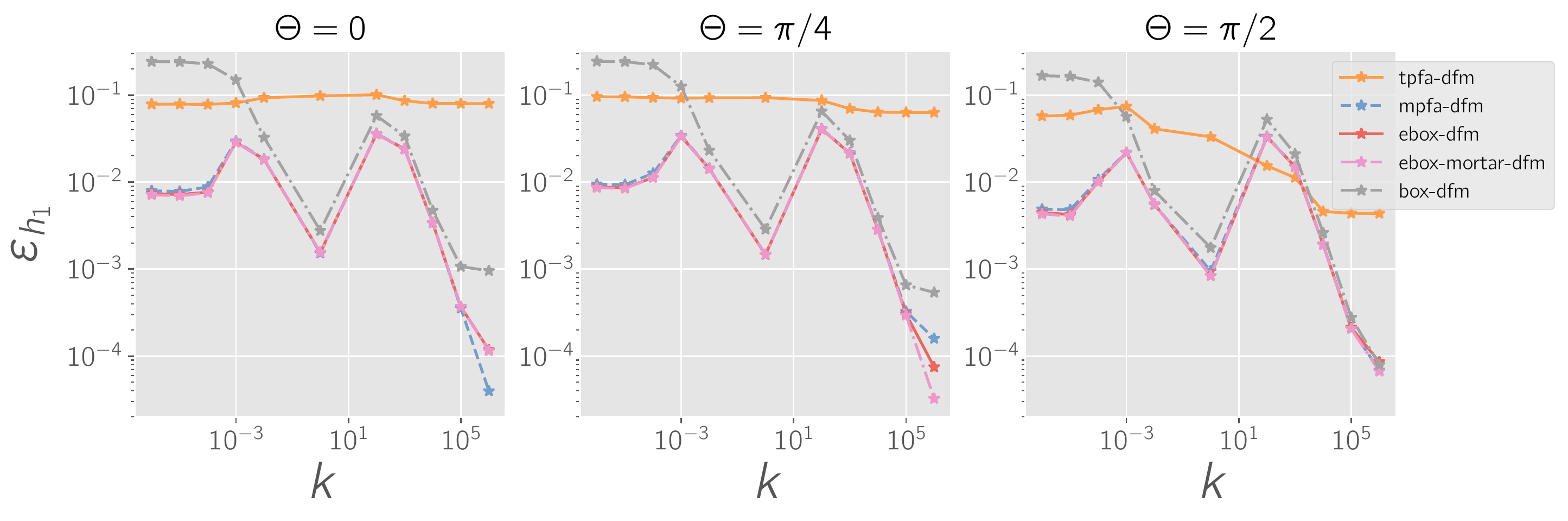}
  \end{subfigure}
  \caption{\textbf{Case 2.2 - error over $k$}. The errors in $\head_2$ and $\head_1$, obtained on the meshes as specified in~\cref{tab:convtest2MeshesAdditional}, are plotted over $k$ for various values of the permeability angle $\permAngle$.}
  \label{fig:convDiscreteN2_kplot}
\end{figure}

\Cref{fig:convDiscreteN2_kplot} depicts the errors in $\head_2$ and $\head_1$
obtained from additional simulations using a mixed-dimensional discretization
as given in~\cref{tab:convtest2MeshesAdditional}. As for $\fracNet_1$
(see~\cref{fig:convDiscreteN1_kplot}), the error
seems to be largest for moderate permeability contrasts, while on this network,
we observe less dependency on the permeability angle $\permAngle$. This can be
explained by the fact that the jump in hydraulic head across the fractures is less
dependent on $\permAngle$ in the case of $\fracNet_2$, where there is always a fracture
oriented in orthogonal direction to the flow field. Consequently, significant
errors can be seen for the \boxDfm scheme and low permeabilities for all values
of $\permAngle$. A noteworthy difference to the results of $\fracNet_1$ is that
the errors in $\head_1$ are about an order of magnitude larger for $k < \num{1e-3}$,
which reflects the lower convergence rates observed in this case, as was
discussed earlier.\\

\begin{figure}[ht]
    \begin{subfigure}{0.3299\textwidth}
        \centering
        \includegraphics[width=0.99\textwidth]{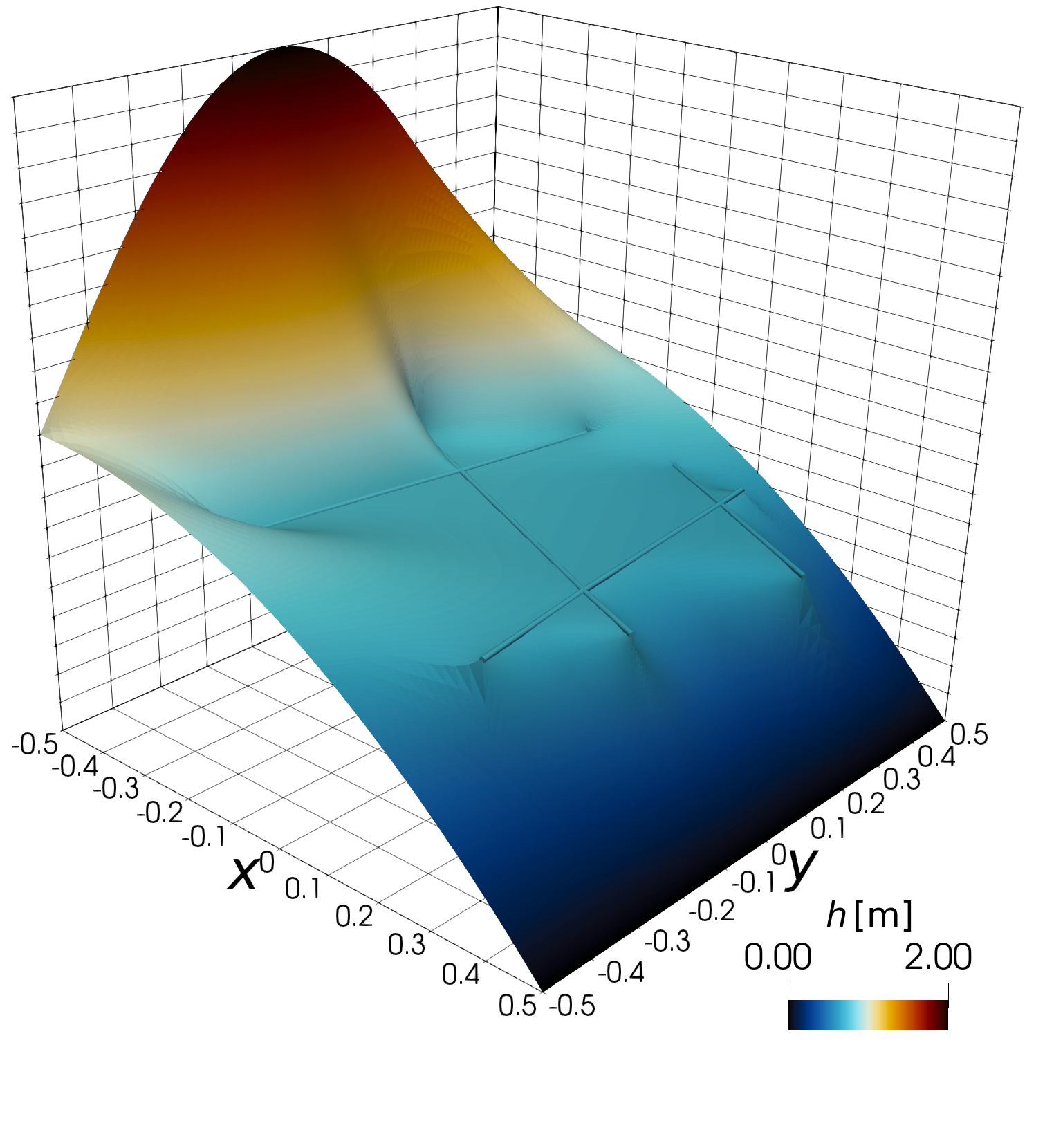}
        \caption{}
        \label{fig:convDiscreteN3ConduitSol}
    \end{subfigure}
    \begin{subfigure}{0.3299\textwidth}
        \centering
        \includegraphics[width=0.99\textwidth]{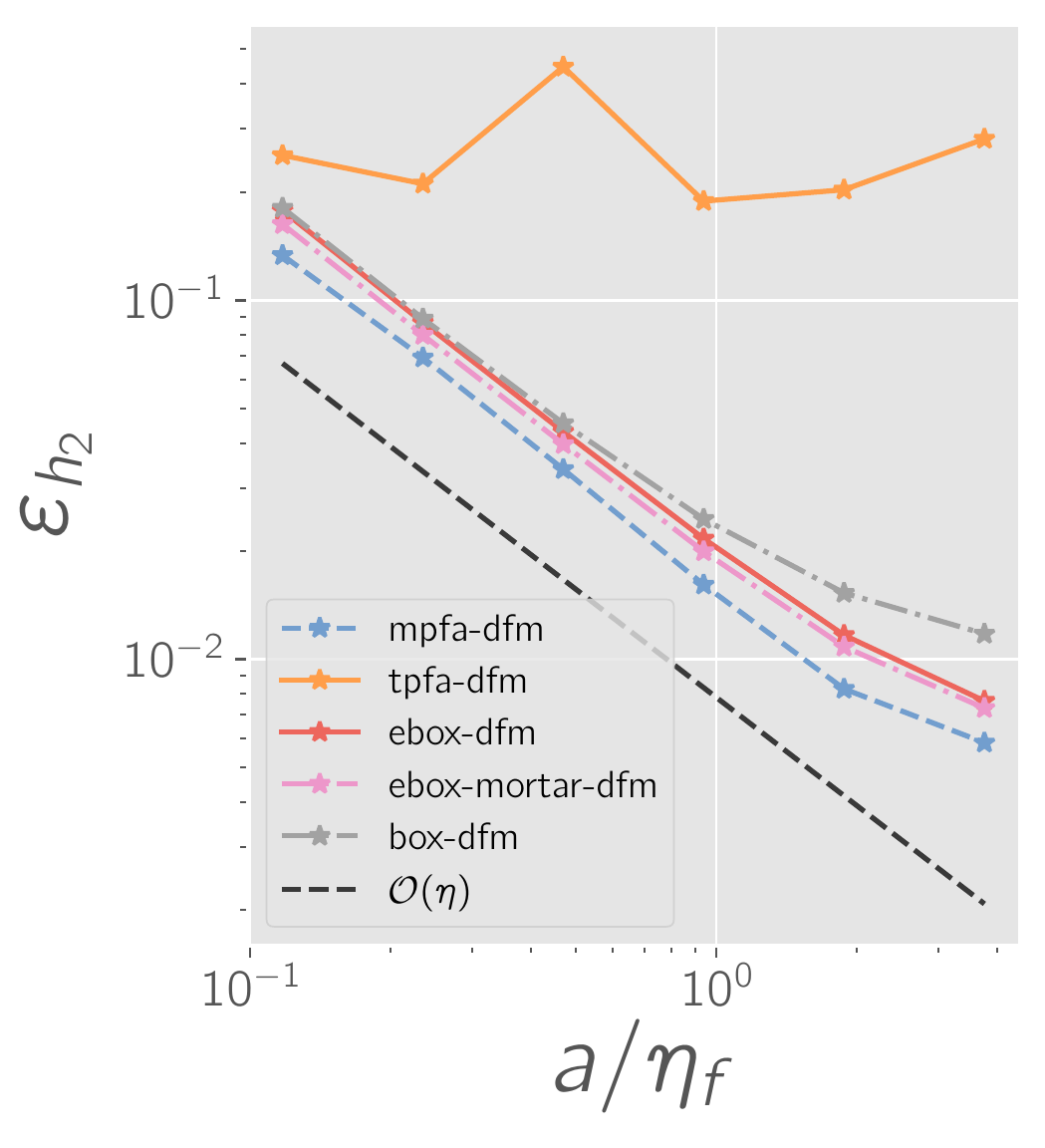}
        \caption{}
        \label{fig:convDiscreteN3ConduitBulk}
    \end{subfigure}
    \begin{subfigure}{0.3299\textwidth}
        \centering
        \includegraphics[width=0.99\textwidth]{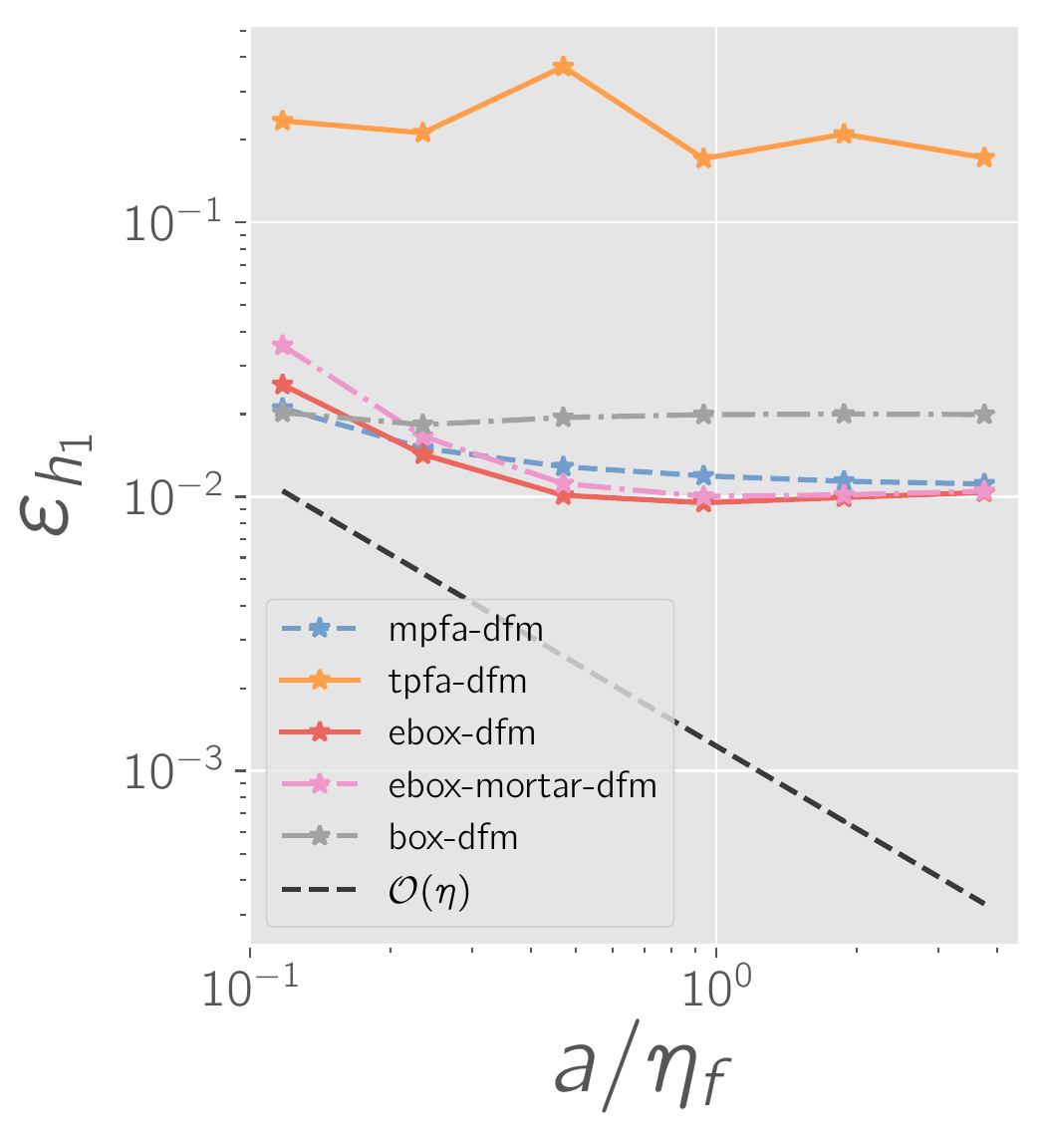}
        \caption{}
        \label{fig:convDiscreteN3ConduitFacet}
    \end{subfigure}
    \begin{subfigure}{0.3299\textwidth}
        \centering
        \includegraphics[width=0.99\textwidth]{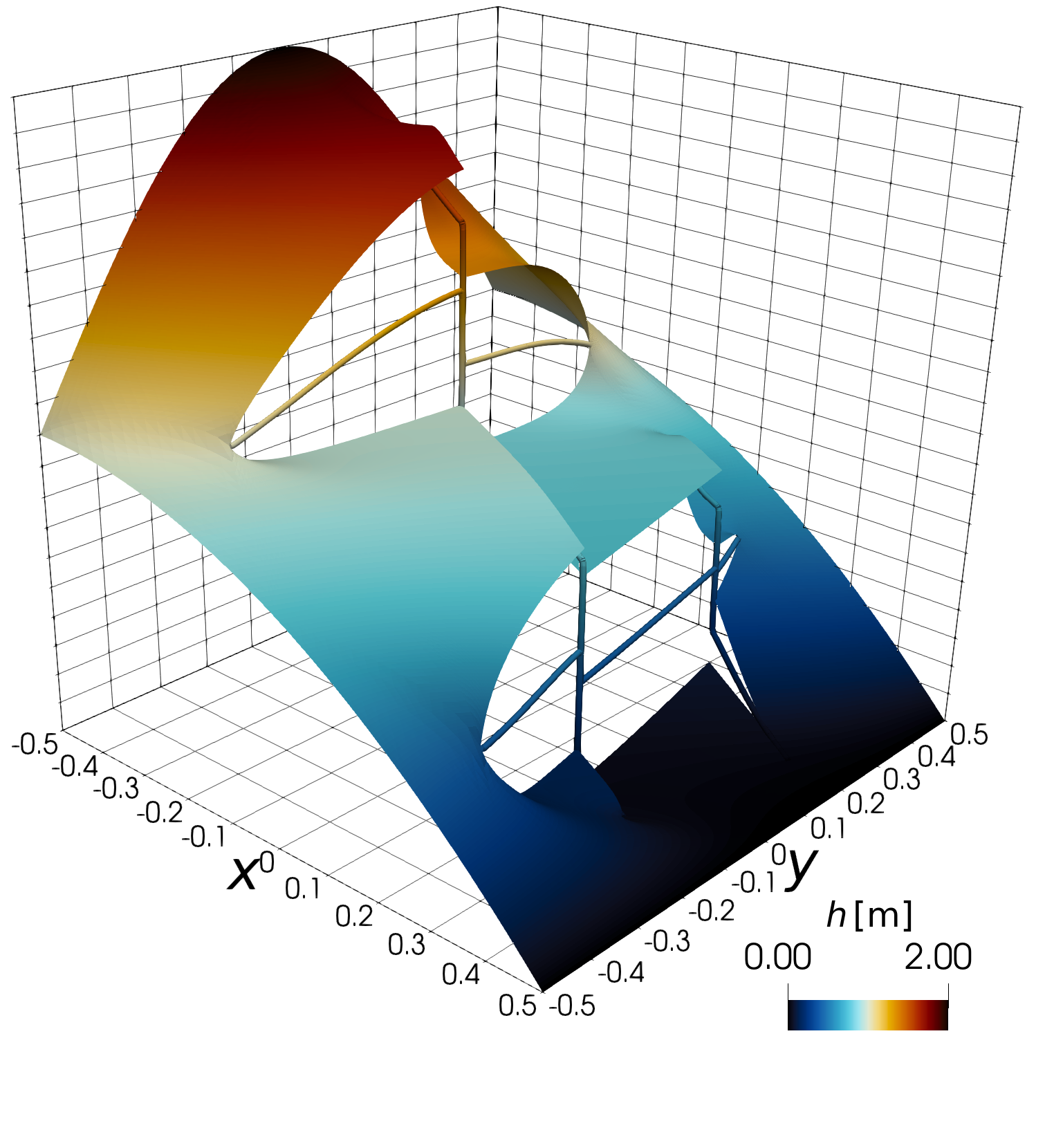}
        \caption{}
        \label{fig:convDiscreteN3BarrierSol}
    \end{subfigure}
    \begin{subfigure}{0.3299\textwidth}
        \centering
        \includegraphics[width=0.99\textwidth]{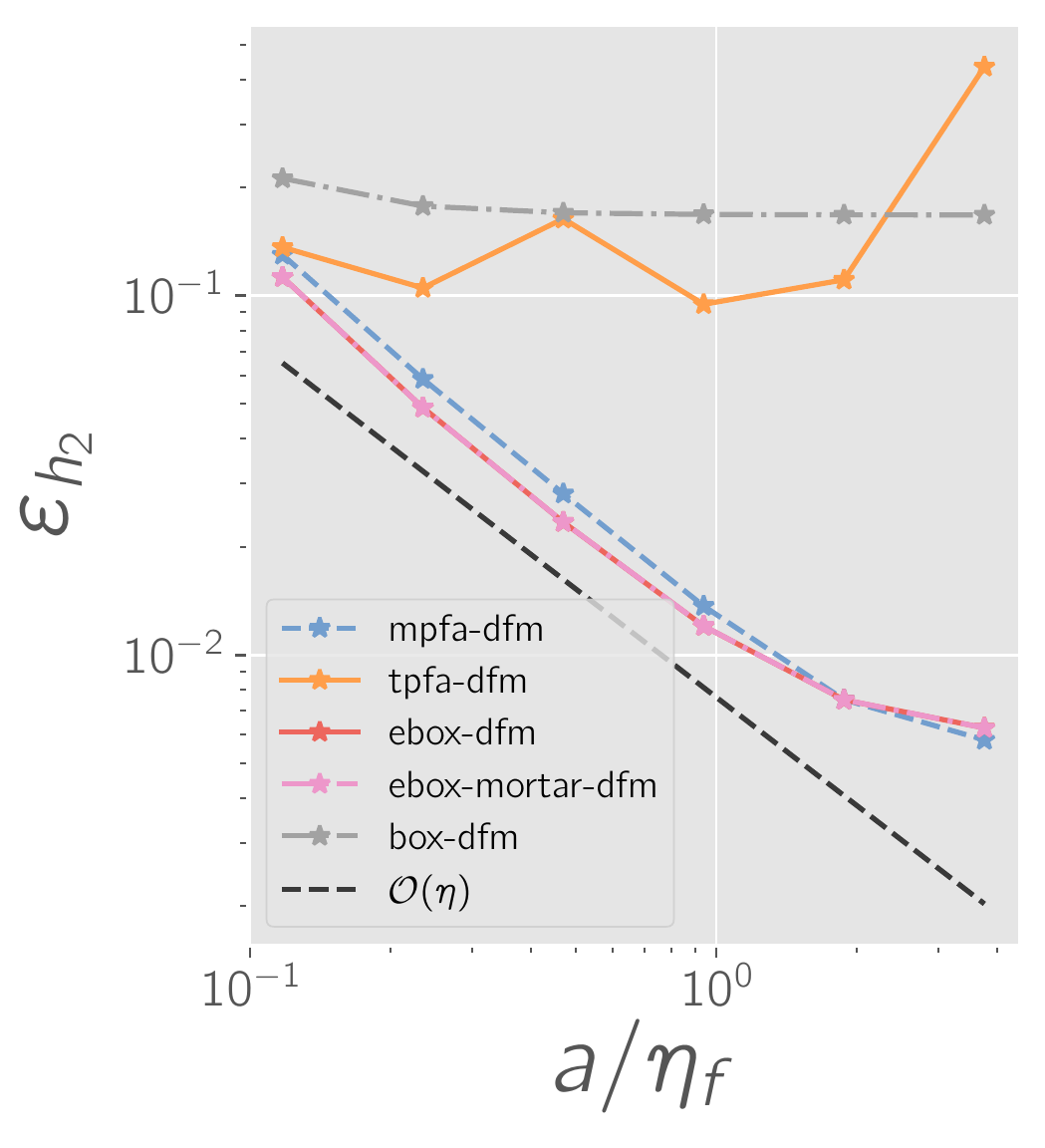}
        \caption{}
        \label{fig:convDiscreteN3BarrierBulk}
    \end{subfigure}
    \begin{subfigure}{0.3299\textwidth}
        \centering
        \includegraphics[width=0.99\textwidth]{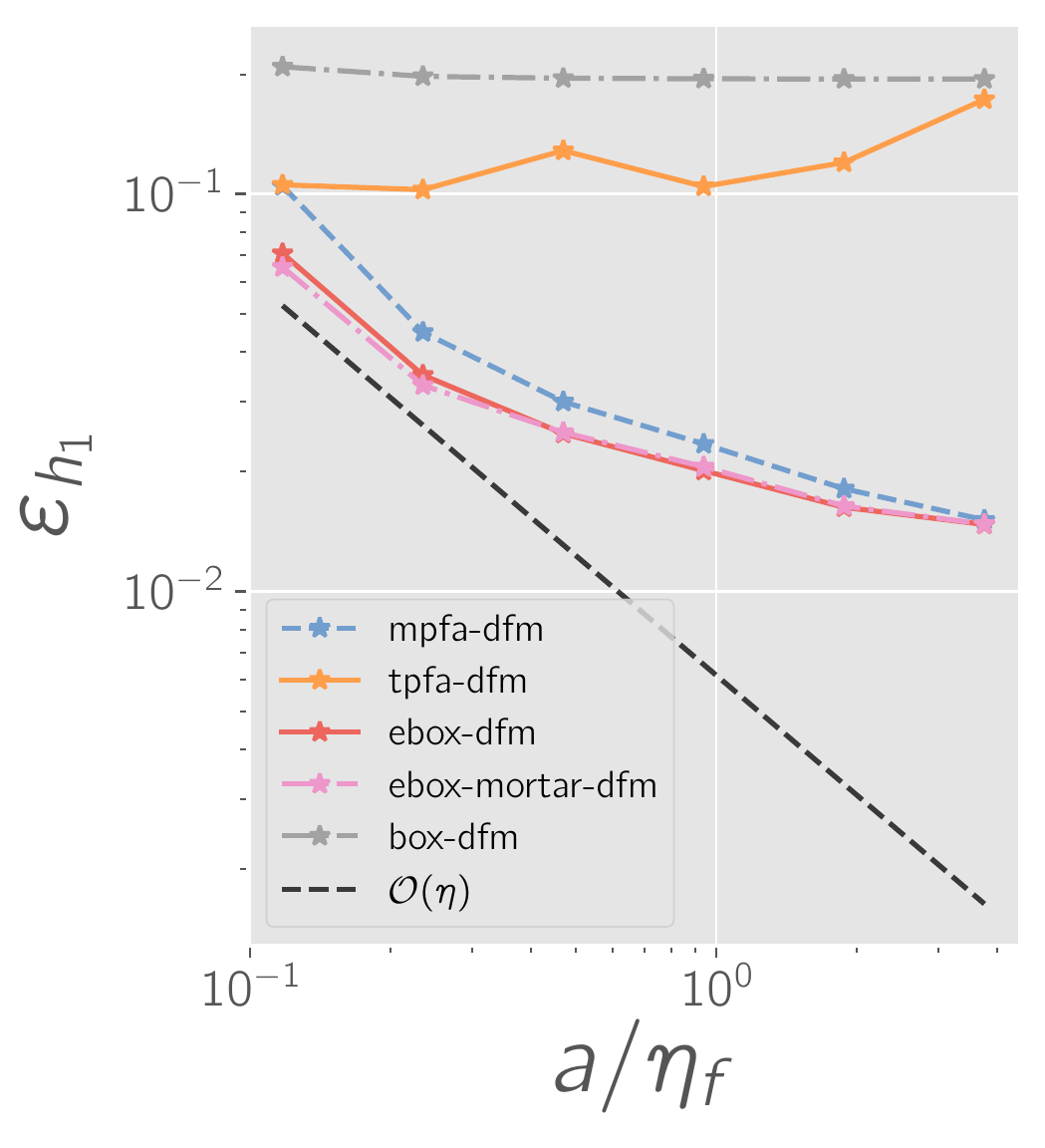}
        \caption{}
        \label{fig:convDiscreteN3BarrierFacet}
    \end{subfigure}
    \caption{\textbf{Case 2.3 - results for $\permAngle = \pi/4$}. The upper row and lower row show the results for $k = \num{1e4}$ and $k = \num{1e-4}$, respectively. The left column depicts the distributions of the hydraulic heads in the domain, obtained with the \eboxDfm scheme on the finest grid,
    while the middle and right columns show the errors in $\head_2$ and $\head_1$ plotted over grid refinement.}
    \label{fig:convDiscreteN3}
\end{figure}

\paragraph{Fracture network $\fracNet_3$}
This is the most complex of the networks considered in this test case,
consisting of four intersecting fractures. In the discussion, we want to
focus again on the results obtained for $\permAngle = \pi/4$, which are depicted
in~\cref{fig:convDiscreteN3ConduitBulk}, while those for $\permAngle = 0$ can be
found in~\cref{tab:convDiscreteErrors_head2_angle0,tab:convDiscreteErrors_head1_angle0}
in the appendix. In contrast to the previous two
networks, the rate of convergence in $\head_2$, for $k = \num{1e4}$, starts
to decrease upon the last refinement. Moreover, only a small reduction of the error
$\errorNorm_{\head_1}$ in the fractures is observed, which stagnates at $\approx \SI{1}{\percent}$.\\

However, despite the significantly increased complexity, the errors obtained on the finest
grid are comparable to those obtained on $\fracNet_2$ (see~\cref{fig:convDiscreteN2}),
which can also be seen in~\cref{fig:convDiscreteN2_kplot}, where the errors for a wider
range of values for $k$, computed on the mesh listed in~\cref{tab:convtest2MeshesAdditional},
is shown. Both qualitatively and quantitatively, very
similar results to those obtained on $\fracNet_2$ (see~\cref{fig:convDiscreteN1_kplot})
are observed, while the main difference is that we do not see an increase in
$\errorNorm_{\head_1}$ for moderately low-permeable fractures on
$\fracNet_3$.

\begin{figure}[ht]
  \begin{subfigure}{\textwidth}
    \centering
    \includegraphics[width=0.99\textwidth]{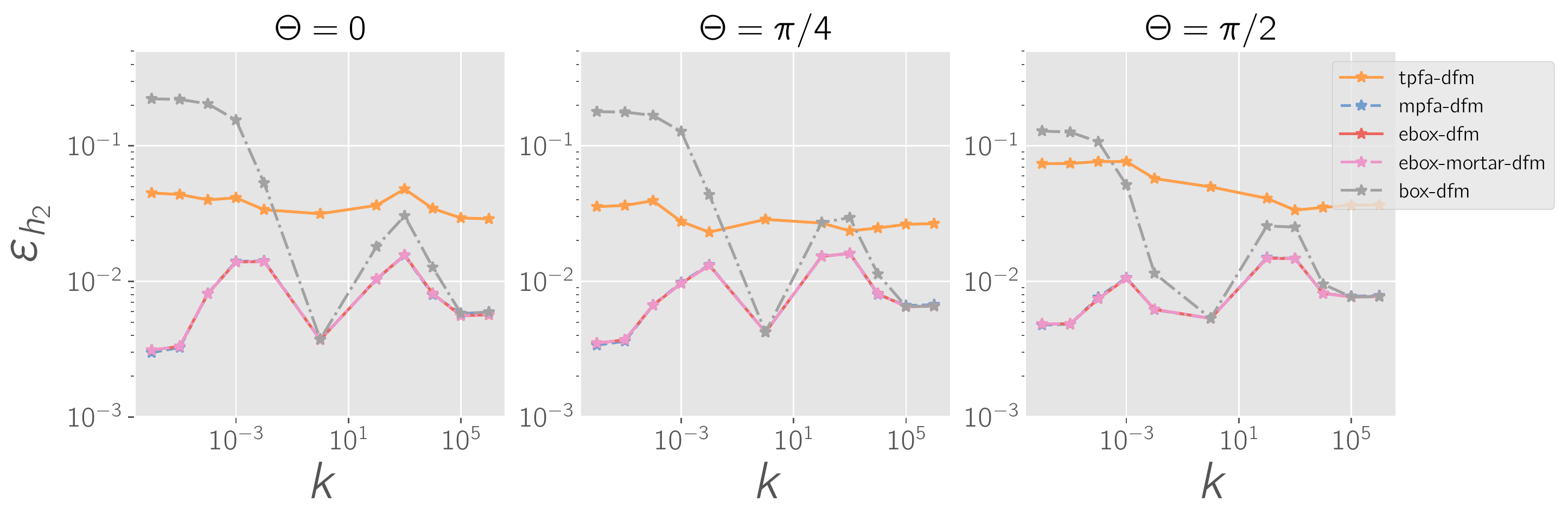}
  \end{subfigure}
  \begin{subfigure}{\textwidth}
    \centering
    \includegraphics[width=0.99\textwidth]{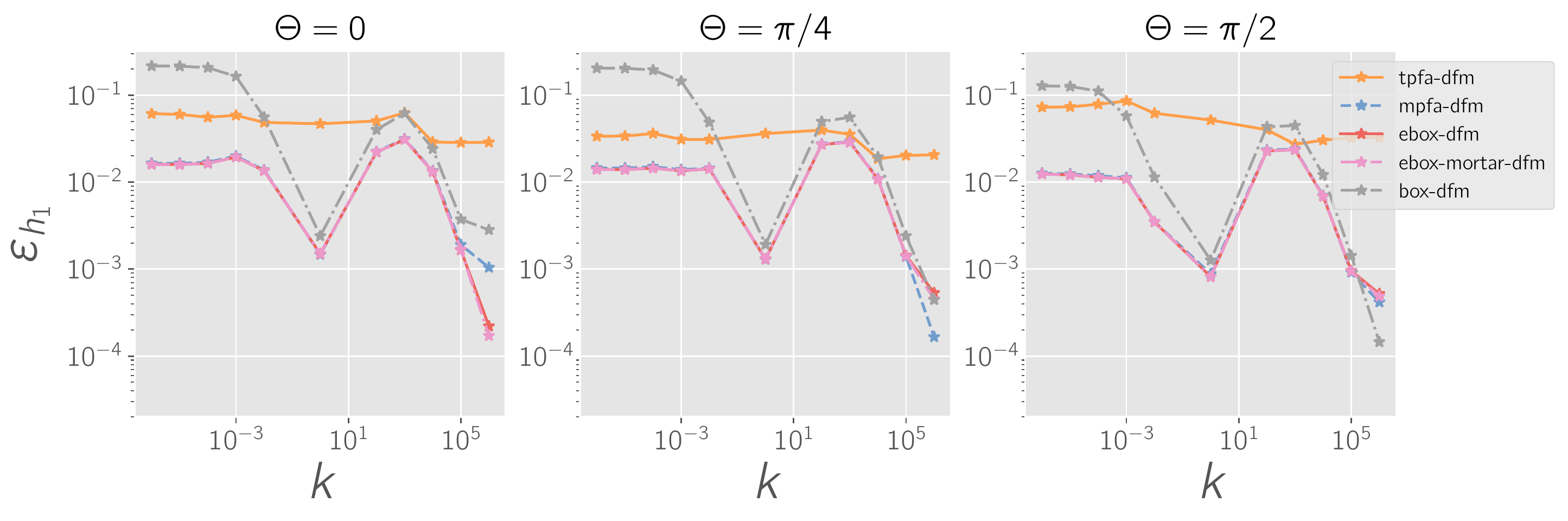}
  \end{subfigure}
  \caption{\textbf{Case 2.3 - error over $k$}. The errors in $\head_2$ and $\head_1$, obtained on the meshes as specified in~\cref{tab:convtest2MeshesAdditional}, are plotted over $k$ for various values of the permeability angle $\permAngle$.}
  \label{fig:convDiscreteN3_kplot}
\end{figure}

\paragraph{Summary}
The results show that the \eboxDfm and \eboxMortarDfm schemes, being applicable
to both low- and highly-permeable fractures and arbitrary bulk permeability angles
$\permAngle$, lead to similar or slightly lower errors as the \mpfaDfm scheme.
Convergence in $\head_2$ was maintained over all grid refinements, while the errors
in $\head_1$ showed to stagnate when the discretization length approaches the
aperture. With increasing complexity of the fracture network, this stagnation
seemed to occur earlier in the refinement cycle, at least for $k = \num{1e4}$,
while for $k = \num{1e-4}$, lower convergence rates were observed already in
the first refinements. This is probably due to the complex jumps in hydraulic head
that develop across intersecting fractures, which the mixed-dimensional model
struggles to capture.\\

The \boxDfm scheme leads to comparable errors for large values of $k$, for which its
model assumptions are fulfilled. However, for moderate permeability contrasts
it led to larger errors when compared to the \mpfaDfm or the other
vertex-centered schemes, and for $k \leq \num{1e-3}$, errors above \SI{10}{\percent}
were observed. Despite the similarity in the results between the \mpfaDfm and the
\eboxDfm and \eboxMortarDfm schemes, the vertex-centered schemes use much
less degrees of freedom on unstructured grids and lead to smaller stencils.
This issue will be investigated in more detail in~\cref{sec:benchmark2d,sec:benchmark3d}.

\subsection{Case 3: interface flux and hydraulic head}
\label{sec:interface}

This test case focuses on the vertex-centered \eboxDfm and \eboxMortarDfm schemes,
and in particular, on the bulk-fracture interface flux and hydraulic head. To this
end, we again consider a unit square bulk domain, and reuse~\cref{eq:convtest2DirichletBcs}
as Dirichlet boundary conditions on the entire bulk domain boundary. A single
fracture is placed in the domain, defined as the connecting segment between the points
$\left(-0.15, 0.35\right)$ and $\left(0.2, -0.35\right)$. As in the previous test
case, we consider an aperture of $\aperture = \num{2e-3}$ and generate an
equi-dimensional discretization of the domain following the procedure described
in~\cref{sec:convDiscrete}, on which a reference solution is computed with the
\mpfa scheme. A rather coarse mixed-dimensional discretization is used, consisting
of \num{28266} triangular and \num{235} line elements, while the reference grid
contains \num{532384} triangular elements. We are particularly interested in
differences among the two schemes in fulfilling the conditions~\eqref{eq:prob_mixed_if_p}
at interfaces to highly-permeable or open fractures. Therefore, we consider
the permeabilities used in~\cref{sec:convDiscrete} with $k = \num{1e4}$.\\
\begin{figure}[ht]
  \begin{subfigure}{\textwidth}
    \centering
    \includegraphics[width=0.99\textwidth]{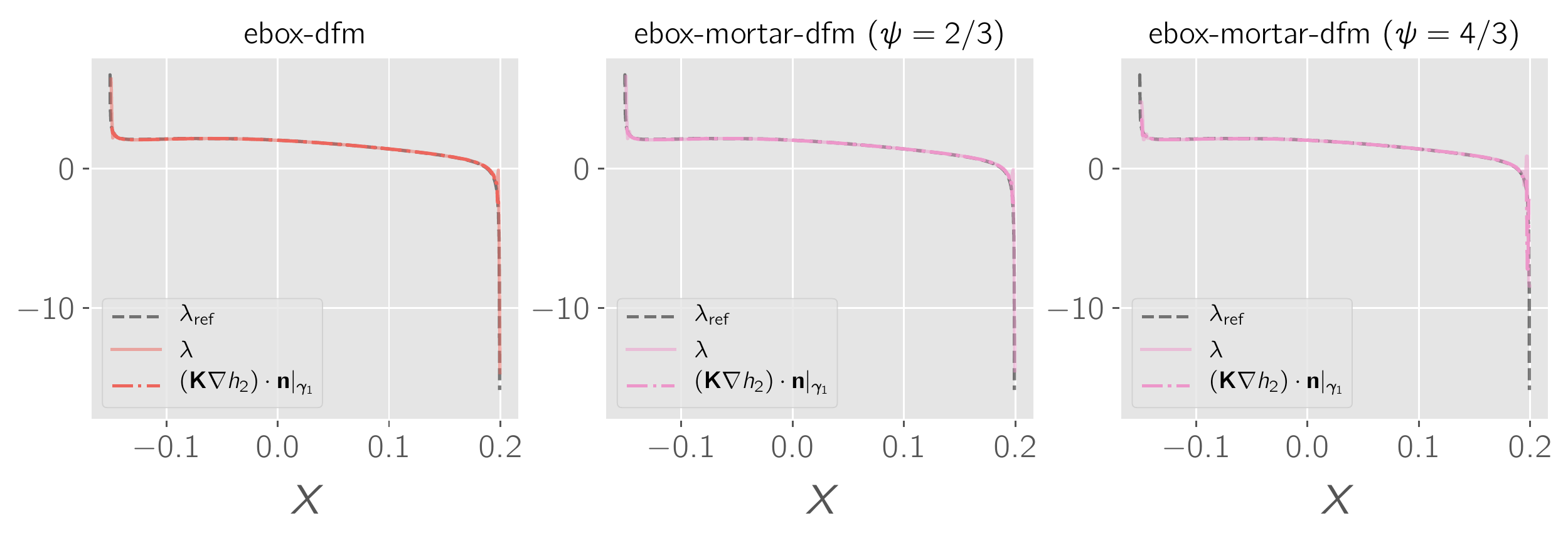}
    \end{subfigure}
  \begin{subfigure}{\textwidth}
    \centering
    \includegraphics[width=0.99\textwidth]{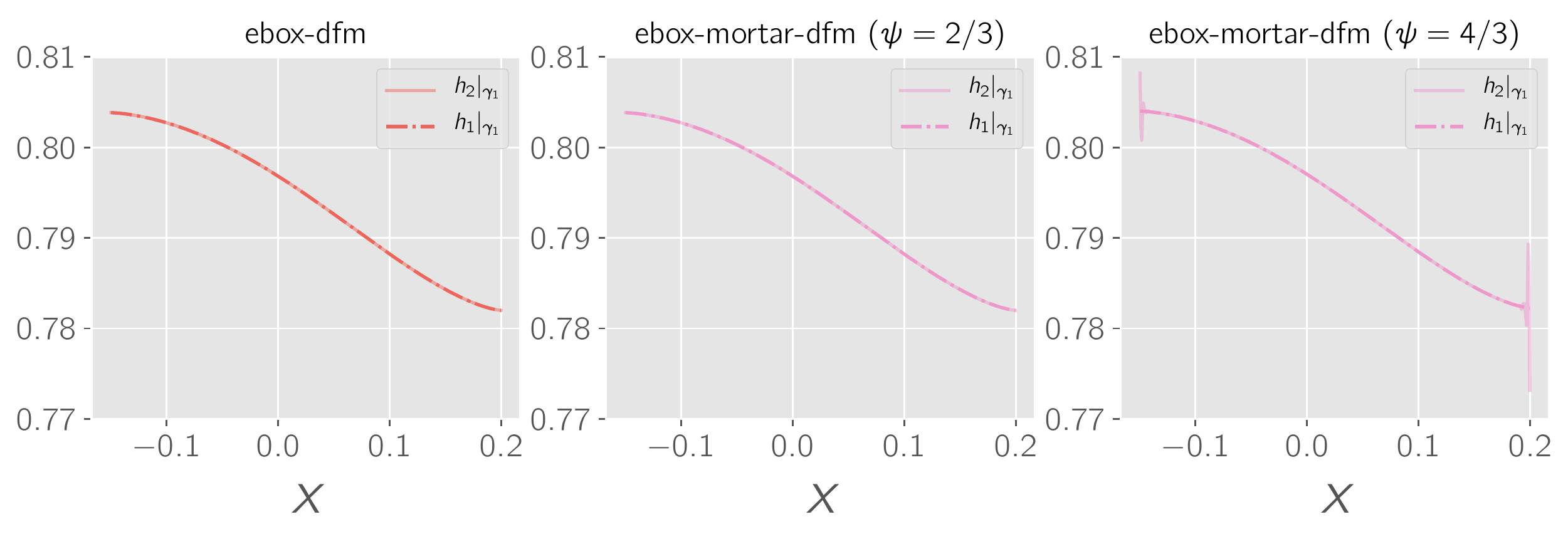}
  \end{subfigure}
  \caption{\textbf{Case 3.1 - interface flux and hydraulic head}. The top and bottom rows show the transfer fluxes and the hydraulic heads along the interface $\interface_1$, respectively.}
  \label{fig:iffluxcase_neum}
\end{figure}
Let us first look at the results for the case of
$\faceSet_{d, \couplIdx} \equiv \faceSet_{d, \blockInterfaceSet}$, that is, using
the conditions~\eqref{eq:prob_mixed_if_flux} at the bulk-fracture interface.~\Cref{fig:iffluxcase_neum}
shows the transfer flux $\lambda$, the term
$\left( \perm \grad \head_2 \right) \scal \n |_{\interface_1}$ and the hydraulic heads
$\head_2$ and $\head_1$ in the bulk medium and the fractures, plotted along the
interface $\interface_1$. Recall that for immersed fractures it is
$\numFractureSides_{1}^{i,j} = 2$, that is, the interface $\interface_1^{i,j}$
consists of both sides of the fracture on which fluxes and hydraulic heads might
have different values. Therefore, the plots in~\cref{fig:iffluxcase_neum} are
restricted to the lower left fracture side on which
$\mathbf{e}_1 \cdot \n |_{\interface_1} > 0$, $\mathbf{e}_1 = \left(1, 0 \right)^T$,
holds for the outer normal vector.\\

From the upper row of~\cref{fig:iffluxcase_neum}
it can be seen that the transfer fluxes produced by the \eboxDfm and the
\eboxMortarDfm schemes (with $\mortarGridFactor \in \{2/3, \, 4/3\}$) agree very
well with the reference solution. A steep increase
in $\meas{\lambda}$ can be seen on the left and right sides of the plot, which
implies large in- and outfluxes concentrated at the fracture tips, resulting
from the fact that the fracture is much more permeable than the surrounding
bulk medium. Small oscillations in $\lambda$ in the discrete solutions can be
seen right at the fracture tips, and we observed these oscillations to stay local
within the first adjacent grid cell, not polluting the solution in the entire
domain. However, the oscillations are not reflected in the term
$\left( \perm \grad \head_2 \right) \scal \n |_{\interface_1}$. Furthermore,
it is noteworthy that for the \eboxMortarDfm scheme, the magnitude of the
oscillation seems to be slightly increased for $\mortarGridFactor = 4/3$ in
comparison to $\mortarGridFactor = 1/3$. Recall that for $\mortarGridFactor = 0.5$,
the \eboxMortarDfm scheme produces the same results as the \eboxDfm scheme, which
is why these results are not included here. Although also not shown here, for
$\mortarGridFactor = 1.0$, we observed spurious oscillations in $\lambda$ across
the entire interface, and we will return to this later in the discussion.\\

Recall that
$\lambda \propto \frac{k}{\aperture} \left( \head_1 - \head_2 \right) |_{\interface_1}$
according to condition~\eqref{eq:prob_mixed_if_flux}. Therefore, with the
permeability in the fracture being much higher than in the bulk medium, the
difference in the hydraulic heads $\head_2$ and $\head_1$ is very small and not
visible in the plots shown in the lower row of~\cref{fig:iffluxcase_neum}. In
particular, due to the above-mentioned proportionality, we expect small oscillations
in the hydraulic head to be present at the fracture tips since they are observed
in $\lambda$.
For the \eboxMortarDfm scheme with $\mortarGridFactor > 1$, noticeable
oscillations of $\head_2$ appear at the fracture tips, as can be seen in the lower
right plot in~\cref{fig:iffluxcase_neum}. However, the projections of the
hydraulic heads used in~\cref{eq:ifConditionsMortarBlock}, that is, the discrete variant of
the condition~\eqref{eq:prob_mixed_if_flux}, seem to level out these oscillations
such that they do not strongly influence the transfer flux $\lambda$ as shown in the
plot directly above.\\

\begin{figure}[ht]
  \begin{subfigure}{\textwidth}
    \centering
    \includegraphics[width=0.99\textwidth]{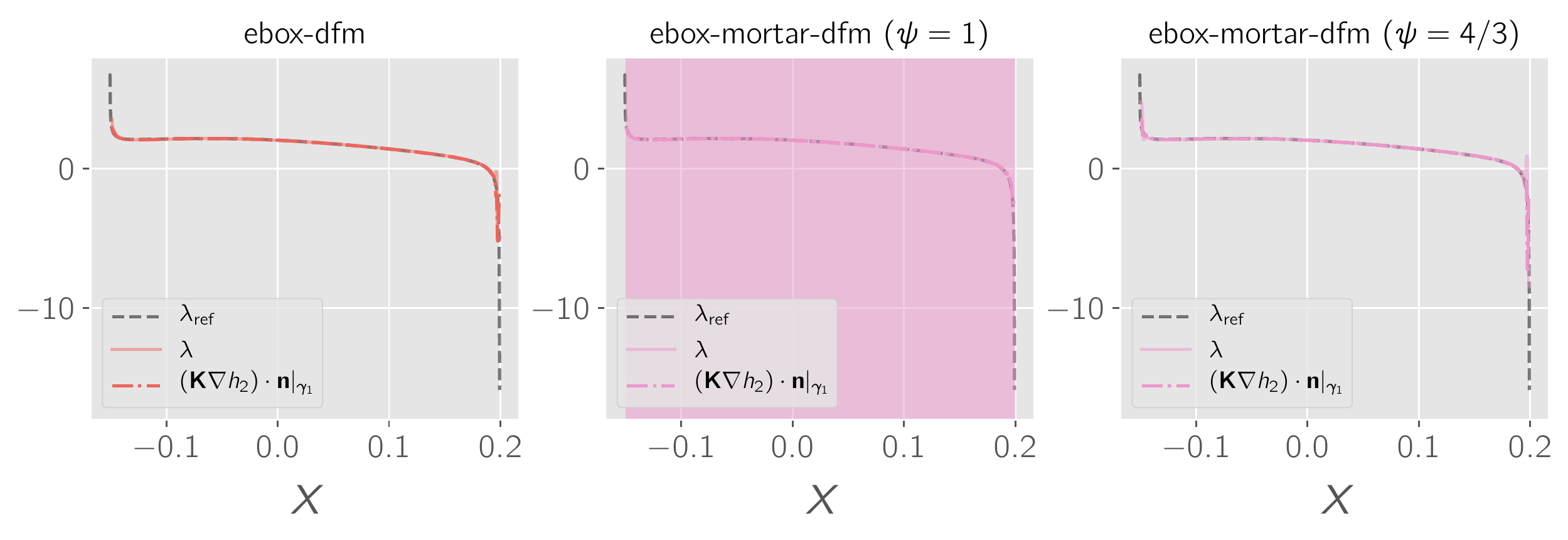}
    \end{subfigure}
  \begin{subfigure}{\textwidth}
    \centering
    \includegraphics[width=0.99\textwidth]{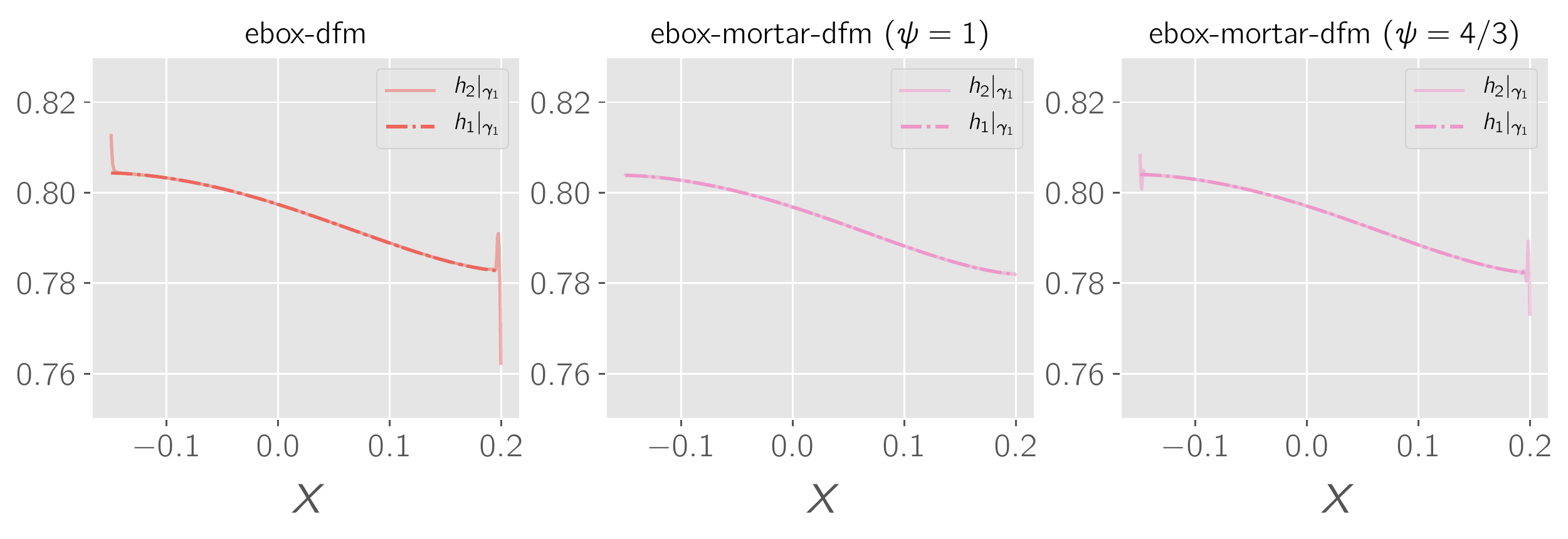}
  \end{subfigure}
  \caption{\textbf{Case 3.2 - interface flux and hydraulic head}. The top and bottom rows show the transfer fluxes and the hydraulic heads along the interface $\interface_1$, respectively.}
  \label{fig:iffluxcase_diri}
\end{figure}
Let us now turn to the case of the continuity conditions~\eqref{eq:prob_mixed_if_p}
being applied at the bulk-fracture interfaces, for which the results are depicted
in~\cref{fig:iffluxcase_diri}. As mentioned before, with these conditions the
\eboxMortarDfm scheme leads to singular system matrices for $\mortarGridFactor < 1$.
As can be seen in the plots, we consider the case of $\mortarGridFactor = 1$, for
which the transfer fluxes $\lambda$ exhibit spurious oscillations along the entire
interface. Since a rather fine grid is used, the oscillations appear as a shaded
area in the upper center plot of~\cref{fig:iffluxcase_diri}. Note that the
y-axis of the plot has been adjusted to match the non-oscillative solution, and the
actual magnitude of the oscillations is much larger. It is noteworthy that
despite the oscillations in $\lambda$, the term
$\left( \perm \grad \head_2 \right) |_{\interface_1}$ is still monotone and fits
well to the reference solution. The same holds for the hydraulic heads which are
shown in the plot directly below.
For a corser mortar grid, using $\mortarGridFactor = 4/3$, we again observe small
oscillations right at the fracture tips of both $\lambda$ and $\head_2$, as can
be seen in the plots in the right column of~\cref{fig:iffluxcase_diri}. The magnitude
of the oscillations seems to be smaller when compared to the case of the
conditions~\eqref{eq:prob_mixed_if_flux} (see~\cref{fig:iffluxcase_neum}).\\

As discussed earlier, the \eboxDfm scheme weakly incorporates the continuity
condition~\eqref{eq:prob_mixed_if_p}, and it can be seen from the lower left
plot of~\cref{fig:iffluxcase_diri} that the scheme is able to capture the continuity
along the entire interface, except for the fracture tips, where we now observe
noticeable differences between $\head_2 |_{\interface_1}$ and $\head_1 |_{\interface_1}$.
Recall that, when the conditions~\eqref{eq:prob_mixed_if_p} are used, the bulk-fracture
transfer fluxes are assembled independent of the interface head
$\head_2 |_{\interface_1}$ after~\cref{eq:discFluxCouplCond}. Therefore, these
deviations do not lead to a deterioration of $\lambda$, which is shown in the
upper left plot of~\cref{fig:iffluxcase_diri}. A very good match with the reference
solution is again obtained, and it seems that the small oscillations of $\lambda$
at the fracture tips are even less pronounced than for
$\faceSet_{d, \couplIdx} \equiv \faceSet_{d, \blockInterfaceSet}$
(see~\cref{fig:iffluxcase_neum}).

\subsection{Two-dimensional benchmark case}
\label{sec:benchmark2d}

\begin{figure}[ht]
  \begin{subfigure}{0.499\textwidth}
    \flushright
    \includegraphics[width=0.9\textwidth]{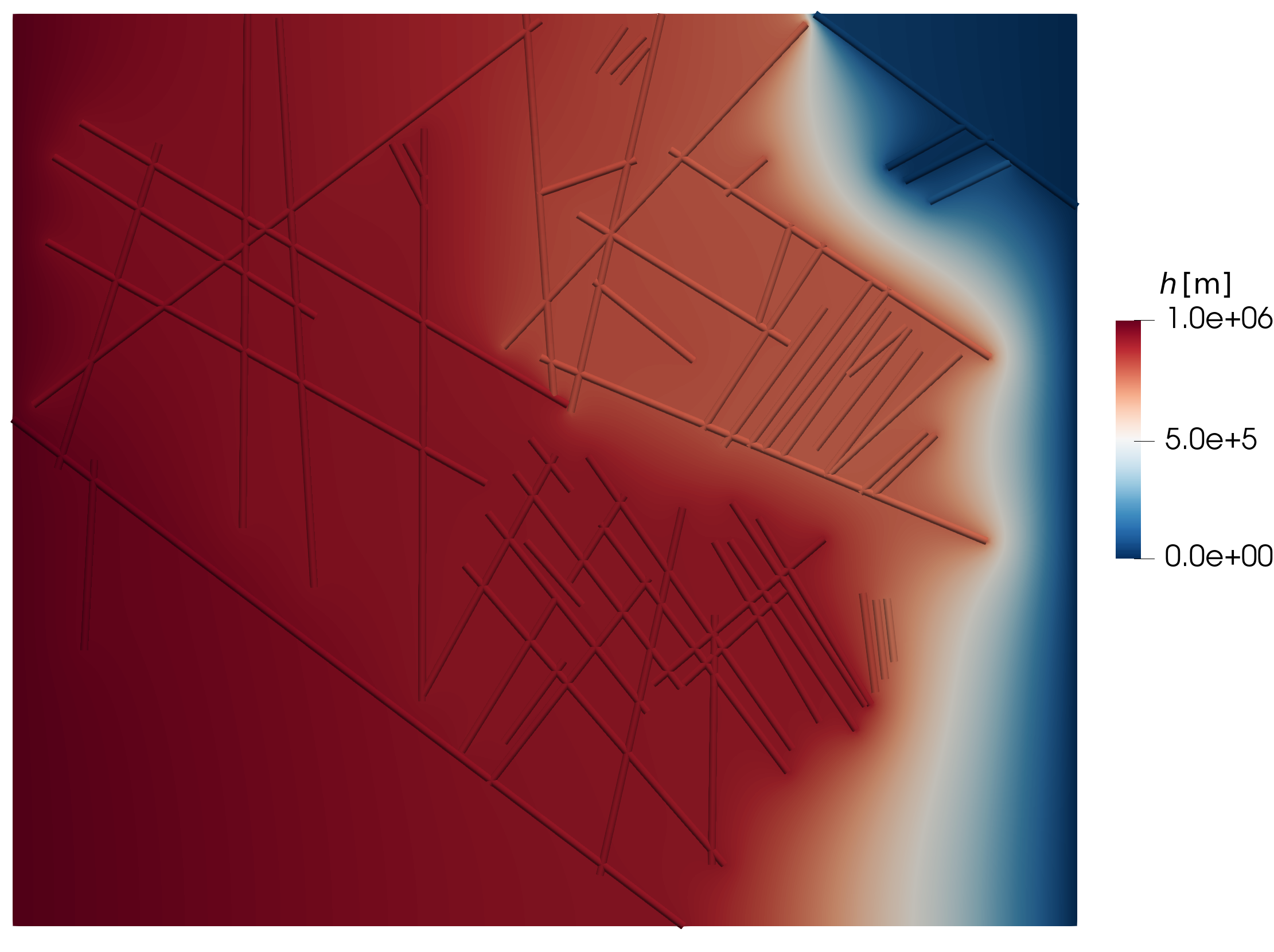}
    \caption{$k = \num{1e-8}$}
    \label{fig:benchmark_2d_conduit_sol}
    \end{subfigure}
  \begin{subfigure}{0.499\textwidth}
    \flushright
    \includegraphics[width=0.9\textwidth]{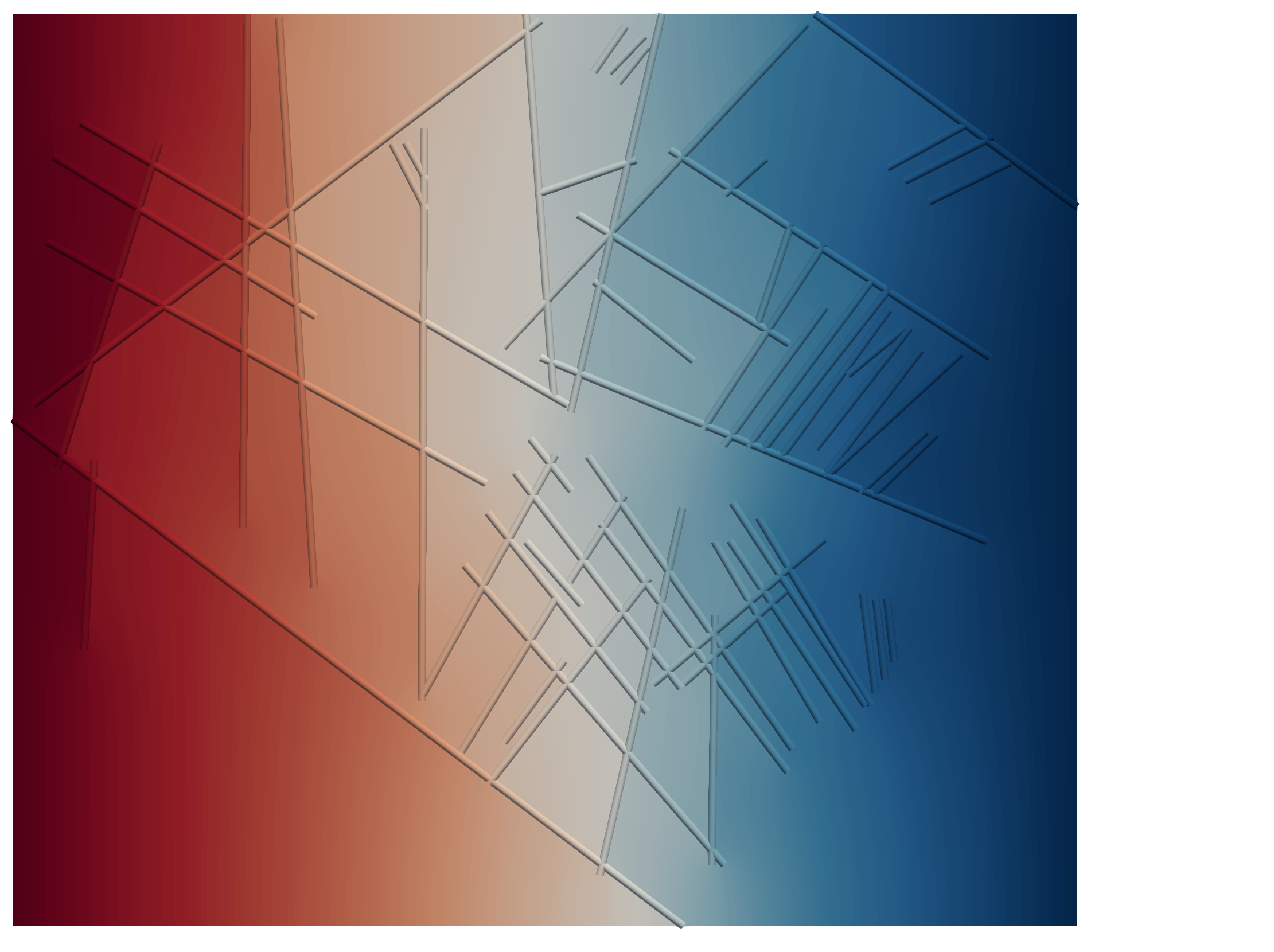}
    \caption{$k = \num{1e-18}$}
    \label{fig:benchmark_2d_barrier_sol}
  \end{subfigure}
  \caption{\textbf{Case 4 - numerical solutions}. Shown are the solutions obtained with the \eboxDfm scheme for the case of highly-permeable fractures (a) and low-permeable fractures (b). The fractures are visualized by tubes, colored with the solution for $\head_1$ in the fractures.}
  \label{fig:benchmark_2d_sol}
\end{figure}

\begin{figure}[ht]
  \begin{subfigure}{0.499\textwidth}
    \centering
    \includegraphics[width=0.99\textwidth]{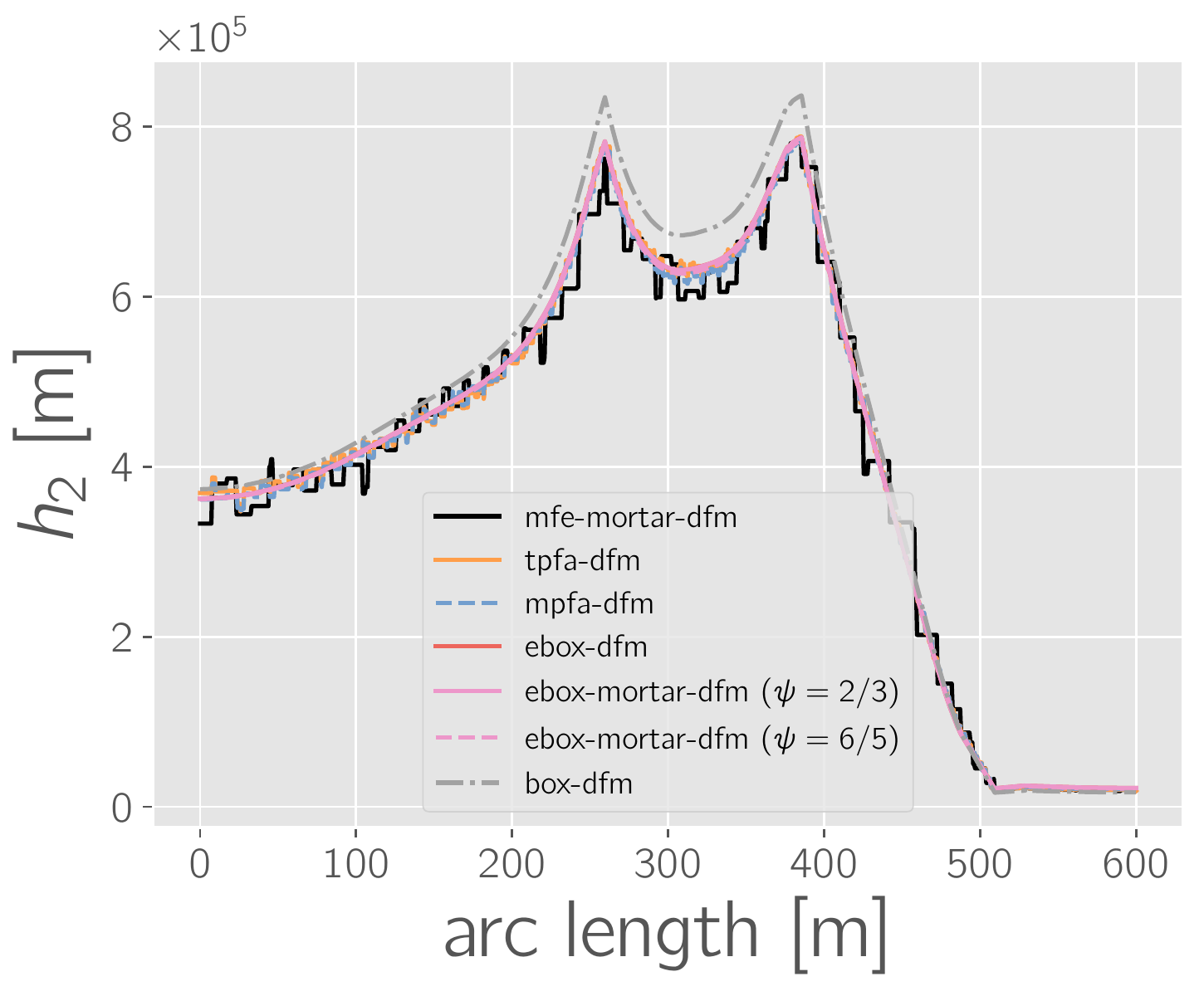}
    \caption{$k = \num{1e-8}$}
    \label{fig:benchmark_2d_conduit_plot}
    \end{subfigure}
  \begin{subfigure}{0.499\textwidth}
    \centering
    \includegraphics[width=0.99\textwidth]{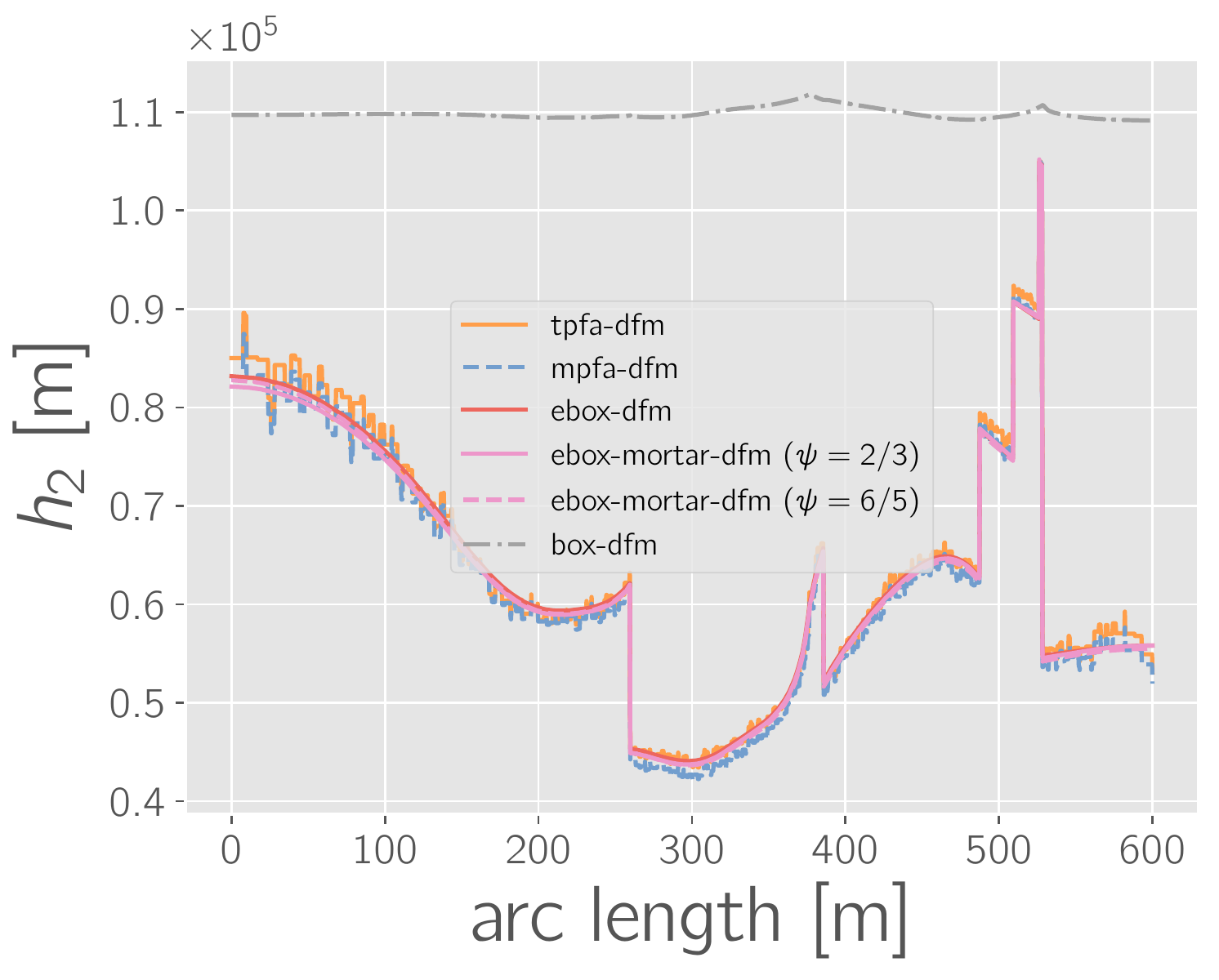}
    \caption{$k = \num{1e-18}$}
    \label{fig:benchmark_2d_barrier_plot}
  \end{subfigure}
  \caption{\textbf{Case 4 - plots over line}. The figures show the hydraulic head $\head_2$ in the bulk medium, plotted along the line $x = \SI{625}{\meter}$. The left plot further depicts the results obtained with the \textit{``Flux-Mortar''} scheme included in the original study \citep{Flemisch2018Benchmarks}, to which it is referred in the legend by \textit{``mfe-mortar-dfm''}.}
  \label{fig:benchmark_2d_plot}
\end{figure}

This test case is taken from \citet{Flemisch2018Benchmarks} and considers 64
fractures whose orientations were deduced from an image of an outcrop near Bergen,
Norway. The domain has a size of \SI{700x600}{\meter}, and a flow from left to
right is induced by prescribing $\head_2 = \head_1 = \SI{1013250}{\meter}$ at
$x = \SI{0}{\meter}$, and $\head_2 = \head_1 = \SI{0}{\meter}$ on the right
boundary at $x = \SI{700}{\meter}$. On the top and bottom boundaries, no-flow
boundary conditions are prescribed. The permeabilities are set to
$\perm_2 = \num{1e-14} \, \I \, \si{\meter\per\second}$ and
$\perm_1 = k \, \I \, \si{\meter\per\second}$, where $k = \num{1e-8}$, and an aperture of
$\aperture = \SI{1e-2}{\meter}$ is used. In this work, we consider a second
variant of the test case with low-permeable fractures using $k = \num{1e-18}$.
An illustration of the domain and the numerical solutions for
both fracture permeabilites can be found in~\cref{fig:benchmark_2d_sol}.\\
As in \citet{Flemisch2018Benchmarks}, we want to compare the numerical solutions
obtained with the different schemes based on plots over line of the hydraulic head.
\Cref{fig:benchmark_2d_plot} shows the hydraulic head plotted along the line
$x = \SI{625}{\meter}$ for both fracture permeabilities, where we added
the numerical solution of the \textit{``Flux-Mortar''} scheme presented in the
original study to the plot for $k = \num{1e-8}$
(The data was taken from the associated repository at
\href{https://git.iws.uni-stuttgart.de/benchmarks/fracture-flow/}{git.iws.uni-stuttgart.de/benchmarks/fracture-flow/}).
The scheme uses a mixed-finite element formulation with flux mortars at the
bulk-fracture interfaces, more details can be found in~\citet{Boon2018}.\\

\begin{table}[bh]
  \centering
  { \small
  \caption{\textbf{Case 4 - matrix characteristics}. Number of degrees of freedom ($N_\mathrm{dof}$) and number of nonzero entries ($N_\mathrm{nnz}$) in the system matrices.}
  \label{tab:benchmark_2d_matrices}
  \begin{tabular}{ *{8}{l} }
  \toprule
  & \tpfaDfm & \mpfaDfm & \eboxDfm & \multicolumn{3}{c}{\eboxMortarDfm} & \boxDfm \\
  &          &          &          & $\mortarGridFactor = 2/3$ & $\mortarGridFactor = 6/5$ & $\mortarGridFactor = 2$ \\ \midrule
  $N_\mathrm{dof}$ & \num{117438} & \num{117438} & \num{64403} & \num{75149} & \num{70473} & \num{57160} \\
  $N_\mathrm{nnz} \, (\faceSet_{d, \couplIdx} \equiv \faceSet_{d, \blockInterfaceSet})$ & \num{473426} & \num{1537488} & \num{476245} & \num{537252} & \num{500547} & - & - \\
  $N_\mathrm{nnz} \, (\faceSet_{d, \couplIdx} \equiv \faceSet_{d, \condInterfaceSet})$ & \num{473426} & \num{1537488} & \num{455750} & - & \num{494465} & \num{473491} & \num{398283} \\
  \bottomrule
  \end{tabular} }
\end{table}
\Cref{fig:benchmark_2d_conduit_plot} shows that the solutions obtained with the
different schemes in the case of highly-permeable fractures agree very
well among each other, except for the one of the \boxDfm scheme, which leads to
slightly higher hydraulic heads in the center part of the plot. This is in
agreement with the observations made in \citet{Flemisch2018Benchmarks}, where the
\boxDfm scheme, although a different implementation, was also included in the
investigations. Moreover, a good match with the results for the
\textit{``Flux-Mortar''} scheme of the original study can be seen. Note that in
that study,
a coarser mesh consisting of \num{7614} triangles and \num{867} line elements
was used, while in this work we use \num{113820} triangular and \num{3618} lines.
In the case of low-permeable fractures we see that the \boxDfm
scheme completely fails to capture relevant features of the solution, while the
results of the remaining schemes agree rather well
(see~\cref{fig:benchmark_2d_barrier_plot}).\\

On unstructured simplex grids, a benefit of vertex-centered schemes is that they
lead to significantly less degrees of freedom. The numbers for the grid used in
this test case are listed in~\cref{tab:benchmark_2d_matrices}, which shows
that approximately two times as many degrees of freedom are used in the cell-centered
schemes when compared to the \eboxDfm scheme.
The \eboxMortarDfm scheme introduces additional degrees of freedom due to discretization
of the mortar variable, and the \boxDfm scheme leads to the
smallest number of unknowns.
Despite the differences in number of unknowns, the number of nonzero entries in the
system matrix for the \tpfaDfm and the \eboxDfm scheme is very similar due to the
small stencils of the \tpfaDfm scheme. However, recall that the
\tpfaDfm scheme is inconsistent on unstructured grids and anisotropic permeabilities,
while \eboxDfm has shown to produce results that are comparable to the \mpfaDfm
scheme for different permeability angles and contrasts (see~\cref{sec:convDiscrete}).\\

Although leading to results of comparable quality, the \mpfaDfm scheme produces
about three times as many nonzero entries in the system matrix as the \eboxDfm
and the \eboxMortarDfm schemes (see~\cref{tab:benchmark_2d_matrices}). However,
we have tested several iterative linear solvers, and they failed to solve the system arising
from the \eboxMortarDfm scheme, which is why the results shown in~\cref{fig:benchmark_2d_plot}
were obtained with the direct solver UMFPack~\citep{Davis2004UmfPack}. Moreover,
for large fracture permeabilities, we have experienced iterative linear solver issues
also for the \eboxDfm scheme. This indicates poor conditioning of the system matrix,
which might be due to the fact that in the \eboxDfm scheme, the bulk-fracture coupling
entries scale with $k/a$ as a consequence of condition~\eqref{eq:prob_mixed_if_flux},
while the diagonal block scales with the bulk permeability. This could explain why we
observed a deterioration of the iterative linear solver performance with increasing fracture
permeability, while this effect was not observed for the cell-centered schemes,
where the coupling entries scale with the harmonic mean of the bulk and fracture
permeabilities. However, one possible remedy seems to be the use of the
conditions~\eqref{eq:prob_mixed_if_p} for highly-permeable fractures, for which
we have observed better iterative linear solver performance with the \eboxDfm scheme, while
visually indistinguishable results compared to those of~\cref{fig:benchmark_2d_plot}
were obtained. Please note that a detailed analysis of the performance of iterative
linear solvers is not within the scope of this work and should be addressed in future
investigations.

\subsection{Three-dimensional benchmark case}
\label{sec:benchmark3d}

\begin{figure}[ht]
    \centering
    \includegraphics[width=0.99\textwidth]{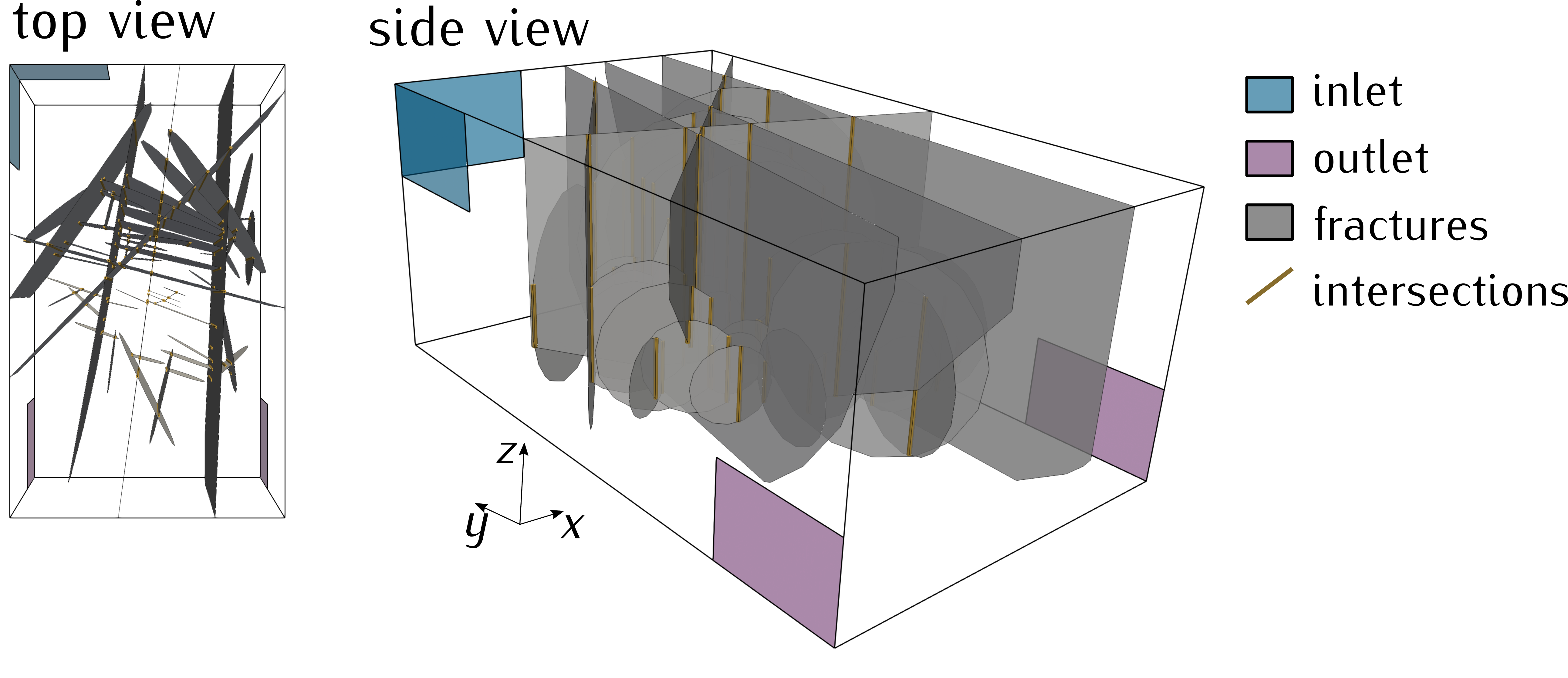}
    \caption{\textbf{Case 5 - computational domain}. The figure visualizes the fracture network considered in the three-dimensional benchmark case, which is taken from \citet{Berre2021Benchmarks}. Moreover, the inlet and outlet boundary segments are depicted.}
    \label{fig:benchmark_3d_domain}
\end{figure}
This benchmark case is taken from \citet{Berre2021Benchmarks}, and it considers
a domain with the dimensions
$\domain = \left( \SI{500}{\meter}, \SI{350}{\meter} \right)
           \times \left( \SI{100}{\meter}, \SI{1500}{\meter} \right)
           \times \left( \SI{-100}{\meter}, \SI{500}{\meter} \right)$,
containing \num{52} fractures with \num{106} intersections. The fracture network
was again created by postprocessing of an outcrop near Bergen in Norway, and was
originally presented in \citet{Fumagalli2019DFN}. \Cref{fig:benchmark_3d_domain} illustrates
the domain and fracture network, and it depicts the inlet and outlet
boundaries, where uniform unit inflow and a hydraulic head of
$\head_2 = \SI{0}{\meter}$ are set, respectively.
The permeabilities
are given by $\perm_2 = \I \, \si{\meter\per\second}$ and
$\perm_1 = \num{1e4} \, \I \, \si{\meter\per\second}$, and an aperture of
$\aperture = \SI{1e-2}{\meter}$ is used.\\

\begin{figure}[ht]
  \centering
  \includegraphics[width=0.8\textwidth]{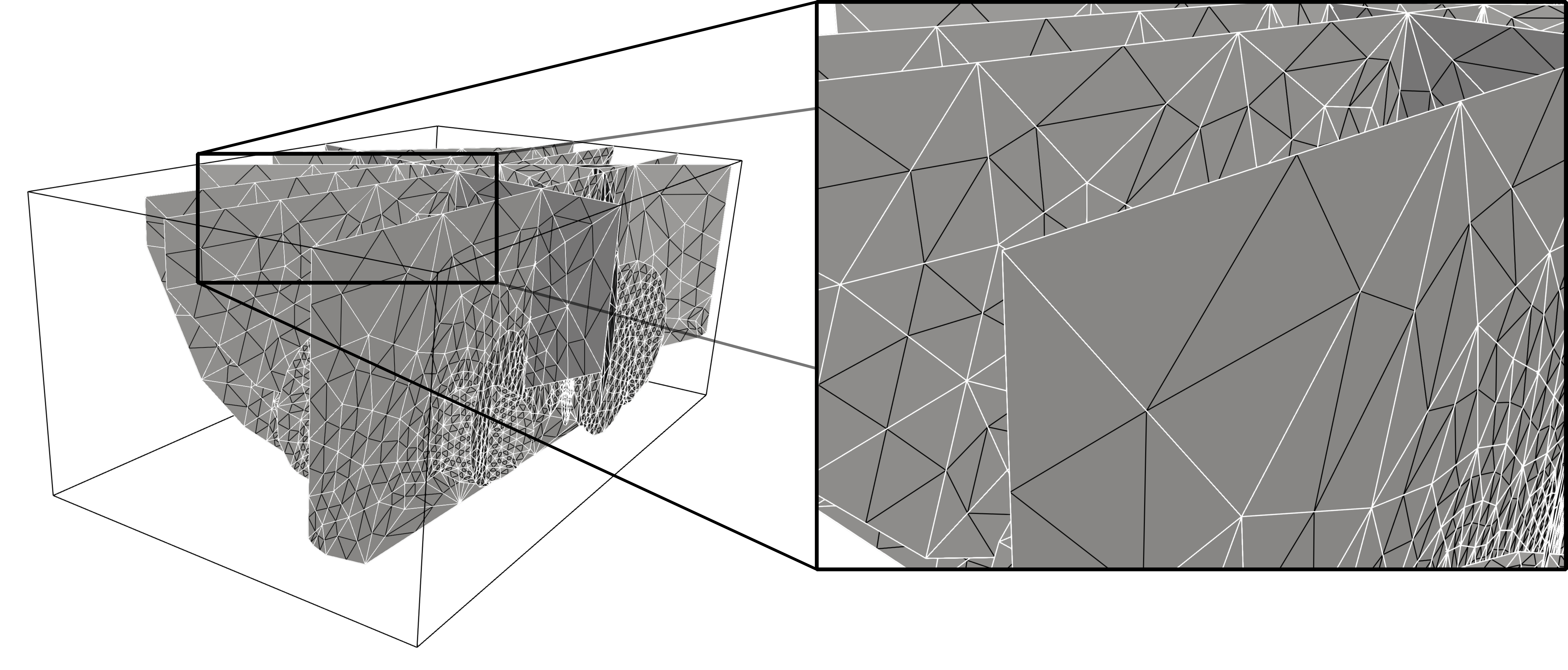}
  \caption{\textbf{Case 5 - fracture and mortar grid}. Closeup on the grids used in the fracture and the mortar domains. The black lines show the grid elements of the fracture grids, while the white lines illustrate those of the mortar grid.}
  \label{fig:benchmark_3d_grid}
\end{figure}

We use the conditions~\eqref{eq:prob_mixed_if_p} at all interfaces, which can be
justified by the small aperture and the comparatively large fracture permeability.
Comparisons of the solutions for both interface conditions
(obtained with the \tpfaDfm scheme)
showed that they lead to visually indistinguishable results.
By choosing the interface condition~\eqref{eq:prob_mixed_if_p}, we must choose
$\mortarGridFactor > 1$, and here we use $\mortarGridFactor = 4$, which means
that each triangle of the discretization
of the mortar domain contains exactly four triangles of the fracture grid.
An illustration of this is given in~\cref{fig:benchmark_3d_grid}. Note that the
implementation is able to handle nonconforming meshes, however, for this benchmark
we have experienced issues with adjusting the geometric tolerance used upon mesh
generation with Gmsh~\citep{Gmsh2009} to that used during mesh intersection in Dumux.
Small differences
in the geometry due to inaccuracies during mesh generation might appear when
generating the meshes used for flow and the mortar domain, respectively.
These in turn might result
in extremely small intersection segments when the used geometric tolerance does not match.
In order avoid such issues, we have generated the mesh used for flow by subdivision
of a coarser grid, which has been used for the discretization of the
mortar domain. Ongoing work adresses the improvement of the robustness of the
intersection algorithms and the choice of optimal geometric tolerances.\\

\begin{figure}[ht]
  \begin{subfigure}{0.499\textwidth}
    \centering
    \includegraphics[width=0.99\textwidth]{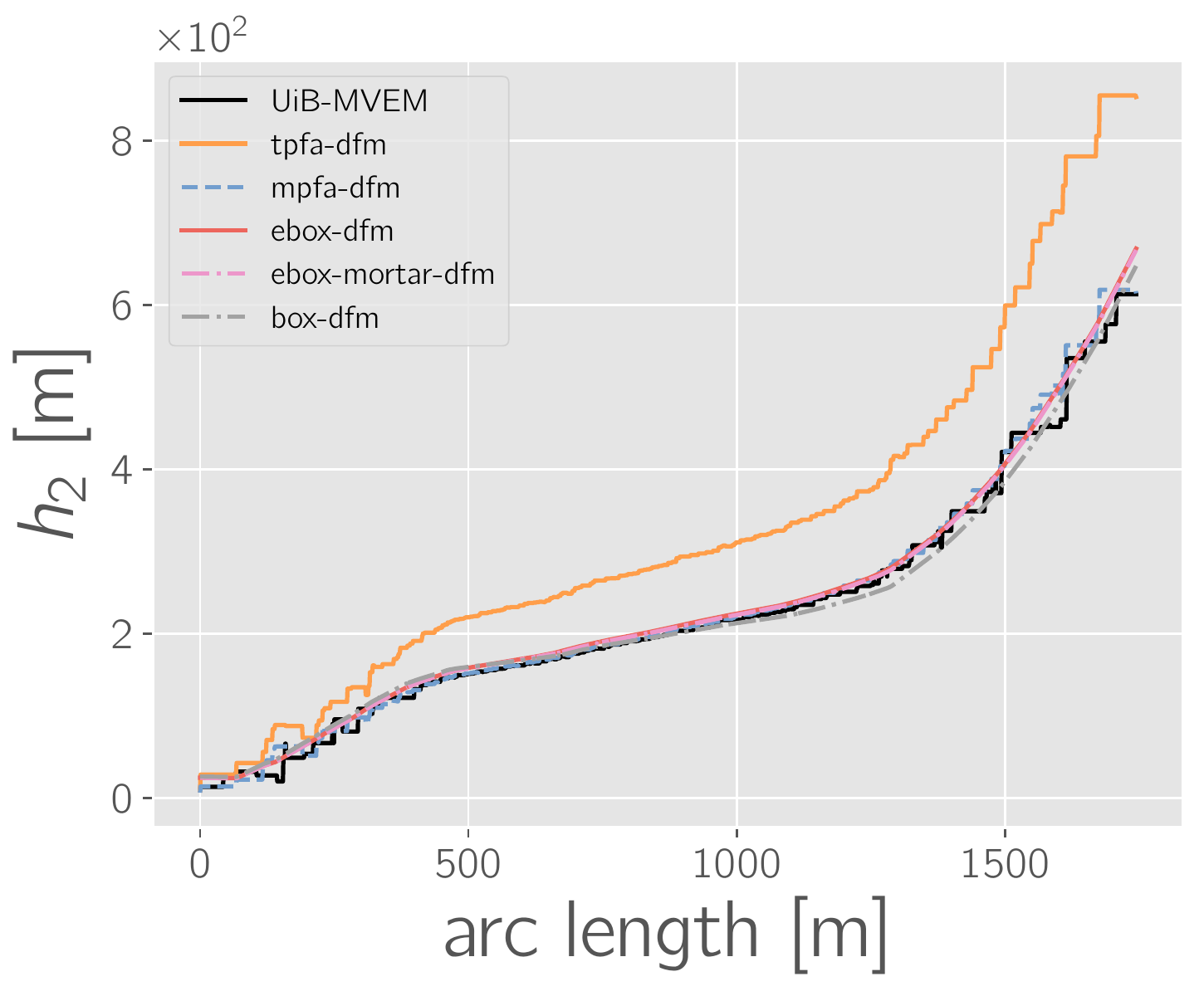}
    \caption{$\left(350, 100, -100 \right) - \left(-500, 1500, 500 \right)$}
    \label{fig:benchmark_3d_plot_1}
    \end{subfigure}
  \begin{subfigure}{0.499\textwidth}
    \centering
    \includegraphics[width=0.99\textwidth]{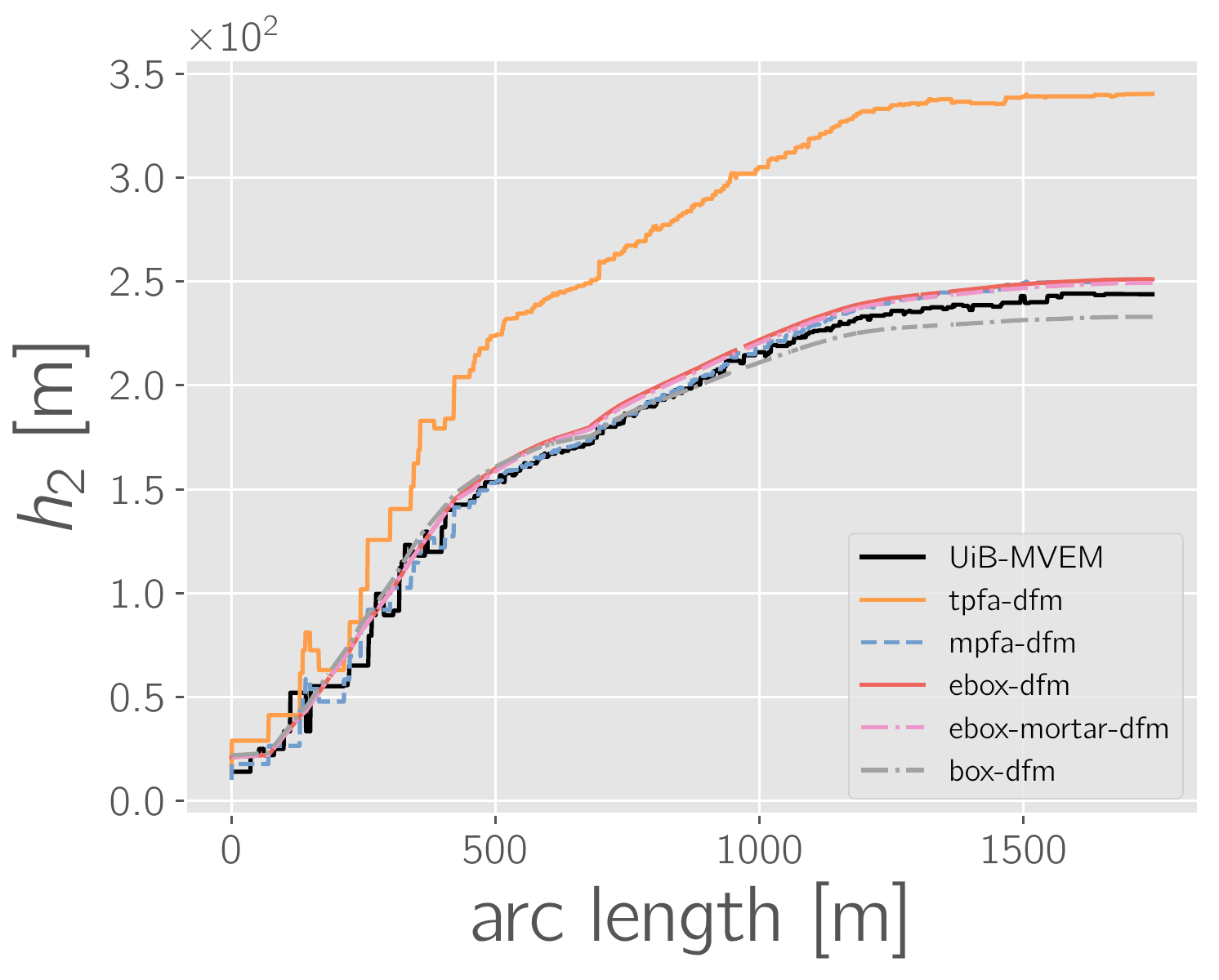}
    \caption{$\left(-500, 100, -100 \right) - \left(350, 1500, 500 \right)$}
    \label{fig:benchmark_3d_plot_2}
  \end{subfigure}
  \caption{\textbf{Case 5 - plots over line}. The figures show the hydraulic head $\head_2$ in the bulk medium, plotted along two different lines as given in the captions. The results shown here correspond to those simulations in which flow along fracture intersections was neglected. Moreover, the plots include the results obtained with the \textit{``UiB-MVEM''} scheme presented in the benchmark study \citet{Berre2021Benchmarks}.}
  \label{fig:benchmark_3d_plot}
\end{figure}

We simulate the benchmark with all schemes considered in this work, and we perform
simulations both taking into account flow along fracture intersections and
assuming continuity after the conditions~\eqref{eq:prob_mixed_junctions}.
An exception is the \boxDfm scheme, for which the current implementation
in \dumux does not allow for flow along intersections.
\Cref{fig:benchmark_3d_plot} shows plots of the hydraulic head $\head_2$
along two different lines, both of which are used as metrics of comparison in
the original study \citet{Berre2021Benchmarks}. The depicted results correspond to the simulations
neglecting intersection flow, however, visually
indistinguishable plots were obtained when taking it into account, which is
why we do not present them here. In this benchmark, the fracture intersections
are all oriented vertically, and thus, orthogonal to the direction of flow. For
other orientations or for larger permeability contrasts, taking into account flow
along intersections might have a significant impact on the overall
flow field (see, for example, \citet{glaeser_2020}).
The plots further show the results of the \textit{``UiB-MVEM''} scheme, which
uses a formulation on the basis of mixed virtual elements and is
one of the schemes considered in the original study \citet{Berre2021Benchmarks}.
For further details on the scheme we refer to~\cite{Fumagalli2019dual}.
The results were taken from the data repository associated with the publication
(\href{https://git.iws.uni-stuttgart.de/benchmarks/fracture-flow-3d}{git.iws.uni-stuttgart.de/benchmarks/fracture-flow-3d}),
and they were originally produced with an implementation provided in the open-source numerical
simulation tool PorePy~\citep{keilegavlen2020porepy}.\\

In the first plot (\cref{fig:benchmark_3d_plot_1}), we can see that a good
match of the hydraulic head was obtained among the different schemes, and that
they fit well to the reference curve of the original study. However, significantly
larger hydraulic heads are obtained with the \tpfaDfm scheme. Note that the
\tpfaDfm scheme, using the implementation provided in PorePy,
was also taken into account in~\citet{Berre2021Benchmarks}, where the respective curve was much closer to
those of the other schemes. But, here we are using a different grid, which is constructed
by refinement of an initial discretization in order to obtain the desired refinement ratio
of $\mortarGridFactor = 4$ between the grids used for flow and the mortar grid.
The properties of the resulting mesh might again promote the inconsistency of the
scheme, as was the case also in~\cref{sec:convDiscrete}.\\

A rather good match is obtained with the \boxDfm scheme, which leads to only
slightly lower hydraulic heads. However, the deviations are a bit more pronounced
in the second plot shown in~\cref{fig:benchmark_3d_plot_2}, where we again observe
that the schemes presented in this work lead to very similar results,
although slightly larger
deviations to the reference curve are visible. But as mentioned, we do
not use the same grid as in the original study and the deviations might decrease
upon further grid refinement.\\

A comparison of the computational cost of the different schemes is given
in~\cref{tab:benchmark3dMatrices}. It is noteworthy that on the three-dimensional
unstructured grid used here, the ratio of degrees of freedom for cell-centered and
vertex-centered schemes is significantly higher than in the two-dimensional case
(see~\cref{sec:benchmark2d}).
Besides that, the \mpfaDfm scheme leads to substantially more nonzero entries
in the linear system matrix, when compared to the other schemes, as a result of
its large stencils on unstructured simplex grids. In particular, the system matrix
of the \mpfaDfm scheme contains about 15 times more nonzero entries than those
of the \eboxDfm and \eboxMortarDfm schemes, which has strong implications on the
computational efficiency of the schemes. Despite these differences, they
produce very similar results (see~\cref{fig:benchmark_3d_plot}).


\begin{table}[h]
  \centering
  \caption{\textbf{Case 5 - matrix characteristics}. Number of degrees of freedom ($N_\mathrm{dof}$) and number of nonzero entries ($N_\mathrm{nnz}$) in the system matrices of the different schemes.
  Results are shown for the case of considering ($\domain_1 \neq \emptyset$) and neglecting ($\domain_1 = \emptyset$) flow along intersections.}
  \label{tab:benchmark3dMatrices}
  \begin{tabular}{ *{6}{l} }
  \toprule
  & \tpfaDfm & \mpfaDfm & \eboxDfm & \eboxMortarDfm & \boxDfm \\ \midrule
  \multicolumn{6}{c}{$\domain_1 = \emptyset$} \\
$N_\mathrm{dof}$ & \num{320004} & \num{320004} & \num{86409} & \num{105835} & \num{47700} \\
$N_\mathrm{nnz}$ & \num{1640238} & \num{21639021} & \num{1390232} & \num{1492260} & \num{707930} \\
\midrule
\multicolumn{6}{c}{$\domain_1 \neq \emptyset$} \\
$N_\mathrm{dof}$ & \num{321028} & \num{321028} & \num{90315} & \num{111603} \\
$N_\mathrm{nnz}$ & \num{1640490} & \num{21677870} & \num{1433077} & \num{1527509} \\
  \bottomrule
  \end{tabular}
\end{table}

\section{Summary and outlook}
\label{sec:conclusions}

In this work we have presented a comparison, by means of numerical experiments,
of different finite volume schemes for the simulation of flow
in fractured porous media. The underlying mathematical model is based on a mixed-dimensional
formulation in which the bulk medium, the fractures and intersections of fractures
are described by geometries of decreasing dimensionality. On each of these
subdomains, an individual system of partial differential equations is solved and
the interaction between them is described by fluxes across interior
boundaries and via source terms on the lower-dimensional features.\\

Several well-established schemes were considered, namely
the cell-centered \tpfaDfm and \mpfaDfm schemes
\citep{Karimi2004TpfaDfm,Sandve2012MpfaDfm,Ahmed2015MpfaDfm}
and the vertex-centered \boxDfm scheme \citep{Reichenberger2006mixed}.
However, the mathematical formulation of the latter differs from the one
presented in this work in several aspects, as it assumes continuity
of the hydraulic head across the fractures and does not explicitly
describe the bulk-fracture mass exchange. Therefore, we have introduced
an extension of the latter in this work, the \eboxDfm scheme,
which is able to capture the discontinuities in the hydraulic head and the
fluxes that potentially occur across the fractures. In this scheme,
the bulk-fracture transfer fluxes are expressed in terms of the hydraulic heads
in the two subdomains and the coupling conditions are incorporated weakly.
We have furthermore introduced the \eboxMortarDfm scheme, in which
the transfer fluxes appear as additional unknowns, allowing for a strong
incorporation of the desired coupling conditions at subdomain interfaces.\\

We have studied the performance of the different schemes on various two- and
three-dimensional test cases. The convergence behaviour
against analytical and equi-dimensional reference solutions was investigated,
and we have applied the schemes to two- and three-dimensional benchmark cases
taken from the literature. Optimal convergence was observed for the
newly-introduced schemes and we have found a good agreement of the numerical
results to those presented in the original benchmark studies. Moreover, we
observed that the solutions
obtained with \eboxDfm and \eboxMortarDfm schemes were comparable to those of
the \mpfaDfm scheme while being computationally much more efficient. In general,
cell-centered schemes lead to many more degrees of freedom than vertex-centered
schemes on unstructured simplex grids. Furthermore, the \mpfaDfm scheme involves
very large stencils, especially on three-dimensional grids, which leads to a large
number of non-zero entries in the system matrices. In this regard, the presented
vertex-centered schemes show to be a promising alternative, producing consistent
results with less computational effort.\\

However, we have observed that several iterative linear solvers struggled to solve
the systems arising from the \eboxDfm and \eboxMortarDfm schemes, especially for
large fracture permeabilities. Therefore, we have used a direct solver in this work.
Future investigations could be targetted at developing iterative linear solvers
that are well-suited for solving the systems arising from these schemes.
Moreover, the performance of the schemes in the context of two-phase flow could be
studied in future work, which is of high relevance in a variety of geotechnical applications.
In multiphasic settings, we expect that the importance of the numerical schemes
being able to capture discontinuities across the fractures is even more pronounced,
as jumps in saturations and capillary pressures occur at material interfaces.
While the cell-centered schemes and the \boxDfm scheme have already been applied
to two-phase flow (see \eg \citet{Reichenberger2006mixed,Glaeser2017MpfaDfmTwoP}),
it might be more difficult to derive a formulation for the \eboxMortarDfm scheme.
Two-phase flow models typically involve non-linearities, and therefore,
it would be of particular interest to study the convergence behaviour of non-linear
solvers for the different numerical schemes.

\section*{Acknowledgements}
We thank the Deutsche Forschungsgemeinschaft (DFG, German Research Foundation)
for supporting this work by funding SFB 1313, Project Number 327154368.

\bibliography{dfm_fv_comparison.bib}

\newpage
\section*{Appendix}
\begin{table}[h]
  {\footnotesize
  \centering
  \caption{\textbf{Case 2 - errors and rates of $\head_2$ for $\permAngle = 0$}. Listed are the errors $\errorNorm_{\head_2}$ and the corresponding rates $r_{\head_2}$ over grid refinement, expressed in terms of $\discLength_\fracIdx$.}
  \begin{tabular}{ *{2}{l} | *{2}{l} | *{2}{l} | *{2}{l} | *{2}{l} | *{2}{l} }
  \toprule
  & & \multicolumn{2}{l}{\tpfaDfm} & \multicolumn{2}{l}{\mpfaDfm} & \multicolumn{2}{l}{\eboxDfm} & \multicolumn{2}{l}{\eboxMortarDfm} & \multicolumn{2}{l}{\boxDfm} \\
  &  $\discLength_\fracIdx$ & $\errorNorm_{\head_2}$ & $r_{\head_2}$ & $\errorNorm_{\head_2}$ & $r_{\head_2}$ & $\errorNorm_{\head_2}$ & $r_{\head_2}$ & $\errorNorm_{\head_2}$ & $r_{\head_2}$ & $\errorNorm_{\head_2}$ & $r_{\head_2}$ \\ \midrule
   \multicolumn{12}{c}{$k = \num{1e4}$} \\ \midrule
\multirow{6}{*}{$\fracNet_1$} & 2.43e-02 & 2.12e-01 &  & 1.40e-01 &  & 1.55e-01 &  & 1.52e-01 &  & 1.55e-01 &  \\
 & 1.21e-02 & 1.64e-01 & 3.67e-01 & 6.78e-02 & 1.04e+00 & 7.25e-02 & 1.10e+00 & 7.05e-02 & 1.11e+00 & 7.26e-02 & 1.10e+00 \\
 & 6.07e-03 & 1.44e-01 & 1.94e-01 & 3.22e-02 & 1.07e+00 & 3.51e-02 & 1.04e+00 & 3.40e-02 & 1.05e+00 & 3.53e-02 & 1.04e+00 \\
 & 3.04e-03 & 1.32e-01 & 1.19e-01 & 1.47e-02 & 1.13e+00 & 1.65e-02 & 1.09e+00 & 1.58e-02 & 1.10e+00 & 1.67e-02 & 1.08e+00 \\
 & 1.52e-03 & 1.36e-01 & -3.73e-02 & 6.26e-03 & 1.23e+00 & 7.34e-03 & 1.17e+00 & 6.91e-03 & 1.20e+00 & 7.56e-03 & 1.15e+00 \\
 & 7.59e-04 & 1.33e-01 & 3.35e-02 & 2.95e-03 & 1.08e+00 & 3.49e-03 & 1.07e+00 & 3.26e-03 & 1.08e+00 & 3.73e-03 & 1.02e+00 \\
 \midrule
\multirow{6}{*}{$\fracNet_2$} & 2.25e-02 & 2.42e-01 &  & 1.53e-01 &  & 2.01e-01 &  & 1.86e-01 &  & 2.01e-01 &  \\
 & 1.12e-02 & 1.72e-01 & 4.95e-01 & 7.82e-02 & 9.71e-01 & 9.71e-02 & 1.05e+00 & 9.10e-02 & 1.03e+00 & 9.75e-02 & 1.04e+00 \\
 & 5.62e-03 & 5.40e-01 & -1.65e+00 & 3.93e-02 & 9.93e-01 & 4.87e-02 & 9.95e-01 & 4.56e-02 & 9.98e-01 & 4.92e-02 & 9.88e-01 \\
 & 2.81e-03 & 1.33e-01 & 2.03e+00 & 1.93e-02 & 1.02e+00 & 2.41e-02 & 1.02e+00 & 2.24e-02 & 1.03e+00 & 2.45e-02 & 1.01e+00 \\
 & 1.41e-03 & 4.94e-01 & -1.90e+00 & 9.17e-03 & 1.07e+00 & 1.15e-02 & 1.07e+00 & 1.06e-02 & 1.08e+00 & 1.19e-02 & 1.04e+00 \\
 & 7.03e-04 & 2.45e-01 & 1.01e+00 & 4.35e-03 & 1.08e+00 & 5.37e-03 & 1.10e+00 & 4.89e-03 & 1.12e+00 & 5.81e-03 & 1.04e+00 \\
 \midrule
\multirow{6}{*}{$\fracNet_3$} & 1.70e-02 & 2.16e-01 &  & 1.33e-01 &  & 1.79e-01 &  & 1.65e-01 &  & 1.75e-01 &  \\
 & 8.52e-03 & 1.76e-01 & 3.01e-01 & 6.69e-02 & 9.86e-01 & 8.63e-02 & 1.05e+00 & 7.88e-02 & 1.07e+00 & 8.46e-02 & 1.05e+00 \\
 & 4.26e-03 & 1.50e-01 & 2.29e-01 & 3.28e-02 & 1.03e+00 & 4.63e-02 & 8.98e-01 & 4.20e-02 & 9.10e-01 & 4.67e-02 & 8.57e-01 \\
 & 2.13e-03 & 2.40e-01 & -6.79e-01 & 1.62e-02 & 1.02e+00 & 2.60e-02 & 8.33e-01 & 2.36e-02 & 8.33e-01 & 2.81e-02 & 7.31e-01 \\
 & 1.06e-03 & 1.09e+00 & -2.19e+00 & 8.93e-03 & 8.60e-01 & 1.54e-02 & 7.58e-01 & 1.41e-02 & 7.39e-01 & 1.90e-02 & 5.65e-01 \\
 & 5.32e-04 & 1.52e-01 & 2.84e+00 & 6.69e-03 & 4.17e-01 & 1.01e-02 & 6.03e-01 & 9.56e-03 & 5.62e-01 & 1.49e-02 & 3.57e-01 \\ \midrule
 \multicolumn{12}{c}{$k = \num{1e-4}$} \\ \midrule
 \multirow{6}{*}{$\fracNet_1$} & 2.43e-02 & 2.35e-01 &  & 1.58e-01 &  & 1.45e-01 &  & 1.45e-01 &  & 1.96e-01 &  \\
 & 1.21e-02 & 1.84e-01 & 3.50e-01 & 8.29e-02 & 9.28e-01 & 6.66e-02 & 1.12e+00 & 6.66e-02 & 1.12e+00 & 1.50e-01 & 3.81e-01 \\
 & 6.07e-03 & 1.60e-01 & 2.02e-01 & 4.20e-02 & 9.79e-01 & 3.14e-02 & 1.08e+00 & 3.14e-02 & 1.08e+00 & 1.40e-01 & 1.09e-01 \\
 & 3.04e-03 & 1.44e-01 & 1.52e-01 & 2.09e-02 & 1.01e+00 & 1.43e-02 & 1.14e+00 & 1.43e-02 & 1.14e+00 & 1.37e-01 & 2.59e-02 \\
 & 1.52e-03 & 1.35e-01 & 9.82e-02 & 1.08e-02 & 9.61e-01 & 6.11e-03 & 1.22e+00 & 6.11e-03 & 1.22e+00 & 1.37e-01 & 5.00e-03 \\
 & 7.59e-04 & 1.28e-01 & 7.88e-02 & 6.48e-03 & 7.31e-01 & 3.47e-03 & 8.17e-01 & 3.47e-03 & 8.17e-01 & 1.37e-01 & 4.65e-04 \\
 \midrule
\multirow{6}{*}{$\fracNet_2$} & 2.25e-02 & 2.68e-01 &  & 2.06e-01 &  & 1.49e-01 &  & 1.49e-01 &  & 2.13e-01 &  \\
 & 1.12e-02 & 2.41e-01 & 1.55e-01 & 1.05e-01 & 9.73e-01 & 7.20e-02 & 1.05e+00 & 7.19e-02 & 1.05e+00 & 1.71e-01 & 3.22e-01 \\
 & 5.62e-03 & 2.10e-01 & 1.96e-01 & 5.15e-02 & 1.03e+00 & 3.51e-02 & 1.04e+00 & 3.51e-02 & 1.04e+00 & 1.59e-01 & 1.02e-01 \\
 & 2.81e-03 & 2.12e-01 & -1.15e-02 & 2.55e-02 & 1.01e+00 & 1.67e-02 & 1.07e+00 & 1.67e-02 & 1.07e+00 & 1.56e-01 & 2.79e-02 \\
 & 1.41e-03 & 8.99e-01 & -2.09e+00 & 1.31e-02 & 9.65e-01 & 7.62e-03 & 1.13e+00 & 7.62e-03 & 1.13e+00 & 1.55e-01 & 7.37e-03 \\
 & 7.03e-04 & 2.04e-01 & 2.14e+00 & 7.39e-03 & 8.23e-01 & 3.86e-03 & 9.81e-01 & 3.86e-03 & 9.81e-01 & 1.55e-01 & 1.97e-03 \\
 \midrule
\multirow{6}{*}{$\fracNet_3$} & 1.70e-02 & 1.52e-01 &  & 1.29e-01 &  & 1.19e-01 &  & 1.19e-01 &  & 2.43e-01 &  \\
 & 8.52e-03 & 1.54e-01 & -2.27e-02 & 6.22e-02 & 1.05e+00 & 4.99e-02 & 1.25e+00 & 4.99e-02 & 1.25e+00 & 2.13e-01 & 1.87e-01 \\
 & 4.26e-03 & 1.13e-01 & 4.48e-01 & 2.90e-02 & 1.10e+00 & 2.35e-02 & 1.09e+00 & 2.35e-02 & 1.09e+00 & 2.07e-01 & 4.40e-02 \\
 & 2.13e-03 & 1.12e-01 & 7.10e-03 & 1.48e-02 & 9.75e-01 & 1.21e-02 & 9.57e-01 & 1.21e-02 & 9.56e-01 & 2.05e-01 & 1.13e-02 \\
 & 1.06e-03 & 1.24e-01 & -1.36e-01 & 9.47e-03 & 6.41e-01 & 8.17e-03 & 5.65e-01 & 8.17e-03 & 5.66e-01 & 2.05e-01 & 3.12e-03 \\
 & 5.32e-04 & 5.01e-01 & -2.02e+00 & 8.04e-03 & 2.36e-01 & 7.48e-03 & 1.28e-01 & 7.47e-03 & 1.28e-01 & 2.05e-01 & 9.23e-04 \\

  \bottomrule
  \end{tabular}
  \label{tab:convDiscreteErrors_head2_angle0}
  }
\end{table}

\begin{table}[h]
  {\footnotesize
  \centering
  \caption{\textbf{Case 2 - errors and rates of $\head_1$ for $\permAngle = 0$}. Listed are the errors $\errorNorm_{\head_1}$ and the corresponding rates $r_{\head_1}$ over grid refinement, expressed in terms of $\discLength_\fracIdx$.}
  \begin{tabular}{ *{2}{l} | *{2}{l} | *{2}{l} | *{2}{l} | *{2}{l} | *{2}{l} }
  \toprule
  & & \multicolumn{2}{l}{\tpfaDfm} & \multicolumn{2}{l}{\mpfaDfm} & \multicolumn{2}{l}{\eboxDfm} & \multicolumn{2}{l}{\eboxMortarDfm} & \multicolumn{2}{l}{\boxDfm} \\
  &  $\discLength_\fracIdx$ & $\errorNorm_{\head_1}$ & $r_{\head_1}$ & $\errorNorm_{\head_1}$ & $r_{\head_1}$ & $\errorNorm_{\head_1}$ & $r_{\head_1}$ & $\errorNorm_{\head_1}$ & $r_{\head_1}$ & $\errorNorm_{\head_1}$ & $r_{\head_1}$ \\ \midrule
  \multicolumn{12}{c}{$k = \num{1e4}$} \\ \midrule
\multirow{6}{*}{$\fracNet_1$} & 2.43e-02 & 1.15e-01 &  & 3.75e-02 &  & 6.51e-02 &  & 6.30e-02 &  & 6.02e-02 &  \\
 & 1.21e-02 & 8.54e-02 & 4.32e-01 & 1.47e-02 & 1.35e+00 & 2.09e-02 & 1.64e+00 & 2.00e-02 & 1.65e+00 & 1.80e-02 & 1.74e+00 \\
 & 6.07e-03 & 7.12e-02 & 2.64e-01 & 5.22e-03 & 1.49e+00 & 7.01e-03 & 1.58e+00 & 6.64e-03 & 1.59e+00 & 5.77e-03 & 1.64e+00 \\
 & 3.04e-03 & 6.53e-02 & 1.24e-01 & 3.30e-03 & 6.63e-01 & 3.45e-03 & 1.02e+00 & 3.36e-03 & 9.81e-01 & 4.19e-03 & 4.62e-01 \\
 & 1.52e-03 & 5.39e-02 & 2.77e-01 & 3.10e-03 & 8.78e-02 & 2.92e-03 & 2.42e-01 & 2.92e-03 & 2.03e-01 & 4.24e-03 & -1.81e-02 \\
 & 7.59e-04 & 5.39e-02 & -2.01e-04 & 3.05e-03 & 2.36e-02 & 2.91e-03 & 2.51e-03 & 2.92e-03 & -1.68e-03 & 4.25e-03 & -1.55e-03 \\
 \midrule
\multirow{6}{*}{$\fracNet_2$} & 2.25e-02 & 1.14e-01 &  & 3.25e-02 &  & 6.52e-02 &  & 6.77e-02 &  & 5.22e-02 &  \\
 & 1.12e-02 & 6.41e-02 & 8.25e-01 & 7.58e-03 & 2.10e+00 & 2.47e-02 & 1.40e+00 & 2.70e-02 & 1.33e+00 & 1.54e-02 & 1.76e+00 \\
 & 5.62e-03 & 9.66e-01 & -3.91e+00 & 4.15e-03 & 8.70e-01 & 9.36e-03 & 1.40e+00 & 1.06e-02 & 1.34e+00 & 4.81e-03 & 1.68e+00 \\
 & 2.81e-03 & 5.66e-02 & 4.09e+00 & 4.14e-03 & 3.25e-03 & 4.57e-03 & 1.03e+00 & 5.17e-03 & 1.04e+00 & 4.96e-03 & -4.21e-02 \\
 & 1.41e-03 & 1.02e-02 & 2.48e+00 & 3.80e-03 & 1.23e-01 & 3.46e-03 & 4.03e-01 & 3.67e-03 & 4.96e-01 & 5.24e-03 & -7.93e-02 \\
 & 7.03e-04 & 2.21e-01 & -4.44e+00 & 3.57e-03 & 8.85e-02 & 3.31e-03 & 6.21e-02 & 3.38e-03 & 1.19e-01 & 5.15e-03 & 2.43e-02 \\
 \midrule
\multirow{6}{*}{$\fracNet_3$} & 1.70e-02 & 7.42e-02 &  & 2.20e-02 &  & 1.18e-01 &  & 1.16e-01 &  & 8.49e-02 &  \\
 & 8.52e-03 & 6.27e-02 & 2.45e-01 & 1.94e-02 & 1.79e-01 & 5.99e-02 & 9.81e-01 & 5.66e-02 & 1.03e+00 & 3.08e-02 & 1.46e+00 \\
 & 4.26e-03 & 5.29e-02 & 2.44e-01 & 1.85e-02 & 7.03e-02 & 3.21e-02 & 9.00e-01 & 3.02e-02 & 9.06e-01 & 2.29e-02 & 4.25e-01 \\
 & 2.13e-03 & 8.13e-02 & -6.20e-01 & 1.67e-02 & 1.45e-01 & 1.79e-02 & 8.42e-01 & 1.70e-02 & 8.30e-01 & 2.39e-02 & -6.19e-02 \\
 & 1.06e-03 & 1.19e-01 & -5.54e-01 & 1.51e-02 & 1.50e-01 & 1.22e-02 & 5.61e-01 & 1.21e-02 & 4.87e-01 & 2.46e-02 & -3.91e-02 \\
 & 5.32e-04 & 5.90e-02 & 1.02e+00 & 1.41e-02 & 9.17e-02 & 1.20e-02 & 1.46e-02 & 1.22e-02 & -5.12e-03 & 2.48e-02 & -1.16e-02 \\ \midrule
 \multicolumn{12}{c}{$k = \num{1e-4}$} \\ \midrule
 \multirow{6}{*}{$\fracNet_1$} & 2.43e-02 & 1.93e-01 &  & 5.04e-02 &  & 5.97e-02 &  & 5.96e-02 &  & 5.48e-02 &  \\
 & 1.21e-02 & 1.54e-01 & 3.28e-01 & 2.15e-02 & 1.23e+00 & 2.02e-02 & 1.56e+00 & 2.01e-02 & 1.57e+00 & 1.88e-02 & 1.54e+00 \\
 & 6.07e-03 & 1.33e-01 & 2.09e-01 & 8.55e-03 & 1.33e+00 & 7.58e-03 & 1.41e+00 & 7.56e-03 & 1.41e+00 & 1.58e-02 & 2.54e-01 \\
 & 3.04e-03 & 1.19e-01 & 1.53e-01 & 3.69e-03 & 1.21e+00 & 3.62e-03 & 1.07e+00 & 3.61e-03 & 1.07e+00 & 1.78e-02 & -1.73e-01 \\
 & 1.52e-03 & 1.07e-01 & 1.53e-01 & 1.97e-03 & 9.06e-01 & 2.33e-03 & 6.35e-01 & 2.33e-03 & 6.32e-01 & 1.89e-02 & -9.08e-02 \\
 & 7.59e-04 & 1.01e-01 & 8.78e-02 & 1.74e-03 & 1.77e-01 & 2.04e-03 & 1.87e-01 & 2.05e-03 & 1.87e-01 & 1.95e-02 & -4.28e-02 \\
 \midrule
\multirow{6}{*}{$\fracNet_2$} & 2.25e-02 & 2.40e-01 &  & 1.84e-01 &  & 9.49e-02 &  & 8.44e-02 &  & 2.43e-01 &  \\
 & 1.12e-02 & 2.04e-01 & 2.32e-01 & 9.34e-02 & 9.80e-01 & 5.02e-02 & 9.17e-01 & 4.40e-02 & 9.41e-01 & 2.34e-01 & 5.34e-02 \\
 & 5.62e-03 & 1.97e-01 & 4.93e-02 & 4.87e-02 & 9.39e-01 & 2.82e-02 & 8.32e-01 & 2.60e-02 & 7.60e-01 & 2.31e-01 & 2.02e-02 \\
 & 2.81e-03 & 2.03e-01 & -4.01e-02 & 2.92e-02 & 7.35e-01 & 1.81e-02 & 6.42e-01 & 1.77e-02 & 5.52e-01 & 2.30e-01 & 7.35e-03 \\
 & 1.41e-03 & 2.58e-01 & -3.47e-01 & 1.80e-02 & 7.01e-01 & 1.12e-02 & 6.87e-01 & 1.11e-02 & 6.80e-01 & 2.29e-01 & 3.00e-03 \\
 & 7.03e-04 & 1.88e-01 & 4.59e-01 & 1.06e-02 & 7.56e-01 & 6.97e-03 & 6.89e-01 & 6.81e-03 & 7.00e-01 & 2.29e-01 & 1.33e-03 \\
 \midrule
\multirow{6}{*}{$\fracNet_3$} & 1.70e-02 & 1.31e-01 &  & 9.83e-02 &  & 8.89e-02 &  & 8.25e-02 &  & 2.25e-01 &  \\
 & 8.52e-03 & 9.04e-02 & 5.31e-01 & 5.33e-02 & 8.84e-01 & 4.22e-02 & 1.08e+00 & 3.97e-02 & 1.06e+00 & 2.12e-01 & 8.45e-02 \\
 & 4.26e-03 & 1.24e-01 & -4.59e-01 & 3.43e-02 & 6.35e-01 & 2.91e-02 & 5.34e-01 & 2.91e-02 & 4.48e-01 & 2.09e-01 & 2.32e-02 \\
 & 2.13e-03 & 1.28e-01 & -4.56e-02 & 2.65e-02 & 3.71e-01 & 2.33e-02 & 3.24e-01 & 2.36e-02 & 3.03e-01 & 2.08e-01 & 7.11e-03 \\
 & 1.06e-03 & 1.48e-01 & -2.05e-01 & 2.05e-02 & 3.71e-01 & 1.86e-02 & 3.25e-01 & 1.86e-02 & 3.45e-01 & 2.07e-01 & 2.73e-03 \\
 & 5.32e-04 & 5.12e-01 & -1.79e+00 & 1.71e-02 & 2.65e-01 & 1.65e-02 & 1.68e-01 & 1.65e-02 & 1.70e-01 & 2.07e-01 & 1.20e-03 \\
  \bottomrule
  \end{tabular}
  \label{tab:convDiscreteErrors_head1_angle0}
  }
\end{table}

\begin{table}[h]
  {\footnotesize
  \centering
  \caption{\textbf{Case 2 - errors and rates of $\head_2$ for $\permAngle = \pi/4$}. Listed are the errors $\errorNorm_{\head_2}$ and the corresponding rates $r_{\head_2}$ over grid refinement, expressed in terms of $\discLength_\fracIdx$.}
  \begin{tabular}{ *{2}{l} | *{2}{l} | *{2}{l} | *{2}{l} | *{2}{l} | *{2}{l} }
  \toprule
  & & \multicolumn{2}{l}{\tpfaDfm} & \multicolumn{2}{l}{\mpfaDfm} & \multicolumn{2}{l}{\eboxDfm} & \multicolumn{2}{l}{\eboxMortarDfm} & \multicolumn{2}{l}{\boxDfm} \\
  &  $\discLength_\fracIdx$ & $\errorNorm_{\head_2}$ & $r_{\head_2}$ & $\errorNorm_{\head_2}$ & $r_{\head_2}$ & $\errorNorm_{\head_2}$ & $r_{\head_2}$ & $\errorNorm_{\head_2}$ & $r_{\head_2}$ & $\errorNorm_{\head_2}$ & $r_{\head_2}$ \\ \midrule
  \multicolumn{12}{c}{$k = \num{1e4}$} \\ \midrule
  \multirow{6}{*}{$\fracNet_1$} & 2.43e-02 & 1.82e-01 &  & 1.53e-01 &  & 1.65e-01 &  & 1.62e-01 &  & 1.65e-01 &  \\
 & 1.21e-02 & 1.01e-01 & 8.53e-01 & 7.25e-02 & 1.08e+00 & 7.86e-02 & 1.07e+00 & 7.63e-02 & 1.08e+00 & 7.89e-02 & 1.07e+00 \\
 & 6.07e-03 & 6.88e-02 & 5.53e-01 & 3.47e-02 & 1.06e+00 & 3.84e-02 & 1.04e+00 & 3.69e-02 & 1.05e+00 & 3.86e-02 & 1.03e+00 \\
 & 3.04e-03 & 5.83e-02 & 2.37e-01 & 1.61e-02 & 1.11e+00 & 1.82e-02 & 1.08e+00 & 1.73e-02 & 1.09e+00 & 1.84e-02 & 1.07e+00 \\
 & 1.52e-03 & 5.60e-02 & 5.90e-02 & 6.92e-03 & 1.22e+00 & 8.05e-03 & 1.18e+00 & 7.51e-03 & 1.20e+00 & 8.22e-03 & 1.16e+00 \\
 & 7.59e-04 & 8.88e-02 & -6.64e-01 & 3.30e-03 & 1.07e+00 & 3.70e-03 & 1.12e+00 & 3.42e-03 & 1.14e+00 & 3.86e-03 & 1.09e+00 \\
 \midrule
\multirow{6}{*}{$\fracNet_2$} & 2.25e-02 & 1.86e-01 &  & 1.65e-01 &  & 1.99e-01 &  & 1.89e-01 &  & 2.00e-01 &  \\
 & 1.12e-02 & 1.18e-01 & 6.59e-01 & 8.30e-02 & 9.92e-01 & 9.79e-02 & 1.03e+00 & 9.38e-02 & 1.01e+00 & 9.85e-02 & 1.02e+00 \\
 & 5.62e-03 & 8.88e-02 & 4.10e-01 & 4.17e-02 & 9.93e-01 & 4.81e-02 & 1.02e+00 & 4.61e-02 & 1.02e+00 & 4.86e-02 & 1.02e+00 \\
 & 2.81e-03 & 9.28e-02 & -6.37e-02 & 2.02e-02 & 1.04e+00 & 2.32e-02 & 1.05e+00 & 2.22e-02 & 1.06e+00 & 2.35e-02 & 1.04e+00 \\
 & 1.41e-03 & 8.63e-02 & 1.05e-01 & 9.37e-03 & 1.11e+00 & 1.08e-02 & 1.11e+00 & 1.02e-02 & 1.11e+00 & 1.11e-02 & 1.08e+00 \\
 & 7.03e-04 & 7.36e-02 & 2.31e-01 & 4.33e-03 & 1.11e+00 & 4.95e-03 & 1.12e+00 & 4.68e-03 & 1.13e+00 & 5.33e-03 & 1.06e+00 \\
 \midrule
\multirow{6}{*}{$\fracNet_3$} & 1.70e-02 & 2.53e-01 &  & 1.33e-01 &  & 1.78e-01 &  & 1.63e-01 &  & 1.81e-01 &  \\
 & 8.52e-03 & 2.11e-01 & 2.64e-01 & 6.92e-02 & 9.46e-01 & 8.66e-02 & 1.04e+00 & 8.00e-02 & 1.02e+00 & 8.88e-02 & 1.03e+00 \\
 & 4.26e-03 & 4.47e-01 & -1.08e+00 & 3.38e-02 & 1.03e+00 & 4.30e-02 & 1.01e+00 & 3.97e-02 & 1.01e+00 & 4.54e-02 & 9.69e-01 \\
 & 2.13e-03 & 1.89e-01 & 1.24e+00 & 1.61e-02 & 1.07e+00 & 2.17e-02 & 9.89e-01 & 2.00e-02 & 9.92e-01 & 2.46e-02 & 8.85e-01 \\
 & 1.06e-03 & 2.03e-01 & -1.05e-01 & 8.27e-03 & 9.64e-01 & 1.17e-02 & 8.93e-01 & 1.08e-02 & 8.82e-01 & 1.53e-02 & 6.84e-01 \\
 & 5.32e-04 & 2.81e-01 & -4.70e-01 & 5.86e-03 & 4.98e-01 & 7.65e-03 & 6.09e-01 & 7.31e-03 & 5.68e-01 & 1.18e-02 & 3.78e-01 \\ \midrule
  \multicolumn{12}{c}{$k = \num{1e-4}$} \\ \midrule
\multirow{6}{*}{$\fracNet_1$} & 2.43e-02 & 1.68e-01 &  & 1.43e-01 &  & 1.49e-01 &  & 1.49e-01 &  & 1.49e-01 &  \\
 & 1.21e-02 & 9.64e-02 & 8.05e-01 & 6.77e-02 & 1.08e+00 & 6.83e-02 & 1.12e+00 & 6.83e-02 & 1.12e+00 & 6.94e-02 & 1.10e+00 \\
 & 6.07e-03 & 6.87e-02 & 4.89e-01 & 3.22e-02 & 1.07e+00 & 3.21e-02 & 1.09e+00 & 3.21e-02 & 1.09e+00 & 3.41e-02 & 1.02e+00 \\
 & 3.04e-03 & 5.94e-02 & 2.09e-01 & 1.47e-02 & 1.13e+00 & 1.45e-02 & 1.14e+00 & 1.45e-02 & 1.14e+00 & 1.84e-02 & 8.92e-01 \\
 & 1.52e-03 & 5.68e-02 & 6.55e-02 & 6.09e-03 & 1.27e+00 & 5.99e-03 & 1.28e+00 & 5.99e-03 & 1.28e+00 & 1.27e-02 & 5.36e-01 \\
 & 7.59e-04 & 6.65e-02 & -2.28e-01 & 2.77e-03 & 1.14e+00 & 2.68e-03 & 1.16e+00 & 2.68e-03 & 1.16e+00 & 1.15e-02 & 1.44e-01 \\
 \midrule
\multirow{6}{*}{$\fracNet_2$} & 2.25e-02 & 1.62e-01 &  & 1.46e-01 &  & 1.49e-01 &  & 1.49e-01 &  & 2.15e-01 &  \\
 & 1.12e-02 & 1.16e-01 & 4.85e-01 & 7.81e-02 & 9.03e-01 & 7.04e-02 & 1.08e+00 & 7.04e-02 & 1.08e+00 & 1.72e-01 & 3.21e-01 \\
 & 5.62e-03 & 1.04e-01 & 1.63e-01 & 4.11e-02 & 9.25e-01 & 3.40e-02 & 1.05e+00 & 3.40e-02 & 1.05e+00 & 1.62e-01 & 9.46e-02 \\
 & 2.81e-03 & 1.19e-01 & -1.97e-01 & 2.13e-02 & 9.47e-01 & 1.59e-02 & 1.09e+00 & 1.59e-02 & 1.09e+00 & 1.59e-01 & 2.46e-02 \\
 & 1.41e-03 & 1.00e-01 & 2.49e-01 & 1.16e-02 & 8.80e-01 & 7.44e-03 & 1.10e+00 & 7.44e-03 & 1.10e+00 & 1.58e-01 & 6.31e-03 \\
 & 7.03e-04 & 1.05e-01 & -7.10e-02 & 7.36e-03 & 6.55e-01 & 4.54e-03 & 7.12e-01 & 4.54e-03 & 7.12e-01 & 1.58e-01 & 1.68e-03 \\
 \midrule
\multirow{6}{*}{$\fracNet_3$} & 1.70e-02 & 1.36e-01 &  & 1.30e-01 &  & 1.13e-01 &  & 1.13e-01 &  & 2.12e-01 &  \\
 & 8.52e-03 & 1.05e-01 & 3.77e-01 & 5.86e-02 & 1.15e+00 & 4.88e-02 & 1.21e+00 & 4.88e-02 & 1.21e+00 & 1.78e-01 & 2.56e-01 \\
 & 4.26e-03 & 1.64e-01 & -6.41e-01 & 2.81e-02 & 1.06e+00 & 2.34e-02 & 1.06e+00 & 2.35e-02 & 1.06e+00 & 1.70e-01 & 6.22e-02 \\
 & 2.13e-03 & 9.47e-02 & 7.91e-01 & 1.37e-02 & 1.04e+00 & 1.20e-02 & 9.61e-01 & 1.20e-02 & 9.61e-01 & 1.68e-01 & 1.51e-02 \\
 & 1.06e-03 & 1.11e-01 & -2.26e-01 & 7.51e-03 & 8.65e-01 & 7.51e-03 & 6.82e-01 & 7.51e-03 & 6.82e-01 & 1.68e-01 & 3.74e-03 \\
 & 5.32e-04 & 4.34e-01 & -1.97e+00 & 5.81e-03 & 3.72e-01 & 6.27e-03 & 2.59e-01 & 6.27e-03 & 2.60e-01 & 1.68e-01 & 9.02e-04 \\
  \bottomrule
  \end{tabular}
  \label{tab:convDiscreteErrors_head2_angle45}
  }
\end{table}

\begin{table}[h]
  {\footnotesize
  \centering
  \caption{\textbf{Case 2 - errors and rates of $\head_1$ for $\permAngle = \pi/4$}. Listed are the errors $\errorNorm_{\head_1}$ and the corresponding rates $r_{\head_1}$ over grid refinement, expressed in terms of $\discLength_\fracIdx$.}
  \begin{tabular}{ *{2}{l} | *{2}{l} | *{2}{l} | *{2}{l} | *{2}{l} | *{2}{l} }
  \toprule
  & & \multicolumn{2}{l}{\tpfaDfm} & \multicolumn{2}{l}{\mpfaDfm} & \multicolumn{2}{l}{\eboxDfm} & \multicolumn{2}{l}{\eboxMortarDfm} & \multicolumn{2}{l}{\boxDfm} \\
  &  $\discLength_\fracIdx$ & $\errorNorm_{\head_1}$ & $r_{\head_1}$ & $\errorNorm_{\head_1}$ & $r_{\head_1}$ & $\errorNorm_{\head_1}$ & $r_{\head_1}$ & $\errorNorm_{\head_1}$ & $r_{\head_1}$ & $\errorNorm_{\head_1}$ & $r_{\head_1}$ \\ \midrule
  \multicolumn{12}{c}{$k = \num{1e4}$} \\ \midrule
  \multirow{6}{*}{$\fracNet_1$} & 2.43e-02 & 2.83e-02 &  & 3.76e-02 &  & 8.00e-02 &  & 7.93e-02 &  & 7.58e-02 &  \\
 & 1.21e-02 & 4.15e-02 & -5.52e-01 & 6.69e-03 & 2.49e+00 & 1.40e-02 & 2.51e+00 & 1.43e-02 & 2.47e+00 & 1.44e-02 & 2.39e+00 \\
 & 6.07e-03 & 4.62e-02 & -1.53e-01 & 3.15e-03 & 1.08e+00 & 2.22e-03 & 2.66e+00 & 2.52e-03 & 2.51e+00 & 4.26e-03 & 1.76e+00 \\
 & 3.04e-03 & 4.91e-02 & -8.77e-02 & 2.88e-03 & 1.30e-01 & 2.12e-03 & 7.16e-02 & 2.22e-03 & 1.84e-01 & 3.76e-03 & 1.80e-01 \\
 & 1.52e-03 & 5.07e-02 & -4.78e-02 & 2.79e-03 & 4.49e-02 & 2.40e-03 & -1.81e-01 & 2.45e-03 & -1.41e-01 & 3.78e-03 & -7.47e-03 \\
 & 7.59e-04 & 5.92e-02 & -2.23e-01 & 2.74e-03 & 2.66e-02 & 2.54e-03 & -8.17e-02 & 2.56e-03 & -6.71e-02 & 3.77e-03 & 2.32e-03 \\
 \midrule
\multirow{6}{*}{$\fracNet_2$} & 2.25e-02 & 1.47e-02 &  & 1.41e-02 &  & 1.14e-02 &  & 2.95e-02 &  & 3.06e-03 &  \\
 & 1.12e-02 & 3.27e-02 & -1.15e+00 & 5.70e-03 & 1.31e+00 & 1.53e-03 & 2.89e+00 & 7.88e-03 & 1.90e+00 & 4.70e-03 & -6.18e-01 \\
 & 5.62e-03 & 3.37e-02 & -4.14e-02 & 3.34e-03 & 7.71e-01 & 2.61e-03 & -7.69e-01 & 2.75e-03 & 1.52e+00 & 4.30e-03 & 1.27e-01 \\
 & 2.81e-03 & 3.00e-02 & 1.67e-01 & 3.15e-03 & 8.52e-02 & 2.87e-03 & -1.36e-01 & 2.49e-03 & 1.44e-01 & 4.03e-03 & 9.20e-02 \\
 & 1.41e-03 & 7.12e-02 & -1.25e+00 & 3.03e-03 & 5.65e-02 & 2.82e-03 & 2.36e-02 & 2.67e-03 & -1.00e-01 & 4.02e-03 & 6.06e-03 \\
 & 7.03e-04 & 3.99e-02 & 8.37e-01 & 2.94e-03 & 4.31e-02 & 2.81e-03 & 6.18e-03 & 2.76e-03 & -4.93e-02 & 4.00e-03 & 7.81e-03 \\
 \midrule
\multirow{6}{*}{$\fracNet_3$} & 1.70e-02 & 2.34e-01 &  & 2.10e-02 &  & 2.56e-02 &  & 3.56e-02 &  & 2.02e-02 &  \\
 & 8.52e-03 & 2.12e-01 & 1.46e-01 & 1.50e-02 & 4.81e-01 & 1.43e-02 & 8.43e-01 & 1.66e-02 & 1.10e+00 & 1.83e-02 & 1.46e-01 \\
 & 4.26e-03 & 3.69e-01 & -8.00e-01 & 1.29e-02 & 2.26e-01 & 1.01e-02 & 4.95e-01 & 1.12e-02 & 5.74e-01 & 1.94e-02 & -8.68e-02 \\
 & 2.13e-03 & 1.71e-01 & 1.11e+00 & 1.19e-02 & 1.11e-01 & 9.50e-03 & 9.08e-02 & 1.00e-02 & 1.55e-01 & 1.99e-02 & -3.51e-02 \\
 & 1.06e-03 & 2.10e-01 & -2.97e-01 & 1.14e-02 & 6.55e-02 & 9.93e-03 & -6.35e-02 & 1.02e-02 & -1.87e-02 & 2.00e-02 & -7.08e-03 \\
 & 5.32e-04 & 1.72e-01 & 2.85e-01 & 1.11e-02 & 3.34e-02 & 1.04e-02 & -6.20e-02 & 1.05e-02 & -4.40e-02 & 1.99e-02 & 5.60e-03 \\ \midrule
  \multicolumn{12}{c}{$k = \num{1e-4}$} \\ \midrule
\multirow{6}{*}{$\fracNet_1$} & 2.43e-02 & 1.04e-01 &  & 4.25e-02 &  & 7.92e-02 &  & 7.92e-02 &  & 8.08e-02 &  \\
 & 1.21e-02 & 9.59e-02 & 1.23e-01 & 9.07e-03 & 2.23e+00 & 1.65e-02 & 2.27e+00 & 1.65e-02 & 2.27e+00 & 1.90e-02 & 2.09e+00 \\
 & 6.07e-03 & 9.51e-02 & 1.26e-02 & 3.13e-03 & 1.54e+00 & 4.06e-03 & 2.02e+00 & 4.06e-03 & 2.02e+00 & 6.33e-03 & 1.58e+00 \\
 & 3.04e-03 & 9.17e-02 & 5.17e-02 & 1.55e-03 & 1.01e+00 & 1.75e-03 & 1.22e+00 & 1.75e-03 & 1.21e+00 & 3.38e-03 & 9.06e-01 \\
 & 1.52e-03 & 8.83e-02 & 5.47e-02 & 1.01e-03 & 6.27e-01 & 1.05e-03 & 7.30e-01 & 1.05e-03 & 7.30e-01 & 2.20e-03 & 6.19e-01 \\
 & 7.59e-04 & 7.09e-02 & 3.17e-01 & 8.15e-04 & 3.03e-01 & 8.22e-04 & 3.58e-01 & 8.22e-04 & 3.59e-01 & 1.63e-03 & 4.29e-01 \\
 \midrule
\multirow{6}{*}{$\fracNet_2$} & 2.25e-02 & 1.23e-01 &  & 8.41e-02 &  & 9.30e-02 &  & 8.50e-02 &  & 2.52e-01 &  \\
 & 1.12e-02 & 1.37e-01 & -1.52e-01 & 5.37e-02 & 6.47e-01 & 4.11e-02 & 1.18e+00 & 3.65e-02 & 1.22e+00 & 2.32e-01 & 1.18e-01 \\
 & 5.62e-03 & 1.50e-01 & -1.33e-01 & 3.78e-02 & 5.06e-01 & 2.52e-02 & 7.09e-01 & 2.41e-02 & 5.96e-01 & 2.26e-01 & 3.57e-02 \\
 & 2.81e-03 & 1.50e-01 & 4.43e-03 & 2.80e-02 & 4.34e-01 & 1.83e-02 & 4.58e-01 & 1.85e-02 & 3.81e-01 & 2.25e-01 & 1.13e-02 \\
 & 1.41e-03 & 1.49e-01 & 9.06e-03 & 2.00e-02 & 4.85e-01 & 1.34e-02 & 4.47e-01 & 1.35e-02 & 4.54e-01 & 2.24e-01 & 4.03e-03 \\
 & 7.03e-04 & 1.55e-01 & -6.24e-02 & 1.41e-02 & 5.05e-01 & 1.04e-02 & 3.66e-01 & 1.04e-02 & 3.78e-01 & 2.24e-01 & 1.60e-03 \\
 \midrule
\multirow{6}{*}{$\fracNet_3$} & 1.70e-02 & 1.05e-01 &  & 1.05e-01 &  & 7.09e-02 &  & 6.55e-02 &  & 2.09e-01 &  \\
 & 8.52e-03 & 1.03e-01 & 4.08e-02 & 4.49e-02 & 1.22e+00 & 3.50e-02 & 1.02e+00 & 3.30e-02 & 9.88e-01 & 1.98e-01 & 8.03e-02 \\
 & 4.26e-03 & 1.29e-01 & -3.27e-01 & 3.00e-02 & 5.81e-01 & 2.50e-02 & 4.88e-01 & 2.51e-02 & 3.94e-01 & 1.96e-01 & 1.55e-02 \\
 & 2.13e-03 & 1.04e-01 & 3.00e-01 & 2.35e-02 & 3.53e-01 & 2.01e-02 & 3.10e-01 & 2.06e-02 & 2.88e-01 & 1.95e-01 & 4.43e-03 \\
 & 1.06e-03 & 1.20e-01 & -2.01e-01 & 1.81e-02 & 3.74e-01 & 1.63e-02 & 3.06e-01 & 1.64e-02 & 3.28e-01 & 1.95e-01 & 1.72e-03 \\
 & 5.32e-04 & 1.73e-01 & -5.27e-01 & 1.51e-02 & 2.61e-01 & 1.47e-02 & 1.44e-01 & 1.48e-02 & 1.52e-01 & 1.95e-01 & 7.49e-04 \\
  \bottomrule
  \end{tabular}
  \label{tab:convDiscreteErrors_head1_angle45}
  }
\end{table}

\end{document}